\PassOptionsToPackage{hypertexnames=false}{hyperref}
\documentclass[pdflatex,sn-mathphys-num]{sn-jnl}

\usepackage{graphicx}%
\usepackage{multirow}%
\usepackage{amsmath,amssymb,amsfonts}%
\usepackage{amsthm}%
\usepackage{mathrsfs}%
\usepackage[title]{appendix}%
\usepackage{xcolor}%
\usepackage{textcomp}%
\usepackage{manyfoot}%
\usepackage{booktabs}%
\usepackage{algorithm}%
\usepackage{algorithmicx}%
\usepackage{algpseudocode}%
\usepackage{listings}%
\usepackage{tikz,enumitem,bm,microtype} 
\usepackage{adjustbox} 
\hypersetup{colorlinks=false,allbordercolors=blue,pdfborderstyle={/S/U/W 1}} 
\usetikzlibrary{arrows.meta,calc,positioning} 
\definecolor{revisionred}{RGB}{180,0,0}

\newcommand{\ii}{\mathrm{i}}

\theoremstyle{thmstyleone}%
\newtheorem{theorem}{Theorem}
\newtheorem{proposition}[theorem]{Proposition}%
\newtheorem{lemma}[theorem]{Lemma}

\theoremstyle{thmstyletwo}%
\newtheorem{remark}{Remark}%

\theoremstyle{thmstylethree}%
\newtheorem{definition}{Definition}%

\begin{document}

\title[Article Title]{Existence and Splitting of Dirac Points in Honeycomb Gyromagnetic Photonic Crystals}


\author[1]{\fnm{Liuyi} \sur{Bao}}\email{baoly24@mails.jlu.edu.cn}

\author*[1]{\fnm{Wei} \sur{Wu}}\email{wei\_wu@jlu.edu.cn}

\affil[1]{\orgdiv{School of Mathematics}, \orgname{Jilin University}, \orgaddress{\street{2699 Qianjin Street}, \city{Changchun}, \postcode{130012}, \state{Jilin}, \country{China}}}

\abstract{We establish the existence of Dirac cones and magnetically induced local band gaps for the transverse magnetic Maxwell system in a two-dimensional honeycomb array of circular yttrium iron garnet (YIG) cylinders. For fixed positive material parameters with nonzero permeability contrast, we prove that sufficiently small cylinders support a Dirac point in the nonmagnetic case $\kappa_g=0$: two Bloch dispersion surfaces meet at a double eigenfrequency and form a nondegenerate cone. For each fixed sufficiently small nonzero gyrotropic bias $\kappa_g$, we show that this degeneracy splits and a local band gap opens near the original Dirac point. Our approach combines a quasi-periodic layer-potential formulation with honeycomb symmetry and rigorously controlled finite-dimensional reductions. The central contribution is a unified derivation of the conical crossing and magnetic gap opening directly from the full Maxwell transmission problem, without assuming an effective discrete model. This establishes the transition from Dirac dispersion to a locally gapped spectrum in a small-inclusion regime with fixed material contrast.}




\maketitle

\section{Introduction}

The study of wave propagation in periodic media brings together spectral theory, asymptotic analysis, and the design of materials with prescribed dispersive properties. The pioneering work of Yablonovitch and John on photonic band-gap structures established the possibility of controlling electromagnetic propagation through the spatial organization of dielectric media~\cite{Yablonovitch1987,John1987}. From the mathematical viewpoint, the relevant information is encoded in the Floquet--Bloch spectrum of a periodic differential operator; see~\cite{Kuchment2016} for a broad account. Of particular interest are Dirac points, at which two dispersion surfaces meet conically. Their existence, the nonvanishing of the cone slope, and their response to symmetry-breaking perturbations are distinct spectral questions that must be resolved at the level of the underlying wave equation.

A rigorous theory of Dirac points in continuous honeycomb structures was developed by Fefferman and Weinstein for periodic Schr\"odinger operators~\cite{FeffermanWeinstein2012,FeffermanWeinstein2014,FeffermanLeeThorpWeinstein2016}. Moreover, magnetic perturbations provide a physically important means of changing this spectral structure. Raghu and Haldane proposed photonic analogues of quantum Hall edge states in time-reversal-breaking periodic media~\cite{RaghuHaldane2008}. Wang and collaborators subsequently studied one-way propagation in gyromagnetic photonic crystals~\cite{Wang2008} and demonstrated unidirectional, backscattering-immune electromagnetic states experimentally~\cite{Wang2009}. More recently, boundary-integral methods have been used to establish interface modes bifurcating from a Dirac point in acoustic waveguides~\cite{QiuLinXieZhang2023} and to study the photonic Hall effect in honeycomb photonic crystals~\cite{LiLinQiuZhang2025, LiLinZhang2023}. These developments motivate a quantitative analysis of the bulk Dirac point and its magnetic splitting in explicit transmission geometries, including the dependence on the size of the inclusions and the material parameters.

In this paper, we study the transverse magnetic Maxwell system in a two-dimensional honeycomb array of circular yttrium iron garnet (YIG) cylinders of radius $\varepsilon$. Using layer-potential techniques, we reformulate the transmission problem as a system of boundary integral equations~\cite{AmmariKosche2024}. The lattice and the positive material parameters are fixed, the permeability contrast is nonzero, and the small parameter is the cylinder radius. This regime differs from the high-contrast subwavelength limit of bubbly crystals: the Dirac frequency considered here approaches the positive background frequency $\omega_{0,m}=|\bm\alpha^*|/\sqrt{\epsilon\mu_{\mathrm m}}$ as $\varepsilon\to0$, where $\bm\alpha^*$ is a fixed Brillouin-zone vertex and $\epsilon$ is the permittivity. We first set the gyrotropic parameter $\kappa_g$ equal to zero and identify an exact double Bloch eigenfrequency. We then determine its conical dispersion and study the splitting induced by the opposite magnetic biases $\kappa_g=\pm\kappa$, with $|\kappa|<\mu$. 

The paper is organized as follows. Section~\ref{sec:2} introduces the geometry, the Maxwell transmission problem, the quasi-periodic layer potentials, and the Green-function decomposition. Section~\ref{sec:3} establishes the rotational sector decomposition and the nonmagnetic inversion--conjugation symmetry. Section~\ref{sec:4} analyzes the candidate frequencies, performs the finite-dimensional reductions, and proves the existence of the exact double frequency and its Dirac cone. Section~\ref{sec:5} treats gyrotropic splitting, the exchange under magnetization reversal, and the common local gap. Appendix~A contains the small-radius operator expansions and the Fourier matrix elements used in the analysis.

\section{Problem formulation and preliminaries}\label{sec:2}

Yttrium iron garnet (YIG) is a synthetic ferrimagnetic garnet used in microwave and optical applications. Its gyrotropic response is modeled by off-diagonal entries in the permeability tensor under a magnetic bias normal to the plane. In this paper, we consider identical, infinitely long circular cylinders parallel to the $x_3$ axis, periodically arranged in a honeycomb pattern in the $(x_1,x_2)$ plane. In the nonmagnetic model where the gyrotropic response is off, time-reversal and lattice-rotation symmetries permit Dirac points near the Brillouin-zone vertices. When the gyrotropic response is present, the additional magnetic bias breaks time-reversal symmetry and can open a local gap near a former Dirac point. We use $\kappa_g\in\{0,+\kappa,-\kappa\}$ to describe the magnetic bias: $0$ describes the nonmagnetic problem, while the nonzero values $+\kappa$ and $-\kappa$ represent opposite magnetizations.

\subsection{Geometry and material model}\label{sec:sec2-1}

Let $a_{\mathrm{lat}}>0$ be the lattice constant. Choose the primitive vectors $\bm e_1=a_{\mathrm{lat}}(\sqrt3/2,1/2)^\top$ and $\bm e_2=a_{\mathrm{lat}}(\sqrt3/2,-1/2)^\top$, and let $\Lambda=\mathbb Z\bm e_1\oplus\mathbb Z\bm e_2$. We choose the primitive cell symmetric about the midpoint of the sublattice sites:
$$
Y=\left\{t_1\bm e_1+t_2\bm e_2: -\frac13\leq t_1,t_2<\frac23\right\}.
$$
The following reciprocal vectors $\bm\kappa_1,\bm\kappa_2$ satisfy $\bm e_i\cdot\bm\kappa_j=2\pi\delta_{ij}$:
$$
\bm\kappa_1=\frac{2\pi}{a_{\mathrm{lat}}} \left(\frac{\sqrt3}{3},1\right)^\top,\qquad  \bm\kappa_2=\frac{2\pi}{a_{\mathrm{lat}}} \left(\frac{\sqrt3}{3},-1\right)^\top.
$$
The reciprocal lattice is $\Lambda^*=\mathbb Z\bm\kappa_1\oplus\mathbb Z\bm\kappa_2$. Denote the first Brillouin zone by $\mathcal B$. We fix the vertex
$$
\bm\alpha^*:=\frac{\bm\kappa_1+2\bm\kappa_2}{3} =\frac{2\pi}{a_{\mathrm{lat}}} \left(\frac1{\sqrt3},-\frac13\right)^\top.
$$
The subsequent local spectral analysis is performed near $\bm\alpha^*$.

\begin{figure}[t]
\centering
\begin{adjustbox}{max width=\textwidth}
\begin{tikzpicture}[scale=1.0, line cap=round, >=Stealth, every node/.style={font=\small}]
\def\a{2.2} \coordinate (O) at (0,0); \coordinate (Eone) at ({sqrt(3)*\a/2},{\a/2}); \coordinate (Etwo) at ({sqrt(3)*\a/2},{-\a/2}); \coordinate (Bbasis) at ({\a/sqrt(3)},0);
\coordinate (Cmid) at ({\a/(2*sqrt(3))},0); \coordinate (Yright) at ({2*\a/sqrt(3)},0); \coordinate (Ytop) at ({\a/(2*sqrt(3))},{\a/2}); \coordinate (Yleft) at ({-\a/sqrt(3)},0); \coordinate (Ybottom) at ({\a/(2*sqrt(3))},{-\a/2}); \foreach \i in {-1,0,1}{ \foreach \j in {-1,0,1}{ \coordinate (Apt) at ({sqrt(3)*\a*(\i+\j)/2},{\a*(\i-\j)/2}); \coordinate (Bpt) at ($(Apt)+({\a/sqrt(3)},0)$); \draw[gray!35] (Apt) circle (0.12); \fill[blue!65] (Apt) circle (1.7pt); \draw[gray!35] (Bpt) circle (0.12); \fill[red!70] (Bpt) circle (1.7pt); } } \draw[very thick,revisionred] (Yright) -- (Ytop) -- (Yleft) -- (Ybottom) -- cycle; \node[revisionred] at ({0.28*\a},{0.36*\a}) {$Y$}; \draw[densely dashed,revisionred] (O) -- (Bbasis); \fill[revisionred] (Cmid) circle (1.6pt); \node[below=3pt,revisionred] at (Cmid) {$\bm c$}; \draw[-Stealth, very thick, blue!70!black] (O) -- (Eone) node[pos=0.82, above left] {$\bm e_1$}; \draw[-Stealth, very thick, blue!70!black] (O) -- (Etwo) node[pos=0.82, below left] {$\bm e_2$}; \node[below=3pt] at (O) {$A$}; \node[below=3pt] at (Bbasis) {$B$};
\end{tikzpicture}
\end{adjustbox}
\caption{Geometry of the honeycomb YIG array. The primitive translations are $\bm e_1,\bm e_2$, and $A,B$ are the two basis sites. The red parallelogram centered at $\bm c$ represents the fundamental cell $Y$, which contains both complete disks symmetrically.} \label{fig:lattice-geometry}
\end{figure}
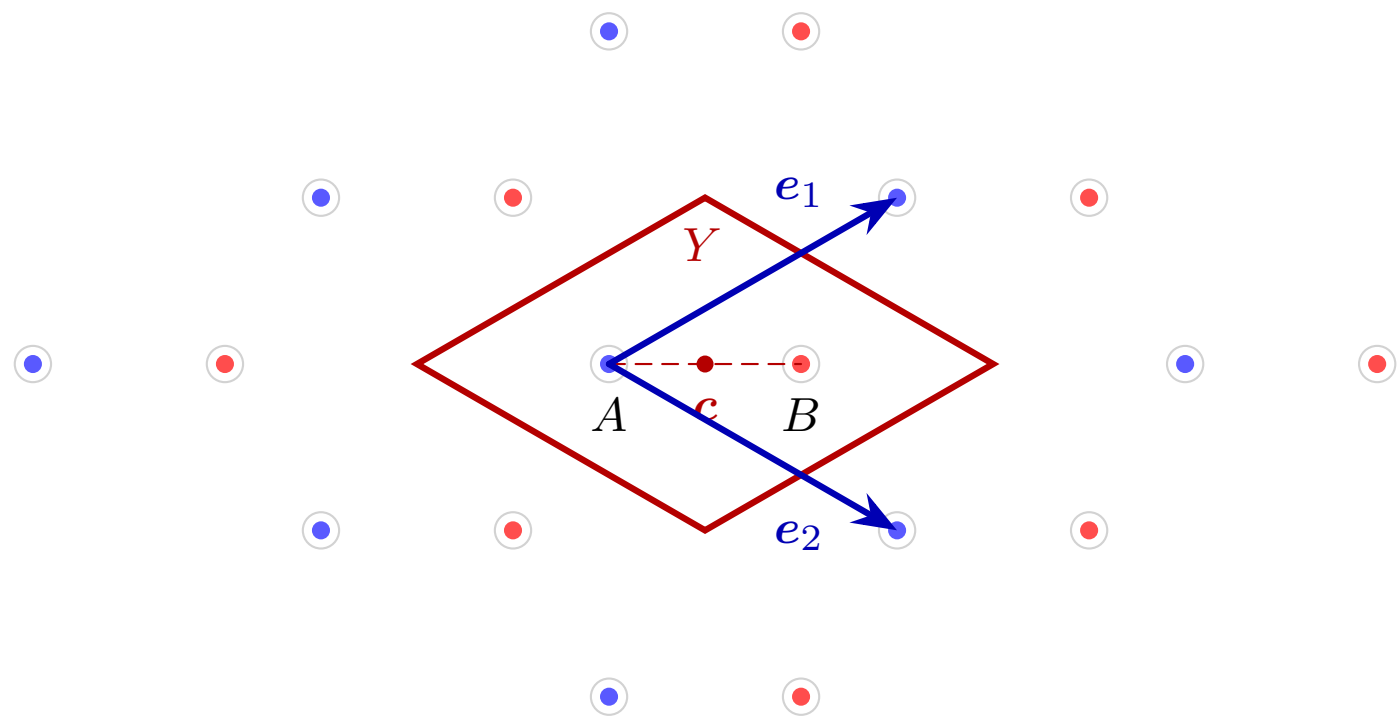

As illustrated in Figure~\ref{fig:lattice-geometry}, the honeycomb array is obtained by placing two distinct basis sites in each cell of the same lattice. We label these sites $A$ and $B$. Choose $\bm z_A=\bm0$ and $\bm z_B=(\bm e_1+\bm e_2)/3$. Their midpoint, also the center of $Y$, is
$$
\bm c:=\frac{\bm z_A+\bm z_B}{2} =\frac{\bm e_1+\bm e_2}{6},\qquad  \bm z_A-\bm c=-(\bm z_B-\bm c).
$$
The two sublattices are $\mathcal L_A=\{\bm z_A+\bm\ell:\bm\ell\in\Lambda\}$ and $\mathcal L_B=\{\bm z_B+\bm\ell:\bm\ell\in\Lambda\}$. Let $D_1=\{\bm x\in\mathbb R^2:|\bm x|<1\}$ be the unit disk. The cylinder of type $p$ centered at $\bm z_p+\bm\ell$ has cross-section $D_{p,\bm\ell}^{\varepsilon} =\bm z_p+\bm\ell+\varepsilon D_1$, where $p=A,B$. The full periodic inclusion set is
$$
D_{\mathrm{per}}^\varepsilon = \bigcup_{\bm\ell\in\Lambda} \left( D_{A,\bm\ell}^{\varepsilon} \cup D_{B,\bm\ell}^{\varepsilon} \right).
$$

For one primitive cell, write $D^\varepsilon=D_A^\varepsilon\cup D_B^\varepsilon$, where $D_A^\varepsilon=\bm z_A+\varepsilon D_1$, $D_B^\varepsilon=\bm z_B+\varepsilon D_1$, with total boundary $\partial D^\varepsilon=\partial D_A^\varepsilon\cup\partial D_B^\varepsilon$. The perpendicular distance between opposite sides of $Y$ is $\sqrt3a_{\mathrm{lat}}/2$. Each site is one third of this distance from its nearest sides, so
$$
\operatorname{dist}(\bm z_A,\partial Y)=\operatorname{dist}(\bm z_B,\partial Y)=\frac{\sqrt3}{6}a_{\mathrm{lat}}.
$$
To ensure $\partial D^\varepsilon\cap\partial Y=\emptyset$ and $\partial D^\varepsilon_A\cap\partial D^\varepsilon_B=\emptyset$, we assume the geometric bound $0<\varepsilon<\sqrt{3}a_{\mathrm{lat}}/6$. For $\bm x\in\partial D_p^\varepsilon$, use the parametrization $\bm x=\bm z_p+\varepsilon\bm r(\theta)$, where $\bm r(\theta)=(\cos\theta,\sin\theta)^\top$, $0\leq\theta<2\pi$. The outward unit normal is $\bm\nu(\theta)=(\cos\theta,\sin\theta)^\top$, and the clockwise unit tangent is $\bm\tau(\theta)=(\sin\theta,-\cos\theta)^\top$. Thus $\partial_{\bm\tau}u=\bm\tau\cdot\nabla u =\nu_2\partial_{x_1}u-\nu_1\partial_{x_2}u$.

The background material is homogeneous and isotropic, with permeability $\mu_{\mathrm m}$ and permittivity $\epsilon$. Inside the cylinders, the permeability tensor and its inverse are
$$
\bm\mu(\kappa_g)=
\begin{pmatrix}
\mu&\ii\kappa_g&0\\
-\ii\kappa_g&\mu&0\\
0&0&\mu_0
\end{pmatrix}
,  \qquad \bm\mu(\kappa_g)^{-1}= \frac1{\mu^2-\kappa_g^2}
\begin{pmatrix}
\mu&-\ii\kappa_g&0\\
\ii\kappa_g&\mu&0\\
0&0&\mu_0^{-1}(\mu^2-\kappa_g^2)
\end{pmatrix},
$$
where the parameter $\kappa_g$ is constant, with value $0$, $+\kappa$, or $-\kappa$. Write $\widetilde\mu^{-1}=\mu/(\mu^2-\kappa_g^2)$ and $\eta=-\kappa_g/(\mu^2-\kappa_g^2)$. We require permittivity $\epsilon>0$, $\mu_{\mathrm m}>0$, $\mu>0$, $\mu_0>0$, $\mu_{\mathrm m}\neq\mu$ and $\mu>|\kappa_g|$. The material coefficients $\epsilon$, $\mu_{\mathrm m}$, $\mu$, and $\kappa_g$ are regarded as frequency-independent parameters. It is immediate that $\bm\mu(\kappa_g)$ is Hermitian positive definite.

We consider time-harmonic TM-polarized electromagnetic waves in this paper. The electric-field amplitude is $\bm E=(0,0,u)^\top$, and the three-dimensional Maxwell system reduces to the following scalar Helmholtz equations in the $(x_1,x_2)$ plane:
\begin{align}\label{eq:system}
\left\{
\begin{aligned}
&\Delta u+\omega_{\mathrm m}(\omega)^2u=0 && \text{in } \mathbb{R}^2\setminus\overline{D_{\mathrm{per}}^\varepsilon}, \\
&\Delta u+\omega_{\mathrm c}(\omega)^2u=0 && \text{in } D_{\mathrm{per}}^\varepsilon, \\
&u|_+=u|_- && \text{on } \partial D_{\mathrm{per}}^\varepsilon, \\
&\frac{\mu}{\mu^2-\kappa_g^2}\partial_{\bm\nu} u|_- +\frac{\ii\kappa_g}{\mu^2-\kappa_g^2}\partial_{\bm\tau} u|_- =\mu_{\mathrm m}^{-1}\partial_{\bm\nu} u|_+ && \text{on } \partial D_{\mathrm{per}}^\varepsilon, \\
& u(\bm x+\bm\ell) = e^{\ii\bm\alpha\cdot\bm\ell}u(\bm x),\quad  \bm\ell\in\Lambda,
\end{aligned}
\right.
\end{align}
where $\omega_{\mathrm{m}}(\omega)^2=\omega^2\epsilon\mu_{\mathrm{m}}, \omega_{\mathrm{c}}(\omega)^2=\omega^2\epsilon\widetilde{\mu}$.

The fourth equation of \eqref{eq:system} can be derived as follows. With time dependence $e^{-\ii\omega t}$, $\bm E$ and $\bm H$ satisfy
$$
\nabla\times\bm E=\ii\omega\bm\mu\bm H, \qquad  \nabla\times\bm H=-\ii\omega\epsilon\bm E.
$$
For $\bm E=(0,0,u)^\top$,
$$
\bm H =\frac1{\ii\omega}\bm\mu^{-1}
\begin{pmatrix}
\partial_yu\\
-\partial_xu\\
0
\end{pmatrix}
=\frac1{\ii\omega}
\begin{pmatrix}
{\widetilde\mu}^{-1} \partial_yu-\ii\eta \partial_xu\\
-\ii\eta \partial_yu-{\widetilde\mu}^{-1} \partial_xu\\
0
\end{pmatrix}.
$$
The third component of the second Maxwell equation gives $ -{\widetilde\mu}^{-1} \Delta u=\omega^2\epsilon u, $ and therefore
$$
\Delta u+\omega^2\epsilon\widetilde\mu u=0,  \qquad \omega_{\mathrm c}(\omega)^2=\omega^2\epsilon\widetilde\mu.
$$
Continuity of the tangential magnetic field is equivalent to continuity of $\bm\nu\times\bm H$. Since
$$
(\bm\nu\times\bm H)_3 = \frac1{\ii\omega} \left( -{\widetilde\mu}^{-1} \partial_{\bm\nu}u +\ii\eta \partial_{\bm\tau}u \right)
$$
and, in the isotropic background, $(\bm\nu\times\bm H)_3 =-(\ii\omega\mu_{\mathrm m})^{-1}\partial_{\bm\nu}u$, we obtain
$$
\frac{\mu}{\mu^2-\kappa_g^2}\partial_{\bm\nu}u|_- +\frac{\ii\kappa_g}{\mu^2-\kappa_g^2}\partial_{\bm\tau}u|_- =\mu_{\mathrm m}^{-1}\partial_{\bm\nu}u|_+,
$$
which is the fourth equation in~\eqref{eq:system}.

\subsection{Layer-potential operators}

We use quasi-periodic Green functions to represent the fields and convert~\eqref{eq:system} into an integral equation on $\partial D^\varepsilon$. For $\lambda\in\{\omega_{\mathrm m},\omega_{\mathrm c}\}$, let $G_\lambda^{\bm\alpha}(\bm x,\bm y)=G_\lambda^{\bm\alpha}(\bm x-\bm y)$ satisfy
$$
(\Delta+\lambda^2)G_\lambda^{\bm\alpha}(\bm x,\bm y)=\sum_{\bm R\in\Lambda}e^{\ii\bm\alpha\cdot\bm R}\delta(\bm x-\bm y-\bm R).
$$
Poisson summation and the Bloch--Fourier expansion give~\cite{Moroz2006}
$$
G_\lambda^{\bm\alpha}(\bm x,\bm y)=\frac{1}{|Y|}\sum_{\bm q\in\Lambda^*}\frac{e^{\ii(\bm q+\bm\alpha)\cdot(\bm x-\bm y)}}{\lambda^2-(\bm q+\bm\alpha)\cdot(\bm q+\bm\alpha)}.
$$

\begin{definition}
For $\lambda\in\{\omega_{\mathrm m},\omega_{\mathrm c}\}$, define the single-layer operator on the total boundary by
$$
\mathcal{S}_{\lambda}^{\bm\alpha,\varepsilon}[\rho](\bm x):=\int_{\partial D^\varepsilon}G_\lambda^{\bm\alpha}(\bm x,\bm y)\rho(\bm y)d\sigma(\bm y),\qquad  \bm x\in\partial D^\varepsilon.
$$
Define the Neumann--Poincar\'e operator by
$$
(\mathcal{K}_{\lambda}^{-\bm\alpha,\varepsilon})^*[\rho](\bm x):=\operatorname{p.v.}\int_{\partial D^\varepsilon} \frac{\partial G_\lambda^{\bm\alpha}(\bm x,\bm y)}{\partial \bm\nu_{\bm x}}\rho(\bm y)d\sigma(\bm y),\qquad  \bm x\in\partial D^\varepsilon.
$$
Define the tangential operator by \cite{AmmariKangLee2009,ColtonKress2013}
$$
\mathcal{T}_{\lambda}^{\bm\alpha,\varepsilon}[\rho](\bm x):=\operatorname{p.v.}\int_{\partial D^\varepsilon} \frac{\partial G_\lambda^{\bm\alpha}(\bm x,\bm y)}{\partial \bm\tau_{\bm x}}\rho(\bm y)d\sigma(\bm y),\qquad  \bm x\in\partial D^\varepsilon.
$$
\end{definition}

The boundary traces of the single-layer potential satisfy the following jump relations.

\begin{lemma}\label{lem:lem1}
For $\rho\in H^{-1/2}(\partial D^\varepsilon)$, the following jump relations hold; see~\cite{AmmariFitzpatrickKangRuizYuZhang2018}.
\begin{enumerate}
\item The single-layer traces are continuous across the boundary:
$$
\mathcal{S}_{\lambda}^{\bm\alpha,\varepsilon}[\rho]\big|_+=\mathcal{S}_{\lambda}^{\bm\alpha,\varepsilon}[\rho]\big|_-.
$$
\item The normal derivatives satisfy
$$
\frac{\partial}{\partial\bm\nu}\mathcal{S}_{\lambda}^{\bm\alpha,\varepsilon}[\rho]\bigg|_{\pm}=\left(\pm\frac12 I+(\mathcal{K}_{\lambda}^{-\bm\alpha,\varepsilon})^*\right)[\rho].
$$
\item The tangential derivatives are continuous and satisfy
$$
\frac{\partial}{\partial\bm\tau}\mathcal{S}_{\lambda}^{\bm\alpha,\varepsilon}[\rho]\bigg|_+=\frac{\partial}{\partial\bm\tau}\mathcal{S}_{\lambda}^{\bm\alpha,\varepsilon}[\rho]\bigg|_-=\mathcal{T}_{\lambda}^{\bm\alpha,\varepsilon}[\rho].
$$
\end{enumerate}
\end{lemma}

For densities $\phi,\psi\in H^{-1/2}(\partial D^\varepsilon)$, define
$$
u_{\mathrm{out}}(\bm x)=\mathcal{S}_{\omega_{\mathrm m}}^{\bm\alpha,\varepsilon}[\phi](\bm x),  \qquad u_{\mathrm{in}}(\bm x) = \mathcal{S}_{\omega_{\mathrm c}}^{\bm\alpha,\varepsilon}[\psi](\bm x).
$$
Lemma~\ref{lem:lem1} and the interface conditions in~\eqref{eq:system} imply that for $\bm x\in \partial D^\varepsilon$, $(\phi,\psi)$ satisfies
\begin{align}
&\mathcal{S}_{\omega_{\mathrm m}}^{\bm\alpha,\varepsilon}[\phi](\bm x)-\mathcal{S}_{\omega_{\mathrm c}}^{\bm\alpha,\varepsilon}[\psi](\bm x)=0, \label{eq:bie-continuity}\\
&\frac{\mu}{\mu^2-\kappa_g^2}\left(-\frac12 I+(\mathcal{K}_{\omega_{\mathrm c}}^{-\bm\alpha,\varepsilon})^*\right)[\psi](\bm x)+\frac{\ii\kappa_g}{\mu^2-\kappa_g^2}\mathcal{T}_{\omega_{\mathrm c}}^{\bm\alpha,\varepsilon}[\psi](\bm x)-\mu_{\mathrm m}^{-1}\left(\frac12 I+(\mathcal{K}_{\omega_{\mathrm m}}^{-\bm\alpha,\varepsilon})^*\right)[\phi](\bm x)=0.\label{eq:bie-flux}
\end{align}

Accordingly, define the boundary integral operator {\small
$$
\mathcal A^\varepsilon(\bm\alpha,\omega):=
\begin{pmatrix}
\mathcal{S}_{\omega_{\mathrm m}}^{\bm\alpha,\varepsilon}&-\mathcal{S}_{\omega_{\mathrm c}}^{\bm\alpha,\varepsilon} \\
-\mu_{\mathrm m}^{-1}\left(\dfrac12 I+(\mathcal{K}_{\omega_{\mathrm m}}^{-\bm\alpha,\varepsilon})^*\right) & \displaystyle\frac{\mu}{\mu^2-\kappa_g^2}\left(-\dfrac12 I+(\mathcal{K}_{\omega_{\mathrm c}}^{-\bm\alpha,\varepsilon})^*\right)+\frac{\ii\kappa_g}{\mu^2-\kappa_g^2}\mathcal{T}_{\omega_{\mathrm c}}^{\bm\alpha,\varepsilon}
\end{pmatrix}.
$$
}

For fixed $\varepsilon>0$, let
$$
\mathcal X_\varepsilon:=H^{-1/2}(\partial D^\varepsilon)\times H^{-1/2}(\partial D^\varepsilon),\qquad  \mathcal Y_\varepsilon:=H^{1/2}(\partial D^\varepsilon)\times H^{-1/2}(\partial D^\varepsilon).
$$

For $p\in\{A,B\}$, define the scaled $L^2$ inner product
$$
\langle f,g\rangle_{p,\varepsilon} :=\frac1\varepsilon\int_{\partial D_p^\varepsilon}\overline{f(\bm x)}g(\bm x)d\sigma(\bm x) =\int_0^{2\pi}\overline{f(\bm z_p+\varepsilon\bm r(\theta))} g(\bm z_p+\varepsilon\bm r(\theta)) d\theta,
$$
and set $\langle f,g\rangle_{\partial D^\varepsilon,\varepsilon}:=\sum_{p=A,B}\langle f,g\rangle_{p,\varepsilon}$.

\begin{proposition}\label{prop:pde-bie-equivalence}
Suppose the quasi-periodic Green functions for $\lambda\in\{\omega_{\mathrm m}(\omega),\omega_{\mathrm c}(\omega)\}$ are defined and the operators
$$
\mathcal{S}_{\lambda}^{\bm\alpha,\varepsilon}:H^{-1/2}(\partial D^\varepsilon)\rightarrow H^{1/2}(\partial D^\varepsilon)
$$
are boundedly invertible. Then~\eqref{eq:system} has a nontrivial solution if and only if there is $(\phi,\psi)\in\mathcal X_\varepsilon\setminus\{(0,0)\}$ such that
\begin{equation}\label{eq:integral}
\mathcal A^\varepsilon(\bm\alpha,\omega)
\begin{pmatrix}
\phi\\
\psi
\end{pmatrix}
=0.
\end{equation}
\end{proposition}

\begin{proof}
Let $u$ be a nontrivial solution of \eqref{eq:system}, and set $h=u|_{\partial D^\varepsilon}$. By Lemma~2.92 \cite{AmmariFitzpatrickKangRuizYuZhang2018}, invertibility of the single-layer operators gives unique densities $\phi,\psi\in H^{-1/2}(\partial D^\varepsilon)$ satisfying
$$
\mathcal{S}_{\omega_{\mathrm m}}^{\bm\alpha,\varepsilon}[\phi]=h,\quad  \mathcal{S}_{\omega_{\mathrm c}}^{\bm\alpha,\varepsilon}[\psi]=h.
$$
The corresponding potentials solve the exterior and interior Helmholtz equations and have the same boundary traces as $u$. Hence
$$
\mathcal A^\varepsilon(\bm\alpha,\omega)
\begin{pmatrix}
\phi\\
\psi
\end{pmatrix}
=0.
$$
Since $u\ne0$ and both single-layer operators are injective, $(\phi,\psi)\neq(0,0)$.

Conversely, suppose
$$
(\phi,\psi)\in\ker\mathcal A^\varepsilon(\bm\alpha,\omega)\setminus\{(0,0)\},
$$
and define
$$
u(\bm x)=
\begin{cases}
\mathcal{S}_{\omega_{\mathrm m}}^{\bm\alpha,\varepsilon}[\phi](\bm x), & \bm x\in Y\setminus\overline{D^\varepsilon},\\
\mathcal{S}_{\omega_{\mathrm c}}^{\bm\alpha,\varepsilon}[\psi](\bm x), & \bm x\in D^\varepsilon.
\end{cases}
$$
Both potentials satisfy their respective Helmholtz equations and quasi-periodic conditions. Since solution is continuous across the boundary, $u(\bm x)\big|_+=u(\bm x)\big|_-$ for $\bm x\in\partial D^\varepsilon$. Lemma~\ref{lem:lem1} gives
$$
\left.\partial_{\bm\nu}u\right|_+=\left(\frac12 I+(\mathcal{K}_{\omega_{\mathrm m}}^{-\bm\alpha,\varepsilon})^*\right)[\phi],\qquad  \left.\partial_{\bm\nu}u\right|_-=\left(-\frac12 I+(\mathcal{K}_{\omega_{\mathrm c}}^{-\bm\alpha,\varepsilon})^*\right)[\psi],
$$
and
$$
\left.\partial_{\bm\tau}u\right|_-=\mathcal{T}_{\omega_{\mathrm c}}^{\bm\alpha,\varepsilon}[\psi].
$$
Thus for $\bm x\in\partial D^\varepsilon$, \eqref{eq:bie-flux} is equivalent to
$$
\mu_{\mathrm m}^{-1}\left.\partial_{\bm\nu}u\right|_+=\frac{\mu}{\mu^2-\kappa_g^2}\left.\partial_{\bm\nu}u\right|_-+\frac{\ii\kappa_g}{\mu^2-\kappa_g^2}\left.\partial_{\bm\tau}u\right|_-
$$
Consequently, $u$ satisfies~\eqref{eq:system}. If $u\equiv0$, its boundary traces give $\mathcal{S}_{\omega_{\mathrm m}}^{\bm\alpha,\varepsilon}[\phi](\bm x)=\mathcal{S}_{\omega_{\mathrm c}}^{\bm\alpha,\varepsilon}[\psi](\bm x)=0$ for $\bm x\in\partial D^\varepsilon$. Injectivity of the two single-layer operators would imply $\phi=\psi=0$, which is a contradiction. Therefore $u$ is a nontrivial solution of \eqref{eq:system}.
\end{proof}
The standard mapping properties of the three layer-potential operators give $\mathcal A^\varepsilon(\bm\alpha,\omega):\mathcal X_\varepsilon\to\mathcal Y_\varepsilon$.

Define
$$
\varphi_n(\theta):=\frac{1}{\sqrt{2\pi}}e^{\ii n\theta},\quad  n\in\mathbb Z.
$$
Then $\{\varphi_n\}$ is an orthonormal basis of $L^2(\partial D_1)$, equivalently $L^2(0,2\pi)$. On $\partial D^\varepsilon$, set
$$
\varphi_{p,n}(\bm x):=
\begin{cases}
\varphi_n(\theta), & \bm x=\bm z_p+\varepsilon\bm r(\theta)\in\partial D_p^\varepsilon,\\
0, & \bm x\in\partial D^\varepsilon\setminus \partial D_p^\varepsilon.
\end{cases}
$$

The family $\{\varphi_{p,n}:p=A,B, n\in\mathbb Z\}$ is orthonormal for the scaled inner product $\langle\cdot,\cdot\rangle_{\partial D^\varepsilon,\varepsilon}$ with $\langle\varphi_{p,m},\varphi_{q,n}\rangle_{\partial D^\varepsilon,\varepsilon}=\delta_{pq}\delta_{mn}$.

Expand the exterior and interior densities as
$$
\phi=\sum_{p=A,B}\sum_{n\in\mathbb Z}c_{1,p,n}\varphi_{p,n},\qquad  \psi=\sum_{p=A,B}\sum_{n\in\mathbb Z}c_{2,p,n}\varphi_{p,n},
$$
where $c_{1,p,n}:=\langle\varphi_{p,n},\phi\rangle_{p,\varepsilon}, \qquad  c_{2,p,n}:=\langle\varphi_{p,n},\psi\rangle_{p,\varepsilon}$. For distributional densities, these Fourier pairings are understood by duality with the smooth Fourier modes.

Apply the operator $\mathcal A^\varepsilon(\bm\alpha,\omega)=\bigl(\mathcal A_{r\sigma}^\varepsilon(\bm\alpha,\omega)\bigr)_{r,\sigma=1}^2$ to the Fourier basis $\{\varphi_{p,n}\}$. For $r,\sigma\in\{1,2\}$, $p,q\in\{A,B\}$, and $m,n\in\mathbb Z$, define its Fourier matrix elements by
$$
\mathcal A^{\varepsilon,r,\sigma}_{p,m;q,n}:=\left\langle\varphi_{p,m},(\mathcal A_{r\sigma}^\varepsilon)^{pq}\varphi_{q,n}\right\rangle_{p,\varepsilon},
$$
where $(\mathcal A_{r\sigma}^\varepsilon)^{pq}$ denotes the action of $\mathcal A_{r\sigma}^\varepsilon$ from $\partial D_q^\varepsilon$ to $\partial D_p^\varepsilon$. Substitute the Fourier expansions into~\eqref{eq:integral} and take the coefficient of $\varphi_{p,m}$ in equation $r$ to obtain
\begin{equation}\label{eq:coefficient-system}
\sum_{\sigma=1}^{2} \sum_{q=A,B} \sum_{n\in\mathbb Z} \mathcal A^{\varepsilon,r,\sigma}_{p,m;q,n} c_{\sigma,q,n} =0, \qquad  r=1,2,\quad  p=A,B,\quad  m\in\mathbb Z.
\end{equation}
Thus the boundary integral eigenproblem is equivalent to the infinite-dimensional system for the Fourier coefficients $\{c_{\sigma,p,n}\}$ given in \eqref{eq:coefficient-system}. A frequency $\omega$ is a characteristic frequency precisely when \eqref{eq:coefficient-system} has a nonzero coefficient sequence.

For $s\in\mathbb R$, define
$$
h^s:=\left\{ \{a_{p,n}\}_{p=A,B, n\in\mathbb Z}: \sum_{p=A,B}\sum_{n\in\mathbb Z}(1+n^2)^s|a_{p,n}|^2<\infty \right\}.
$$
The circle Fourier expansion identifies $H^s(\partial D^\varepsilon)$ with $h^s$ for each fixed $\varepsilon>0$.

Let $\mathfrak H:=\ell^2\left(\{1,2\}\times\{A,B\}\times\mathbb Z\right)$. For $U=(\phi,\psi)\in\mathcal X_\varepsilon$ with coefficients $\{c_{1,p,n}\}$ and $\{c_{2,p,n}\}$, set
$$
(\mathcal U_XU)_{\sigma,p,n} :=(1+n^2)^{-1/4}c_{\sigma,p,n}, \qquad  \sigma=1,2.
$$
For $F=(f,g)\in\mathcal Y_\varepsilon$ with coefficients $f_{p,n}$ and $g_{p,n}$, set
$$
(\mathcal U_YF)_{1,p,n}:=(1+n^2)^{1/4}f_{p,n},  \qquad (\mathcal U_YF)_{2,p,n}:=(1+n^2)^{-1/4}g_{p,n}.
$$
Then $\mathcal U_X:\mathcal X_\varepsilon\to\mathfrak H$ and $\mathcal U_Y:\mathcal Y_\varepsilon\to\mathfrak H$ are bounded isomorphisms.

For $p\in\{A,B\}$ and $n\in\mathbb Z$, we refer to the one-dimensional subspace $\mathcal M_{p,n}:=\operatorname{span}\{\varphi_{p,n}\}$ as the \emph{Fourier mode indexed by $(p,n)$}, and abbreviate it as the \emph{mode $(p,n)$} when no ambiguity arises. In the two-density formulation, this mode corresponds to the coefficient subspace
$$
\mathfrak M_{p,n}:= \operatorname{span}\{\mathbf e_{1,p,n},\mathbf e_{2,p,n}\} \subset\mathfrak H,
$$
Here $\mathbf e_{\sigma,p,n}$ denotes the standard basis vector of $\mathfrak H$, defined by
$$
(\mathbf e_{\sigma,p,n})_{\sigma',p',n'} =\delta_{\sigma,\sigma'}\delta_{p,p'}\delta_{n,n'}.
$$
Thus, the notation $(p,n)$ labels a boundary component and an angular Fourier index. The single-layer matrix elements contain the boundary-measure factor $\varepsilon$, whereas the leading normal and tangential terms are $\mathcal{O}(1)$. We therefore divide the first boundary equation by $\varepsilon$. We set $\widehat{\mathcal{S}}_\lambda^{\bm\alpha,\varepsilon} :=\varepsilon^{-1}\mathcal{S}_\lambda^{\bm\alpha,\varepsilon}$. For $U\in\mathcal X_\varepsilon$, write $\mathcal A^\varepsilon(\bm\alpha,\omega)U=(F_1,F_2)$ and define the normalized Hilbert-space coefficient operator by
\begin{equation}\label{eq:weighted-A-fourier-entry}
\mathbb A^\varepsilon(\bm\alpha,\omega):\mathfrak H\longrightarrow\mathfrak H, \qquad  \mathbb A^\varepsilon(\bm\alpha,\omega) \mathcal U_XU :=\mathcal U_Y\bigl(\varepsilon^{-1}F_1,F_2\bigr).
\end{equation}
Since $\varepsilon>0$ and $\mathcal U_X,\mathcal U_Y$ are isomorphisms,
\begin{equation}\label{eq:normalized-kernel-equivalence}
\dim\ker\mathbb A^\varepsilon(\bm\alpha,\omega)=\dim\ker\mathcal A^\varepsilon(\bm\alpha,\omega).
\end{equation}

For $p=q$, the entry $\mathcal A^{\varepsilon,r,\sigma}_{p,m;p,n}$ describes self-interaction on one boundary component, and its normal traces include the $\pm\frac12 I$ jump terms. For $p\ne q$, the two boundaries have positive separation, and the Green function and its first derivatives are smooth on $\partial D_p^\varepsilon\times\partial D_q^\varepsilon$, so the corresponding normal operator has no jump.

\subsection{Decomposition of Green's function}
Define $\Delta_{\mathrm{m}}:=\omega_{\mathrm{m}}^2-|\bm\alpha^*|^2$ and $\Delta_{\mathrm{c}}:=\omega_{\mathrm{c}}^2-|\bm\alpha^*|^2$. Define the resonant reciprocal set $\Lambda_0^*:=\{\bm q\in\Lambda^*:|\bm\alpha^*+\bm q|=|\bm\alpha^*|\}$ and define the finite-rank resonant part by
$$
G_{\Lambda_0^*}^{\lambda}(\bm x) :=\frac1{|Y|} \sum_{\bm q\in\Lambda_0^*} \frac{e^{\ii(\bm\alpha^*+\bm q)\cdot\bm x}} {\lambda^2-|\bm\alpha^*+\bm q|^2}.
$$
$\Lambda_0^*$ consists of three points $\bm q_1:=\bm 0$, $\bm q_2:=-(\bm\kappa_1+\bm\kappa_2)$ and $\bm q_3:=-\bm\kappa_2$. The free-space Green function for the Helmholtz equation is $H_0(\lambda;\bm x):=-\frac{\ii}{4}H_0^{(1)}(\lambda|\bm x|)$. 

To ensure the well-definedness of $G_{\omega_m}^{\bm\alpha^*}$ and $G_{\omega_c}^{\bm\alpha^*}$, we impose the following requirements, that for $s=r,m$, there exists constant $C_0>0$ such that
\begin{equation}\label{eq:boundary-material-nonresonance}
\inf_{\bm q\in\Lambda^*}\left|\omega_s^2-(\bm\alpha^*+\bm q)\cdot(\bm\alpha^*+\bm q)\right|>0,  \qquad \inf_{\bm q\in\Lambda^*\backslash\Lambda_0^*}\left|\omega_s^2-(\bm\alpha^*+\bm q)\cdot(\bm\alpha^*+\bm q)\right|\geq C_0>0,
\end{equation}

Define
\begin{equation}\label{eq:green-three-part}
\widetilde G^\lambda(\bm x) :=G_\lambda^{\bm\alpha^*}(\bm x) -H_0(\lambda;\bm x) -G_{\Lambda_0^*}^{\lambda}(\bm x),
\end{equation}
then $G_\lambda^{\bm\alpha^*}=H_0(\lambda;\cdot)+G_{\Lambda_0^*}^{\lambda}+\widetilde G^\lambda$.

\begin{lemma}\label{lem:resonant-green-local}
Let $Q=2\pi/a_{\mathrm{lat}}$ and $\bm\xi_j=\bm\alpha^*+\bm q_j$, $j=1,2,3$. Then
$$
\sum_{j=1}^3\bm\xi_j=\bm0,  \qquad \sum_{j=1}^3\bm\xi_j\bm\xi_j^\top=\frac{2Q^2}{3}I_2,
$$
and, as $|\bm\xi|\to0$,
$$
\sum_{j=1}^3e^{\ii\bm\xi_j\cdot\bm\xi} =3-\frac{Q^2}{3}|\bm\xi|^2+\mathcal{O}(|\bm\xi|^3).
$$
Consequently,
\begin{equation}\label{eq:resonant-green-expansion}
G_{\Lambda_0^*}^{\lambda}(\bm\xi) =\frac1{|Y|\Delta_\lambda} \left(3-\frac{Q^2}{3}|\bm\xi|^2\right) +\mathcal{O}\left(\frac{|\bm\xi|^3}{|\Delta_\lambda|}\right).
\end{equation}
\end{lemma}

\begin{proof}
Direct calculation shows that
$$
\bm\xi_1=\frac{2\pi}{a_{\mathrm{lat}}} \left(\frac1{\sqrt3},-\frac13\right)^\top,  \qquad \bm\xi_2=\frac{2\pi}{a_{\mathrm{lat}}} \left(-\frac1{\sqrt3},-\frac13\right)^\top, \qquad  \bm\xi_3=\frac{2\pi}{a_{\mathrm{lat}}} \left(0,\frac23\right)^\top
$$
and $\bm\xi_1+\bm\xi_2+\bm\xi_3=\bm0$, $\bm\xi_1\bm\xi_1^\top+\bm\xi_2\bm\xi_2^\top+\bm\xi_3\bm\xi_3^\top=\frac{2Q^2}{3}I_2$. For $j=1,2,3$,
$$
e^{\ii\bm\xi_j\cdot\bm\xi} =1+\ii\bm\xi_j\cdot\bm\xi -\frac12(\bm\xi_j\cdot\bm\xi)^2+\mathcal{O}(|\bm\xi|^3).
$$
Summation cancels the linear terms, and
$$
\sum_{j=1}^3(\bm\xi_j\cdot\bm\xi)^2 =\bm\xi^\top \left(\sum_{j=1}^3\bm\xi_j\bm\xi_j^\top\right)\bm\xi =\frac{2Q^2}{3}|\bm\xi|^2.
$$
Therefore
$$
\sum_{j=1}^3e^{\ii\bm\xi_j\cdot\bm\xi} =3-\frac{Q^2}{3}|\bm\xi|^2+\mathcal{O}(|\bm\xi|^3).
$$
For $\bm q_j\in\Lambda_0^*$, the common denominator is $\lambda^2-|\bm\alpha^*+\bm q_j|^2=\Delta_\lambda$. Thus
$$
G_{\Lambda_0^*}^{\lambda}(\bm\xi) =\frac{1}{|Y|\Delta_\lambda} \sum_{j=1}^3e^{\ii\bm\xi_j\cdot\bm\xi} =\frac1{|Y|\Delta_\lambda} \left(3-\frac{Q^2}{3}|\bm\xi|^2\right) +\mathcal{O} \left(\frac{|\bm\xi|^3}{|\Delta_\lambda|}\right).
$$
\end{proof}

The small-argument expansion of the zeroth-order Hankel function~\cite{Olver2010} gives, as $|\bm x|\to0$,
\begin{equation}\label{eq:hankel-local-expansion}
H_0(\lambda;\bm x) =\frac1{2\pi}\ln|\bm x|+\frac1{2\pi}(\ln\lambda+\gamma_0) +\mathcal{O}(|\bm x|^2|\ln|\bm x||),
\end{equation}
where $\gamma_0:=\gamma-\ln2-\ii\pi/2$ and $\gamma$ is Euler's constant. Termwise differentiation gives
\begin{equation}\label{eq:hankel-gradient-expansion}
\nabla_{\bm x}H_0(\lambda;\bm x) =\frac1{2\pi}\frac{\bm x}{|\bm x|^2} +\mathcal{O}(|\bm x| |\ln|\bm x||).
\end{equation}

At $\bm\alpha=\bm\alpha^*$, abbreviate the operators defined in Section~2 as $\mathcal{S}_\lambda^\varepsilon$, $\mathcal{K}_\lambda^{\varepsilon,*}$, and $\mathcal{T}_\lambda^\varepsilon$. With \eqref{eq:green-three-part}, we have the following decompositions:
\begin{align}
\mathcal{S}_\lambda^\varepsilon&=\mathcal{S}_{H_0,\lambda}^\varepsilon+\mathcal{S}_{\Lambda_0^*,\lambda}^\varepsilon+\widetilde{\mathcal{S}}_\lambda^\varepsilon, \notag\\
\mathcal{K}_\lambda^{\varepsilon,*}&=\mathcal{K}_{H_0,\lambda}^{\varepsilon,*}+\mathcal{K}_{\Lambda_0^*,\lambda}^{\varepsilon,*}+\widetilde{\mathcal{K}}_\lambda^{\varepsilon,*}, \notag\\
\mathcal{T}_\lambda^\varepsilon&=\mathcal{T}_{H_0,\lambda}^\varepsilon+\mathcal{T}_{\Lambda_0^*,\lambda}^\varepsilon+\widetilde{\mathcal{T}}_\lambda^\varepsilon. \label{eq:operator-three-part}
\end{align}
Later in this paper, we will frequently call $\mathcal{S}_{H_0,\lambda}^\varepsilon$ the \emph{singular term}, $\mathcal{S}_{\Lambda_0^*,\lambda}^\varepsilon$ the \emph{finite-rank term}, and $\widetilde{\mathcal{S}}_\lambda^\varepsilon$ the \emph{remainder term}. The terms for $\mathcal{K}_\lambda^{\varepsilon,*}$ and $\mathcal{T}_\lambda^\varepsilon$ are named similarly.

Appendix~A derives the detailed small-radius expansions and Fourier coefficients of the Green function and the associated layer-potential operators.

\section{Symmetry decomposition modulo three}\label{sec:3}

The coefficient equation obtained in Section~2 is infinite-dimensional, but its entries are constrained by the rotational symmetry of the honeycomb array. We exploit this structure at $\bm\alpha=\bm\alpha^*$ before making any finite-dimensional reduction. The essential point is that rotation acts not only on the angular Fourier index but also on the Bloch phase associated with the boundary component. The invariant index is therefore $n+s_p$ modulo three, rather than $n$ alone. We first establish the rotational invariance of the quasi-periodic Green function and derive the resulting selection rule for all three layer-potential operators. This gives a direct-sum decomposition of the boundary spaces and of $\mathbb A^\varepsilon$ into three invariant sectors. We then use inversion about the midpoint of the two sites, combined with complex conjugation, to relate the sector restrictions in the nonmagnetic case. The resulting antilinear map exchanges sectors $0$ and $1$ and establishes equality of their real characteristic frequencies, while sector $2$ is mapped to itself. These identities separate the symmetry mechanism responsible for the paired frequencies from the estimates needed to prove their existence and multiplicity. They also prepare the magnetic analysis: rotation continues to preserve the sectors, whereas inversion--conjugation will relate opposite magnetizations rather than the two sectors at a fixed nonzero bias.

Let $R$ be counterclockwise rotation with rotation angle $\vartheta:=2\pi/3$:
$$
R=
\begin{pmatrix}
-\frac12&-\frac{\sqrt3}{2}\\
\frac{\sqrt3}{2}&-\frac12
\end{pmatrix}.
$$
The lattice vectors satisfy $ R\bm e_1=-\bm e_2, R\bm e_2=\bm e_1-\bm e_2. $ For $p\in\{A,B\}$, define $\bm\ell_p\in\Lambda$ by $\bm\ell_p:=R\bm z_p-\bm z_p$.
Since $\bm z_A=\bm0$, $\bm z_B=(\bm e_1+\bm e_2)/3$, we have $\bm\ell_A=\bm0$, $\bm\ell_B=-\bm e_2$. Let $s_p\in\{0,1,2\}$, $p=A, B$ be constants satisfying
\begin{equation}\label{eq:sublattice-phase-index}
e^{\ii\bm\alpha^*\cdot\bm\ell_p}=e^{\ii\vartheta s_p}.
\end{equation}
In our specified coordinate system, $s_A=0$, $s_B=1$.

On the reciprocal lattice, $R^{\mathsf \top}\bm\kappa_1=\bm\kappa_2$, $R^{\mathsf\top}\bm\kappa_2=-(\bm\kappa_1+\bm\kappa_2)$, so $R^{\mathsf \top}\bm\alpha^*=\bm\alpha^*-(\bm\kappa_1+\bm\kappa_2)$. Thus the reciprocal shift in Lemma~\ref{lem:green-C3-invariance} is $\bm q_R=-(\bm\kappa_1+\bm\kappa_2)$. Consequently $u(R\bm x)$ has the same $\bm\alpha^*$-Bloch conditions as $u$. Moreover, $RY$ is another fundamental domain of $\Lambda$, so
$$
\int_Y\overline{u(R\bm x)}v(R\bm x) d\bm x =\int_{RY}\overline{u(\bm x)}v(\bm x) d\bm x =\int_Y\overline u v d\bm x.
$$

\begin{lemma}\label{lem:green-C3-invariance}
Suppose $R^\top\bm\alpha^*=\bm\alpha^*+\bm q_R$, with $\bm q_R\in\Lambda^*$, and $\lambda^2\ne|\bm\alpha^*+\bm q|^2$ for every $\bm q\in\Lambda^*$. Then
\begin{equation}\label{eq:green-C3-invariance}
G_\lambda^{\bm\alpha^*}(R\bm z)=G_\lambda^{\bm\alpha^*}(\bm z).
\end{equation}
\end{lemma}

\begin{proof}
Since $R^\top\Lambda^*=\Lambda^*$ and $R^\top\bm\alpha^*=\bm\alpha^*+\bm q_R$, the map $\bm q\mapsto\bm q'=\bm q_R+R^\top\bm q$ is a bijection of $\Lambda^*$. Therefore
\begin{align*}
G_\lambda^{\bm\alpha^*}(R\bm z) &=\frac1{|Y|}\sum_{\bm q\in\Lambda^*} \frac{e^{\ii(\bm\alpha^*+\bm q)\cdot R\bm z}} {\lambda^2-|\bm\alpha^*+\bm q|^2} =\frac1{|Y|}\sum_{\bm q\in\Lambda^*} \frac{e^{\ii R^\top(\bm\alpha^*+\bm q)\cdot\bm z}} {\lambda^2-|R^\top(\bm\alpha^*+\bm q)|^2}\\
&=\frac1{|Y|}\sum_{\bm q'\in\Lambda^*} \frac{e^{\ii(\bm\alpha^*+\bm q')\cdot\bm z}} {\lambda^2-|\bm\alpha^*+\bm q'|^2} =G_\lambda^{\bm\alpha^*}(\bm z),
\end{align*}
The Fourier sums here are identities of Bloch distributions. Away from the lattice singularities, elliptic regularity gives the same identity for the smooth kernels. Thus no rearrangement of an absolutely convergent pointwise series is presumed. This proves~\eqref{eq:green-C3-invariance}.
\end{proof}

Now we introduce a few important notations, which are crucial for the symmetry decomposition. Using the local basis $\varphi_{p,n}$, define, for $j\in\{0,1,2\}$,
\begin{equation}\label{eq:Hj-definition}
\mathcal H_s^{(j)}:=\overline{\operatorname{span}}\{\varphi_{p,n}:n+s_p\equiv j\pmod3\}, \quad  p\in\{A,B\},\quad  n\in\mathbb Z.
\end{equation}
The closure is taken with respect to the $H^s(\partial D^\varepsilon)$ norm. The associated two-density domain and range spaces are
$$
\mathcal V^{(j)}:=\mathcal H_{-1/2}^{(j)}\oplus\mathcal H_{-1/2}^{(j)}, \qquad  \mathcal W^{(j)}:=\mathcal H_{1/2}^{(j)}\oplus\mathcal H_{-1/2}^{(j)}.
$$
We call $\mathcal H_s^{(j)}\subset H^s(\partial D^\varepsilon)$ \emph{boundary-function sector $j$}, $\mathcal V^{(j)}\subset\mathcal X_\varepsilon$ \emph{two-density domain sector $j$}, and $\mathcal W^{(j)}\subset\mathcal Y_\varepsilon$ \emph{two-density range sector $j$}. In subsequent discussion we often suppress $s$ and only write $\mathcal H^{(j)}$, when there is no ambiguity. It is immediate that $\varphi_{A,n}$ belongs to boundary-function sector $n\pmod3$, whereas $\varphi_{B,n}$ belongs to boundary-function sector $n+1\pmod3$.

For every $s\in\mathbb R$, the physical rotation is the bounded map
\begin{align*}
&\mathcal R_s :H^s(\partial D^\varepsilon)\longrightarrow H^s(\partial D^\varepsilon),\\
&(\mathcal R_s f)(\bm x_p(\theta)) :=e^{\ii\vartheta s_p}f(\bm x_p(\theta+\vartheta)).
\end{align*}
It is immediate that $\mathcal R_s\varphi_{p,n}=e^{\ii\vartheta(n+s_p)}\varphi_{p,n}$, $\mathcal R_s^3=I$. In particular,
$$
\operatorname{ran}\left(\sum_{r=0}^2e^{-\ii\vartheta jr}\mathcal R_s^r\right) =\ker(\mathcal R_s-e^{\ii\vartheta j} I)=\mathcal H_s^{(j)}.
$$

On the coefficient Hilbert space $\mathfrak H$ of Section~2, define $\mathcal R:\mathfrak H\to\mathfrak H$ by
\begin{equation}\label{eq:coefficient-rotation-action}
(\mathcal RC)_{\sigma,p,n}:=e^{\ii\vartheta(n+s_p)}C_{\sigma,p,n}.
\end{equation}
Call $\mathfrak H^{(j)}$ \emph{coefficient sector $j$}, which represents the closed eigenspace of $\mathcal R$ for eigenvalue $e^{\ii\vartheta j}$. With $\mathcal U_X$ and $\mathcal U_Y$ we have
$$
\mathcal U_X\mathcal V^{(j)} = \mathfrak H^{(j)} = \mathcal U_Y\mathcal W^{(j)}.
$$
Therefore for brevity, we will simply call $\mathcal H_s^{(j)}, \mathcal V^{(j)}, \mathcal W^{(j)}$ and $\mathfrak H^{(j)}$ \emph{sector $j$}, when the underlying domain or range space is clear from context.

\begin{theorem}[Invariant subspaces modulo three]\label{thm:mod3-invariant-subspaces}
For $j\in\{0,1,2\}$, the three layer-potential operators satisfy $ \mathcal{S}_\lambda^{\bm\alpha^*,\varepsilon}\mathcal H_{-1/2}^{(j)} \subset\mathcal H_{1/2}^{(j)}, (\mathcal{K}_\lambda^{-\bm\alpha^*,\varepsilon})^*\mathcal H_{-1/2}^{(j)} \subset\mathcal H_{-1/2}^{(j)}, \mathcal{T}_\lambda^{\bm\alpha^*,\varepsilon}\mathcal H_{-1/2}^{(j)} \subset\mathcal H_{-1/2}^{(j)}. $ Consequently, the two-density boundary operator satisfies
$$
\mathcal A^\varepsilon(\bm\alpha^*,\omega)\mathcal V^{(j)} \subset\mathcal W^{(j)}, \qquad  j=0,1,2.
$$
\end{theorem}

\begin{proof}
Recall $\bm\ell_p:=R\bm x_p(\theta)-\bm x_p(\theta+\vartheta)$. Lemma~\ref{lem:green-C3-invariance} and Bloch quasi-periodicity give
$$
G_\lambda^{\bm\alpha^*}(\bm x_p(\theta)-\bm x_q(t)) =e^{\ii\bm\alpha^*\cdot(\bm\ell_p-\bm\ell_q)} G_\lambda^{\bm\alpha^*}(\bm x_p(\theta+\vartheta)-\bm x_q(t+\vartheta)).
$$
For the normal and tangential kernels, use $\bm\nu_p(\theta+\vartheta)=R\bm\nu_p(\theta)$, $\bm\tau_p(\theta+\vartheta)=R\bm\tau_p(\theta)$ and the chain rule
$$
\nabla_{\bm x_p(\theta)}G_\lambda^{\bm\alpha^*} (\bm x_p(\theta)-\bm x_q(t)) =e^{\ii\bm\alpha^*\cdot(\bm\ell_p-\bm\ell_q)} R^\top \nabla_{\bm x_p(\theta+\vartheta)}G_\lambda^{\bm\alpha^*} (\bm x_p(\theta+\vartheta)-\bm x_q(t+\vartheta)).
$$
Thus the kernels of $\mathcal{S}$, $\mathcal{K}^*$, and $\mathcal{T}$ from component $q$ to component $p$ obey the same phase relation. For any Fourier entry $(F^{pq,\varepsilon})_{m,n}$ defined in \eqref{eq:matrix-element-S}, the change of variables $\theta'=\theta+\vartheta$, $t'=t+\vartheta$ gives
\begin{equation}\label{eq:selection-phase-equation}
(F^{pq,\varepsilon})_{m,n}=e^{\ii\bm\alpha^*\cdot(\bm\ell_p-\bm\ell_q)}e^{\ii(m-n)\vartheta}(F^{pq,\varepsilon})_{m,n}.
\end{equation}
When $(F^{pq,\varepsilon})_{m,n}\ne0$, equations~\eqref{eq:sublattice-phase-index} and~\eqref{eq:selection-phase-equation} imply
$$
m+s_p\equiv n+s_q\pmod3.
$$
If $\varphi_{q,n}\in\mathcal H_{-1/2}^{(j)}$, then $n+s_q\equiv j\pmod3$. Every nonzero output coefficient therefore satisfies $m+s_p\equiv j\pmod3$. Thus $ \mathcal{S}_\lambda^{\bm\alpha^*,\varepsilon}\mathcal H_{-1/2}^{(j)} \subset\mathcal H_{1/2}^{(j)}, (\mathcal{K}_\lambda^{-\bm\alpha^*,\varepsilon})^*\mathcal H_{-1/2}^{(j)} \subset\mathcal H_{-1/2}^{(j)}, \mathcal{T}_\lambda^{\bm\alpha^*,\varepsilon}\mathcal H_{-1/2}^{(j)} \subset\mathcal H_{-1/2}^{(j)}. $

Every entry of $\mathcal A^\varepsilon(\bm\alpha^*,\omega)$ is a linear combination of $\mathcal{S}$, $\mathcal{K}^*$, $\mathcal{T}$, and the identity with constant coefficients. The identity also preserves $p,n$, and hence $n+s_p\bmod3$. Therefore
$$
\mathcal A^\varepsilon(\bm\alpha^*,\omega)\mathcal V^{(j)} \subset \mathcal W^{(j)},\qquad  j=0,1,2.
$$
\end{proof}

The modulo-three decomposition converts rotational symmetry at $K$ point into block diagonality, which reduces \eqref{eq:coefficient-system} from the infinite-dimensional boundary equation to finite-dimensional equation system.

By~\eqref{eq:Hj-definition}, the functional spaces of interest can be decomposed as
\begin{align}
&H^s(\partial D^\varepsilon) =\mathcal H_s^{(0)}\oplus\mathcal H_s^{(1)}\oplus\mathcal H_s^{(2)}, \notag\\
&\mathcal X_\varepsilon =\mathcal V^{(0)}\oplus\mathcal V^{(1)}\oplus\mathcal V^{(2)},  \qquad \mathcal Y_\varepsilon =\mathcal W^{(0)}\oplus\mathcal W^{(1)}\oplus\mathcal W^{(2)}. \label{eq:H-direct-sum}
\end{align}
Define the restrictions
$$
\mathcal A_j^\varepsilon(\omega) :=\left.\mathcal A^\varepsilon(\bm\alpha^*,\omega)\right|_{\mathcal V^{(j)}}: \mathcal V^{(j)}\longrightarrow\mathcal W^{(j)}.
$$
Theorem~\ref{thm:mod3-invariant-subspaces} and~\eqref{eq:H-direct-sum} give
$$
\mathcal A^\varepsilon(\bm\alpha^*,\omega) =\mathcal A_0^\varepsilon(\omega)\oplus \mathcal A_1^\varepsilon(\omega)\oplus \mathcal A_2^\varepsilon(\omega),
$$
and the kernel decomposition
\begin{equation}\label{eq:kernel-direct-sum}
\ker\mathcal A^\varepsilon(\bm\alpha^*,\omega) =\ker\mathcal A_0^\varepsilon(\omega)\oplus \ker\mathcal A_1^\varepsilon(\omega)\oplus \ker\mathcal A_2^\varepsilon(\omega).
\end{equation}
Thus a frequency is characteristic precisely when at least one restriction $\mathcal A_j^\varepsilon(\omega)$ has a nonzero kernel.

Now we show that, when $\kappa_g=0$, sector 0 and sector 1 share the same real characteristic frequencies. Using the midpoint $\bm c$ defined in Section~2, set $P\bm x:=2\bm c-\bm x$. For $\bm x_p(\theta)=\bm z_p+\varepsilon\bm r(\theta)$, writing $\overline A=B$ and $\overline B=A$, we have
\begin{equation}\label{eq:P-boundary-param}
P\bm x_p(\theta)=\bm x_{\overline p}(\theta+\pi).
\end{equation}
Define the antilinear operator $\mathcal J:H^s(\partial D^\varepsilon)\to H^s(\partial D^\varepsilon)$ by
\begin{equation}\label{eq:J-definition}
(\mathcal J f)(\bm x):=\overline{f(P\bm x)}, \qquad  \bm x\in\partial D^\varepsilon.
\end{equation}
We suppress $s$ when the space is clear. Equation~\eqref{eq:P-boundary-param} gives
\begin{equation}\label{eq:J-mapping}
\mathcal J\varphi_{p,n}=(-1)^n\varphi_{\overline p,-n}.
\end{equation}
Since $s_A+s_B=1$, for $\varphi_{p,n}\in\mathcal H_s^{(j)}$,
$$
-n+s_{\overline p} \equiv-(n+s_p)+s_p+s_{\overline p} \equiv 1-j\pmod3,
$$
and hence $\mathcal J\mathcal H_s^{(j)}=\mathcal H_s^{(1-j)}$. Define separately
$$
\mathcal J_X:=\mathcal J_{-1/2}\oplus\mathcal J_{-1/2},\qquad  \mathcal J_Y:=\mathcal J_{1/2}\oplus\mathcal J_{-1/2}.
$$
It follows that $\mathcal J_X: \mathcal X_\varepsilon\to \mathcal X_\varepsilon$ and $\mathcal J_Y:\mathcal Y_\varepsilon\to\mathcal Y_\varepsilon$. Thus
\begin{equation}\label{eq:J2-Vj-map}
\mathcal J_X\mathcal V^{(j)}=\mathcal V^{(1-j)}.
\end{equation}

\begin{proposition}\label{prop:conjugate-isospectral}
When $\kappa_g=0$, $\omega\in\mathbb R$,
\begin{equation}\label{eq:J2-commutation}
\mathcal A^\varepsilon(\bm\alpha^*,\omega)\mathcal J_X =\mathcal J_Y\mathcal A^\varepsilon(\bm\alpha^*,\omega).
\end{equation}
If $(\phi,\psi)\in\mathcal V^{(j)}\setminus\{0\}$, $j\in\{0,1,2\}$ satisfies
\begin{equation}\label{eq:j-eigenfunction}
\mathcal A^\varepsilon(\bm\alpha^*,\omega)(\phi,\psi)=0,
\end{equation}
then $\mathcal J_X(\phi,\psi)\in\mathcal V^{(1-j)}\setminus\{0\}$ and
\begin{equation}\label{eq:Sj-eigenfunction}
\mathcal A^\varepsilon(\bm\alpha^*,\omega) \mathcal J_X(\phi,\psi)=0.
\end{equation}
Thus sectors $j$ and $1-j$ (mod 3) have the same real characteristic frequencies.
\end{proposition}

\begin{proof}
For real frequencies, the Green function satisfies
\begin{equation}\label{eq:green-conjugate-inversion}
\overline{G_\lambda^{\bm\alpha^*}(\bm x-\bm y)} =G_\lambda^{\bm\alpha^*}(\bm y-\bm x) =G_\lambda^{\bm\alpha^*}(P\bm x-P\bm y).
\end{equation}
The inversion preserves arc length and maps $\partial D^\varepsilon$ onto itself. Substituting $\bm y=P\bm z$ and using~\eqref{eq:green-conjugate-inversion} gives
$$
\mathcal J \mathcal{S}_\lambda^{\bm\alpha^*,\varepsilon} =\mathcal{S}_\lambda^{\bm\alpha^*,\varepsilon}\mathcal J.
$$
The linear part of $P$ is $-I$, and the outward normal satisfies $\bm\nu(P\bm x)=-\bm\nu(\bm x)$. Differentiating~\eqref{eq:green-conjugate-inversion} in the first variable produces a second minus sign. They cancel, giving
$$
\mathcal J(\mathcal{K}_\lambda^{-\bm\alpha^*,\varepsilon})^* =(\mathcal{K}_\lambda^{-\bm\alpha^*,\varepsilon})^*\mathcal J.
$$
The nonmagnetic operator contains only $\mathcal{S}$, $\mathcal{K}^*$, the identity, and real coefficients, proving~\eqref{eq:J2-commutation}.

Apply $\mathcal J_Y$ to the equation for a nonzero solution of~\eqref{eq:j-eigenfunction} and use~\eqref{eq:J2-commutation}. This gives~\eqref{eq:Sj-eigenfunction}, and from~\eqref{eq:J2-Vj-map} the resulting densities are in $\mathcal V^{(1-j)}$. Since $\mathcal J_X^2=I$, the correspondence is bijective.
\end{proof}

\begin{remark}
For complex frequencies in a domain invariant under conjugation, \eqref{eq:J2-commutation} can be continued as
\begin{equation}\label{eq:J-complex-intertwining}
\mathcal A^\varepsilon(\bm\alpha^*,\overline\omega)\mathcal J_X =\mathcal J_Y\mathcal A^\varepsilon(\bm\alpha^*,\omega).
\end{equation}
With this continuation, bounds and invertibility for sector $0$ transfer to sector $1$ at the conjugate parameter. Thus Proposition~\ref{prop:conjugate-isospectral} pairs sectors $0$ and $1$, while sector $2$ maps to itself. For $\kappa_g\ne0$, the gyrotropic term still preserves all three sectors but generally breaks~\eqref{eq:J2-commutation}, allowing the paired frequencies to split.
\end{remark}

\section{Characteristic frequencies and Dirac cones when $\kappa_g=0$}\label{sec:4}
With the symmetry decomposition established, we now determine the nonmagnetic characteristic frequencies near the two dipole candidates and prove the existence of a Dirac cone. The analysis proceeds from candidate equations to exact spectral statements. We first extract the self-interaction blocks for the modes $(B,-1)$ and $(A,1)$, whose leading determinants identify neighborhoods in which either $\Delta_{\mathrm c}$ or $\Delta_{\mathrm m}$ is of order $\varepsilon^2$. These blocks alone do not determine the spectrum, because their coupling to constant and quadrupole modes can contribute at leading order. We therefore retain a six-dimensional core in each of sectors $0$ and $1$, prove uniform invertibility of the remaining Fourier modes, and perform successive Schur eliminations. The resulting equations exclude the $\Delta_{\mathrm c}$ candidate from the paired sectors and reduce the $\Delta_{\mathrm m}$ candidate to an analytic scalar equation with a simple leading zero. Rouch\'e's theorem, the nonmagnetic symmetry, and a separate analysis of sector $2$ then yield an exact characteristic frequency with a two-dimensional kernel and an explicit small-radius expansion. Finally, we perturb the Bloch vector away from $\bm\alpha^*$. After normalizing the right and left nullvectors by the frequency derivative, we calculate the mixed-sector Bloch derivative of the full reduced equation, including the terms generated by eliminating the other modes. Its nonvanishing leading coefficient produces a two-dimensional system with conical dispersion, completing the passage from the symmetry-protected degeneracy to a local Dirac point.

\subsection{Candidate Dirac frequencies}

Extract the diagonal Fourier blocks of $\mathbb A^\varepsilon(\bm\alpha^*,\omega)$ defined in \eqref{eq:weighted-A-fourier-entry}, with source and target mode $(B,-1)$ in sector $0$ and $(A,1)$ in sector $1$. The resulting two-density blocks are
$$
\mathsf D_0^{(1)} =
\begin{pmatrix}
(\widehat{\mathcal{S}}_{\omega_{\mathrm m}}^{BB,\varepsilon})_{-1,-1} & -(\widehat{\mathcal{S}}_{\omega_{\mathrm c}}^{BB,\varepsilon})_{-1,-1}\\
-\mu_{\mathrm m}^{-1} \left(\dfrac12+(\mathcal{K}_{\omega_{\mathrm m}}^{BB,\varepsilon,*})_{-1,-1}\right) & \mu^{-1} \left(-\dfrac12+(\mathcal{K}_{\omega_{\mathrm c}}^{BB,\varepsilon,*})_{-1,-1}\right)
\end{pmatrix},
$$
$$
\mathsf D_1^{(1)} =
\begin{pmatrix}
(\widehat{\mathcal{S}}_{\omega_{\mathrm m}}^{AA,\varepsilon})_{1,1} & -(\widehat{\mathcal{S}}_{\omega_{\mathrm c}}^{AA,\varepsilon})_{1,1}\\
-\mu_{\mathrm m}^{-1} \left(\dfrac12+(\mathcal{K}_{\omega_{\mathrm m}}^{AA,\varepsilon,*})_{1,1}\right) & \mu^{-1} \left(-\dfrac12+(\mathcal{K}_{\omega_{\mathrm c}}^{AA,\varepsilon,*})_{1,1}\right)
\end{pmatrix}.
$$

\begin{proposition}\label{prop:leading-candidate-determinant}
Consider blocks $\mathsf D_0^{(1)}$ and $\mathsf D_1^{(1)}$ in the region $|\Delta_{\mathrm m}|\ge c\varepsilon^2$, $|\Delta_{\mathrm c}|\ge c\varepsilon^2$. Then $\mathsf D_0^{(1)}$ and $\mathsf D_1^{(1)}$ have the same leading determinant. The zeros of the leading determinant satisfy
\begin{equation}\label{eq:two-leading-candidates}
\Delta_{\mathrm c}=\frac{4\pi Q^2}{3|Y|}\varepsilon^2,  \qquad \Delta_{\mathrm m}=\frac{4\pi Q^2}{3|Y|}\frac{\mu_{\mathrm m}-\mu}{\mu_{\mathrm m}+\mu}\varepsilon^2,
\end{equation}
where $Q:=2\pi/a_{\mathrm{lat}}$ as defined previously.
\end{proposition}

\begin{proof}
By Lemma~\ref{lem:low-self-fourier}, for $j=0,1$,
$$
\mathsf D_j^{(1)} =
\begin{pmatrix}
-\dfrac12+ \dfrac{2\pi Q^2\varepsilon^2}{3|Y|\Delta_{\mathrm m}} & \dfrac12- \dfrac{2\pi Q^2\varepsilon^2}{3|Y|\Delta_{\mathrm c}} \\
-\mu_{\mathrm m}^{-1} \left( \dfrac12+ \dfrac{2\pi Q^2\varepsilon^2}{3|Y|\Delta_{\mathrm m}} \right) & \displaystyle \mu^{-1} \left( -\dfrac12+ \dfrac{2\pi Q^2\varepsilon^2}{3|Y|\Delta_{\mathrm c}} \right)
\end{pmatrix}
+ \mathcal{O}(\varepsilon^2|\ln\varepsilon|).
$$
Set $a_s:=2\pi Q^2\varepsilon^2/(3|Y|\Delta_s)$. The assumed bounds on the two denominators give
$$
\det\mathsf D_j^{(1)} =\left(-\frac12+a_{\mathrm c}\right) \left\{\frac1\mu\left(-\frac12+a_{\mathrm m}\right) -\frac1{\mu_{\mathrm m}}\left(\frac12+a_{\mathrm m}\right)\right\} +\mathcal{O}(\varepsilon^2|\ln\varepsilon|).
$$
The first factor gives $\Delta_{\mathrm c}=4\pi Q^2\varepsilon^2/(3|Y|)$. The second is
$$
-\frac{\mu_{\mathrm m}+\mu}{2\mu\mu_{\mathrm m}} +\frac{\mu_{\mathrm m}-\mu}{\mu\mu_{\mathrm m}}a_{\mathrm m},
$$
whose zero gives the equation for $\Delta_m$ in~\eqref{eq:two-leading-candidates}.
\end{proof}

\begin{remark}
These are only the leading zeros of the dipole two-density determinant. Other modes in the same symmetry sector can modify the effective characteristic equation through off-diagonal entries. Both candidates must therefore be tested in the full operator.
\end{remark}

\subsection{Low Fourier modes and finite-dimensional reduction near the candidates}
Choose a fixed $r>0$ smaller than their absolute values, and consider the two regions
\begin{align}
\left| \Delta_{\mathrm c} -\frac{4\pi Q^2}{3|Y|}\varepsilon^2 \right| &\le r\varepsilon^2, & |\Delta_{\mathrm m}| &\ge c_{\mathrm m}>0, \label{eq:delta-c-candidate-neighborhood}\\
\left| \Delta_{\mathrm m} -\frac{4\pi Q^2}{3|Y|} \frac{\mu_{\mathrm m}-\mu}{\mu_{\mathrm m}+\mu} \varepsilon^2 \right| &\le r\varepsilon^2, & |\Delta_{\mathrm c}| &\ge c_{\mathrm c}>0, \label{eq:delta-m-candidate-neighborhood}
\end{align}
where $c_{\mathrm m}$ and $c_{\mathrm c}$ are fixed constants independent of $\varepsilon$.

For target mode $(p,m)$ and source mode $(q,n)$, the corresponding two-density block of $\mathbb A^\varepsilon(\bm\alpha^*,\omega)$ at $\kappa_g=0$ is
\begin{equation}\label{eq:two-density-fourier-block}
\mathsf B_{p,m;q,n}^{\varepsilon}:=
\begin{pmatrix}
(\widehat{\mathcal{S}}_{\omega_{\mathrm m}}^{pq,\varepsilon})_{m,n}& -(\widehat{\mathcal{S}}_{\omega_{\mathrm c}}^{pq,\varepsilon})_{m,n} \\
-\mu_{\mathrm m}^{-1} \left[ \dfrac12\delta_{pq}\delta_{mn} +(\mathcal{K}_{\omega_{\mathrm m}}^{pq,\varepsilon,*})_{m,n} \right] & \displaystyle \mu^{-1} \left[ -\dfrac12\delta_{pq}\delta_{mn} +(\mathcal{K}_{\omega_{\mathrm c}}^{pq,\varepsilon,*})_{m,n} \right]
\end{pmatrix}.
\end{equation}
Since $\langle\varphi_{p,m},\varphi_{q,n}\rangle =\delta_{pq}\delta_{mn}$, the normal-derivative $\pm\frac12I$ jump terms occur only in diagonal blocks with $(p,m)=(q,n)$.

\begin{lemma}\label{lem:class0-low-cross-blocks}
Write $\rho_s^\varepsilon:= 4\pi^2\varepsilon/(a_{\rm lat}|Y|\Delta_s)$, $s\in\{\mathrm m,\mathrm c\}$. In either~\eqref{eq:delta-c-candidate-neighborhood} or \eqref{eq:delta-m-candidate-neighborhood}, the off-diagonal blocks among the three modes $(A,0)$, $(B,-1)$, and $(B,2)$ in sector $0$ satisfy
\begin{align*}
&\mathsf B_{A,0;B,-1}^{\varepsilon}=
\begin{pmatrix}
-\rho_{\mathrm m}^\varepsilon+\mathcal{O}(\varepsilon)& \rho_{\mathrm c}^\varepsilon+\mathcal{O}(\varepsilon)\\
\mathcal{O}(\varepsilon)&\mathcal{O}(\varepsilon)
\end{pmatrix},\\
&\mathsf B_{B,-1;A,0}^{\varepsilon}=
\begin{pmatrix}
-\rho_{\mathrm m}^\varepsilon+\mathcal{O}(\varepsilon)& \rho_{\mathrm c}^\varepsilon+\mathcal{O}(\varepsilon)\\
\mu_{\mathrm m}^{-1}\rho_{\mathrm m}^\varepsilon+\mathcal{O}(\varepsilon)& -\mu^{-1}\rho_{\mathrm c}^\varepsilon+\mathcal{O}(\varepsilon)
\end{pmatrix},\\
&\mathsf B_{A,0;B,2}^{\varepsilon} =
\begin{pmatrix}
\dfrac{\pi Q^2\varepsilon^2}{3|Y|\Delta_{\mathrm m}} & -\dfrac{\pi Q^2\varepsilon^2}{3|Y|\Delta_{\mathrm c}}\\
0&0
\end{pmatrix}
+\mathcal{O}(\varepsilon^2),\\
&\mathsf B_{B,2;A,0}^{\varepsilon} =
\begin{pmatrix}
\dfrac{\pi Q^2\varepsilon^2}{3|Y|\Delta_{\mathrm m}} & -\dfrac{\pi Q^2\varepsilon^2}{3|Y|\Delta_{\mathrm c}}\\
-\dfrac{2\pi Q^2\varepsilon^2} {3\mu_{\mathrm m}|Y|\Delta_{\mathrm m}} & \dfrac{2\pi Q^2\varepsilon^2} {3\mu|Y|\Delta_{\mathrm c}}
\end{pmatrix}
+\mathcal{O}(\varepsilon^2),\\
&\mathsf B_{B,-1;B,2}^{\varepsilon} =
\begin{pmatrix}
-\dfrac{3\pi k_D^3\varepsilon^3}{8|Y|\Delta_{\mathrm m}} & \dfrac{3\pi k_D^3\varepsilon^3}{8|Y|\Delta_{\mathrm c}}\\
\dfrac{3\pi k_D^3\varepsilon^3} {8\mu_{\mathrm m}|Y|\Delta_{\mathrm m}} & -\dfrac{3\pi k_D^3\varepsilon^3} {8\mu|Y|\Delta_{\mathrm c}}
\end{pmatrix}
+\mathcal{O}(\varepsilon^2|\ln\varepsilon|),\\
&\mathsf B_{B,2;B,-1}^{\varepsilon} =
\begin{pmatrix}
-\dfrac{3\pi k_D^3\varepsilon^3}{8|Y|\Delta_{\mathrm m}} & \dfrac{3\pi k_D^3\varepsilon^3}{8|Y|\Delta_{\mathrm c}}\\
\dfrac{3\pi k_D^3\varepsilon^3} {4\mu_{\mathrm m}|Y|\Delta_{\mathrm m}} & -\dfrac{3\pi k_D^3\varepsilon^3} {4\mu|Y|\Delta_{\mathrm c}}
\end{pmatrix}
+\mathcal{O}(\varepsilon^2|\ln\varepsilon|).
\end{align*}
\end{lemma}

\begin{proof}
First consider the entries between $(A,0)$ and $(B,-1)$. Lemma~\ref{lem:low-resonant-phases} gives $\Xi_{AB}(-1)=-3\ii$, $\Xi_{BA}(1)=3\ii$. Using also $J_0(z_D)=1+\mathcal{O}(\varepsilon^2)$, $J_1(z_D)=z_D/2+\mathcal{O}(\varepsilon^3)$ and $J_1'(z_D)=1/2+\mathcal{O}(\varepsilon^2)$, equations~\eqref{eq:exact-res-S}--\eqref{eq:exact-res-K} give
\begin{align*}
(\widehat{\mathcal{S}}_{\omega_s}^{AB,\varepsilon})_{0,-1} &= (\widehat{\mathcal{S}}_{\omega_s}^{BA,\varepsilon})_{-1,0} = -\frac{6\pi}{|Y|\Delta_s}J_0(z_D)J_1(z_D)+\mathcal{O}(\varepsilon)\\
&= -\frac{3\pi k_D\varepsilon}{|Y|\Delta_s}+\mathcal{O}(\varepsilon) = -\frac{4\pi^2\varepsilon} {a_{\rm lat}|Y|\Delta_s}+\mathcal{O}(\varepsilon) = -\rho_s^\varepsilon+\mathcal{O}(\varepsilon),\\
(\mathcal{K}_{\omega_s}^{BA,\varepsilon,*})_{-1,0} &= -\frac{6\pi\varepsilon k_D}{|Y|\Delta_s} J_1'(z_D)J_0(z_D)+\mathcal{O}(\varepsilon) = -\rho_s^\varepsilon+\mathcal{O}(\varepsilon),\\
(\mathcal{K}_{\omega_s}^{AB,\varepsilon,*})_{0,-1} &=\mathcal{O}(\varepsilon).
\end{align*}
Substituting into~\eqref{eq:two-density-fourier-block} proves the first two blocks.

Next consider $(A,0)$ and $(B,2)$. Using $\Xi_{BA}(-2)=\Xi_{AB}(2)=-3$, $J_2(z_D)=z_D^2/8+\mathcal{O}(\varepsilon^4)$, $J_2'(z_D)=z_D/4+\mathcal{O}(\varepsilon^3)$, and $J_0'(z_D)=-J_1(z_D)$ gives
\begin{align*}
(\widehat{\mathcal{S}}_{\omega_s}^{BA,\varepsilon})_{2,0} &=(\widehat{\mathcal{S}}_{\omega_s}^{AB,\varepsilon})_{0,2} =\frac{6\pi}{|Y|\Delta_s}J_2(z_D)J_0(z_D)+\mathcal{O}(\varepsilon^2)\\
&=\frac{6\pi}{|Y|\Delta_s} \left( \frac{k_D^2\varepsilon^2}{8}+\mathcal{O}(\varepsilon^4) \right) +\mathcal{O}(\varepsilon^2) =\frac{\pi Q^2\varepsilon^2}{3|Y|\Delta_s} +\mathcal{O}(\varepsilon^2),
\end{align*}
\begin{align*}
(\mathcal{K}_{\omega_s}^{BA,\varepsilon,*})_{2,0} &= \frac{6\pi\varepsilon k_D}{|Y|\Delta_s} J_2'(z_D)J_0(z_D)+\mathcal{O}(\varepsilon^2)\\
&= \frac{6\pi\varepsilon k_D}{|Y|\Delta_s} \left( \frac{k_D\varepsilon}{4}+\mathcal{O}(\varepsilon^3) \right) +\mathcal{O}(\varepsilon^2) = \frac{2\pi Q^2\varepsilon^2}{3|Y|\Delta_s} +\mathcal{O}(\varepsilon^2),\\
(\mathcal{K}_{\omega_s}^{AB,\varepsilon,*})_{0,2} &= \frac{6\pi\varepsilon k_D}{|Y|\Delta_s} J_0'(z_D)J_2(z_D)+\mathcal{O}(\varepsilon^2) = \mathcal{O}(\varepsilon^2).
\end{align*}
Substitution into~\eqref{eq:two-density-fourier-block} gives the next two blocks.

Finally, $(B,-1)$ and $(B,2)$ lie on the same boundary component but have different Fourier indices. Thus $\delta_{mn}=0$, and no normal-derivative jump term enters these blocks. Using $\Xi_{BB}(3)=-3\ii$, $\Xi_{BB}(-3)=3\ii$, and $J_{-1}=-J_1$, $J_{-1}'=-J_1'$, we obtain
\begin{align*}
(\widehat{\mathcal{S}}_{\omega_s}^{BB,\varepsilon})_{-1,2} &= (\widehat{\mathcal{S}}_{\omega_s}^{BB,\varepsilon})_{2,-1} = -\frac{6\pi}{|Y|\Delta_s}J_1(z_D)J_2(z_D) +\mathcal{O}(\varepsilon^2|\ln\varepsilon|)\\
&= -\frac{3\pi k_D^3\varepsilon^3} {8|Y|\Delta_s} +\mathcal{O}(\varepsilon^2|\ln\varepsilon|),\\
(\mathcal{K}_{\omega_s}^{BB,\varepsilon,*})_{-1,2} &= -\frac{6\pi\varepsilon k_D}{|Y|\Delta_s} J_1'(z_D)J_2(z_D) +\mathcal{O}(\varepsilon^2|\ln\varepsilon|)\\
&= -\frac{3\pi k_D^3\varepsilon^3} {8|Y|\Delta_s} +\mathcal{O}(\varepsilon^2|\ln\varepsilon|),\\
(\mathcal{K}_{\omega_s}^{BB,\varepsilon,*})_{2,-1} &= -\frac{6\pi\varepsilon k_D}{|Y|\Delta_s} J_2'(z_D)J_1(z_D) +\mathcal{O}(\varepsilon^2|\ln\varepsilon|)\\
&= -\frac{3\pi k_D^3\varepsilon^3} {4|Y|\Delta_s} +\mathcal{O}(\varepsilon^2|\ln\varepsilon|).
\end{align*}
Substitution into~\eqref{eq:two-density-fourier-block} gives the remaining two blocks.
\end{proof}

In either candidate-frequency neighborhood, one of $\Delta_{\mathrm m}$ and $\Delta_{\mathrm c}$ is of order $\varepsilon^2$, while the other is bounded away from zero. It follows from Lemma~\ref{lem:class0-low-cross-blocks} that $\mathsf B_{A,0;B,2}^{\varepsilon}, \mathsf B_{B,2;A,0}^{\varepsilon}$ have nonvanishing order-one entries, whereas $\mathsf B_{B,-1;B,2}^{\varepsilon}$ and $\mathsf B_{B,2;B,-1}^{\varepsilon}$ are $\mathcal{O}(\varepsilon)$. Thus the quadrupole mode $(B,2)$ is coupled at leading order to the constant mode $(A,0)$ and cannot be discarded without accounting for this interaction.

For $j=0,1$, let $\Pi_j$ denote the orthogonal projection in $\mathfrak H$ onto the symmetry sector $\mathfrak H^{(j)}=\ker(\mathcal R-e^{2\pi \ii j/3} I)$. Let $P_j:\mathfrak H\rightarrow\mathfrak H_{{\rm core}}^{(j)}$ denote the orthogonal projection onto the six-dimensional coefficient subspace with
\begin{align*}
\mathfrak H_{\mathrm{core}}^{(0)} &:= \operatorname{span} \left\{ \mathbf e_{\sigma,A,0}, \mathbf e_{\sigma,B,-1}, \mathbf e_{\sigma,B,2} :\sigma=1,2 \right\},\\
\mathfrak H_{\mathrm{core}}^{(1)} &:= \operatorname{span} \left\{ \mathbf e_{\sigma,B,0}, \mathbf e_{\sigma,A,1}, \mathbf e_{\sigma,A,-2} :\sigma=1,2 \right\}.
\end{align*}
Set $Q_j:=\Pi_j-P_j$. Then $\mathfrak H^{(j)}=P_j\mathfrak H\oplus Q_j\mathfrak H$, where $Q_j\mathfrak H$ is the orthogonal complement of the core within sector $j$. All remaining Fourier modes outside the core have $|n|\ge 3$, and their coupling to the core is $\mathcal{O}(\varepsilon)$.

\begin{proposition}\label{prop:six-dimensional-reduction}
Let $\omega$ satisfy \eqref{eq:delta-c-candidate-neighborhood} or~\eqref{eq:delta-m-candidate-neighborhood}. There are constants $C>0$ and $\varepsilon_0>0$, independent of $\varepsilon$, such that, for $0<\varepsilon<\varepsilon_0$,
$$
Q_j\mathbb A^\varepsilon(\bm\alpha^*,\omega)Q_j:Q_j\mathfrak H\longrightarrow Q_j\mathfrak H
$$
is invertible and
$$
\left\| \bigl( Q_j\mathbb A^\varepsilon(\bm\alpha^*,\omega)Q_j \bigr)^{-1} \right\| \le C.
$$

Relative to the decomposition $\Pi_j\mathfrak H=P_j\mathfrak H\oplus Q_j\mathfrak H$, the equation $\mathbb A^\varepsilon(\bm\alpha^*,\omega)C=0$ with $C\in\Pi_j\mathfrak H$ is equivalent to the six-dimensional equation
$$
\left[ P_j\mathbb A^\varepsilon P_j - P_j\mathbb A^\varepsilon Q_j \bigl( Q_j\mathbb A^\varepsilon Q_j \bigr)^{-1} Q_j\mathbb A^\varepsilon P_j \right]C_P=0, \qquad  C_P\in P_j\mathfrak H,
$$
where all operators are evaluated at $(\bm\alpha^*,\omega)$. The complementary component is uniquely determined by
$$
C_Q = - \bigl( Q_j\mathbb A^\varepsilon Q_j \bigr)^{-1} Q_j\mathbb A^\varepsilon P_j C_P.
$$

Moreover,
$$
\left\| P_j\mathbb A^\varepsilon Q_j \bigl( Q_j\mathbb A^\varepsilon Q_j \bigr)^{-1} Q_j\mathbb A^\varepsilon P_j \right\| \le C\varepsilon^2.
$$
Thus the kernel of the full sector restriction is isomorphic to the Schur kernel. In particular,
$$
\dim\ker \mathbb A^\varepsilon(\bm\alpha^*,\omega)\big|_{\mathfrak H^{(j)}} = \dim\ker \left[ P_j\mathbb A^\varepsilon P_j - P_j\mathbb A^\varepsilon Q_j \bigl( Q_j\mathbb A^\varepsilon Q_j \bigr)^{-1} Q_j\mathbb A^\varepsilon P_j \right].
$$
\end{proposition}
\begin{proof}
By the sector rule, after removing $(A,0)$, $(B,-1)$, and $(B,2)$ from sector $0$, all remaining modes have $|n|\ge3$. The same is true in sector $1$ after removing $(B,0)$, $(A,1)$, and $(A,-2)$.

Lemmas~\ref{lem:appendix-a-1} and~\ref{lem:appendix-a-3} give $\widehat{\mathcal{S}}_{H_0,\omega_s}^{pp,\varepsilon}\varphi_{p,n} =-(2|n|)^{-1}\varphi_{p,n} +\mathcal{O}(\varepsilon^2|\ln\varepsilon|)$ and $\mathcal{K}_{H_0,\omega_s}^{pp,\varepsilon,*}\varphi_{p,n} =\mathcal{O}(\varepsilon^2|\ln\varepsilon|)$.
For each $|n|\ge3$, the leading free-space contribution to $\mathsf B_{p,m;q,n}^{\varepsilon}$ is
$$
M_n^{(0)} =
\begin{pmatrix}
-\dfrac1{2|n|} & \dfrac1{2|n|}\\
-\dfrac1{2\mu_{\mathrm m}} & -\dfrac1{2\mu}
\end{pmatrix}.
$$

Applying $\mathcal U_X$ and $\mathcal U_Y$ and writing $x_{\sigma,p,n}=(1+n^2)^{-1/4}c_{\sigma,p,n}$ gives
$$
\begin{pmatrix}
(\mathcal U_YF)_{1,p,n}\\
(\mathcal U_YF)_{2,p,n}
\end{pmatrix}
=
\begin{pmatrix}
(1+n^2)^{1/4}f_{p,n}\\
(1+n^2)^{-1/4}g_{p,n}
\end{pmatrix}
=
\begin{pmatrix}
-\dfrac{\sqrt{1+n^2}}{2|n|} & \dfrac{\sqrt{1+n^2}}{2|n|} \\
-\dfrac1{2\mu_{\mathrm m}} & -\dfrac1{2\mu}
\end{pmatrix}
\begin{pmatrix}
x_{1,p,n}\\
x_{2,p,n}
\end{pmatrix}.
$$
With $\eta_n:=\sqrt{1+n^2}/|n|$, the leading block is
$$
\mathsf D_n =
\begin{pmatrix}
-\dfrac{\eta_n}{2}&\dfrac{\eta_n}{2}\\
-\dfrac1{2\mu_{\mathrm m}}&-\dfrac1{2\mu}
\end{pmatrix}
.
$$
For $|n|\ge3$, $1\le\eta_n\le\sqrt{10}/3$, and $\det\mathsf D_n =\eta_n(\mu_{\mathrm m}+\mu)/(4\mu\mu_{\mathrm m})$. Thus $\mathsf D_n$ is invertible, with $\|\mathsf D_n^{-1}\|$ uniformly bounded over the indicated modes. Define
\begin{align*}
\mathbb D_j := \bigoplus_{\substack{p\in\{A,B\},|n|\ge3\\
n+s_p\equiv j   \pmod3}} \mathsf D_n .
\end{align*}
This operator is boundedly invertible, with $\|\mathbb D_j^{-1}\|\le C$.

Use the sectorwise projections $P_j$ and $Q_j=\Pi_j-P_j$ just defined. Equations~\eqref{eq:exact-res-S}--\eqref{eq:exact-res-T}, together with
\begin{equation}\label{eq:bessel-estimate}
|J_n(z_D)| \le C\frac{(C\varepsilon)^{|n|}}{|n|!}, \qquad  2J_n'(z_D)=J_{n-1}(z_D)-J_{n+1}(z_D),
\end{equation}
imply that, for $|m|,|n|\ge3$, the weighted resonant entries satisfy
$$
\left| (\mathcal U_Y\mathcal A_{\Lambda_0^*}^\varepsilon \mathcal U_X^{-1})_{p,m;q,n} \right| \le C(1+|m|+|n|)^3 \frac{(C\varepsilon)^{|m|+|n|-2}} {|m|!|n|!}.
$$
Hence
\begin{align*}
\sum_{\substack{p,q\in\{A,B\}\\
|m|,|n|\ge3}} \left| (\mathcal U_Y\mathcal A_{\Lambda_0^*}^\varepsilon \mathcal U_X^{-1})_{p,m;q,n} \right|^2 &\le C \sum_{|m|,|n|\ge3} (1+|m|+|n|)^6 \frac{(C\varepsilon)^{2(|m|+|n|-2)}} {(|m|!)^2(|n|!)^2} &\le C\varepsilon^8,
\end{align*}
so the resonant operator on $Q_j\mathfrak H$ is $\mathcal{O}(\varepsilon^4)$.

On one boundary component, the $H_0$ kernel depends only on $\theta-t$. The Hankel expansion gives, for $|n|\ge3$,
$$
\left| (\widehat{\mathcal{S}}_{H_0,\omega_s}^{pp,\varepsilon})_{n,n} +\frac1{2|n|} \right| \le C\varepsilon^2(1+|n|)^{-3},
$$
and
$$
\left| (\mathcal{K}_{H_0,\omega_s}^{pp,\varepsilon,*})_{n,n} \right| + \left| (\mathcal{T}_{H_0,\omega_s}^{pp,\varepsilon})_{n,n} -\frac{\ii}{2}\operatorname{sgn}(n) \right| \le C\varepsilon^2|\ln\varepsilon|.
$$
For $\widetilde G^{\omega_s}$ and the smooth cross-component kernels, the weighted Fourier entries satisfy, for every $N\ge1$, the estimate
$$
\left| a_{p,m;q,n}^{\varepsilon} \right| \le C_N\varepsilon^2 (1+|m|+|n|)^{-N}, \qquad  |m|,|n|\ge3.
$$
Taking $N$ sufficiently large and summing gives
$$
\left\| Q_j\mathbb A^\varepsilon(\bm\alpha^*,\omega)Q_j -\mathbb D_j \right\| \le C\varepsilon^2|\ln\varepsilon|.
$$
Combining this with $\|\mathbb D_j^{-1}\|\le C$,
$$
\left\| \mathbb D_j^{-1} \left( Q_j\mathbb A^\varepsilon(\bm\alpha^*,\omega)Q_j -\mathbb D_j \right) \right\| \le C\varepsilon^2|\ln\varepsilon|<1
$$
gives the stated bound for sufficiently small $\varepsilon$. Therefore
\begin{align*}
Q_j\mathbb A^\varepsilon Q_j &= \mathbb D_j \left[ I+ \mathbb D_j^{-1} \left( Q_j\mathbb A^\varepsilon Q_j-\mathbb D_j \right) \right],\\
(Q_j\mathbb A^\varepsilon Q_j)^{-1} &= \left[ I+ \mathbb D_j^{-1} \left( Q_j\mathbb A^\varepsilon Q_j-\mathbb D_j \right) \right]^{-1} \mathbb D_j^{-1},
\end{align*}
and
$$
\left\| (Q_j\mathbb A^\varepsilon Q_j)^{-1} \right\| \le C.
$$

We next estimate the coupling between $P_j\mathfrak H$ and $Q_j\mathfrak H$. Equations~\eqref{eq:exact-res-S}--\eqref{eq:exact-res-T} give
$$
\left\| P_j\mathbb A^\varepsilon(\bm\alpha^*,\omega)Q_j \right\| + \left\| Q_j\mathbb A^\varepsilon(\bm\alpha^*,\omega)P_j \right\| \le C\varepsilon.
$$
Consequently,
$$
\left\| P_j\mathbb A^\varepsilon Q_j (Q_j\mathbb A^\varepsilon Q_j)^{-1} Q_j\mathbb A^\varepsilon P_j \right\| \le \left\|P_j\mathbb A^\varepsilon Q_j\right\| \left\|(Q_j\mathbb A^\varepsilon Q_j)^{-1}\right\| \left\|Q_j\mathbb A^\varepsilon P_j\right\| \le C\varepsilon^2.
$$

Using the decomposition $\Pi_j\mathfrak H=P_j\mathfrak H\oplus Q_j\mathfrak H$, the kernel equation becomes
$$
\begin{pmatrix}
P_j\mathbb A^\varepsilon P_j & P_j\mathbb A^\varepsilon Q_j\\
Q_j\mathbb A^\varepsilon P_j & Q_j\mathbb A^\varepsilon Q_j
\end{pmatrix}
\begin{pmatrix}
C_P\\
C_Q
\end{pmatrix}
=0.
$$
Invertibility of $Q_j\mathbb A^\varepsilon Q_j$ gives
$$
C_Q = -(Q_j\mathbb A^\varepsilon Q_j)^{-1} Q_j\mathbb A^\varepsilon P_j C_P,
$$
and therefore
$$
\left[ P_j\mathbb A^\varepsilon P_j - P_j\mathbb A^\varepsilon Q_j (Q_j\mathbb A^\varepsilon Q_j)^{-1} Q_j\mathbb A^\varepsilon P_j \right]C_P =0.
$$
Thus the kernel equation is equivalent to the six-dimensional Schur equation on $\mathfrak H_{\mathrm{core}}^{(j)}$ and its correction satisfies
$$
\left\| P_j\mathbb A^\varepsilon Q_j (Q_j\mathbb A^\varepsilon Q_j)^{-1} Q_j\mathbb A^\varepsilon P_j \right\| =\mathcal{O}(\varepsilon^2).
$$
\end{proof}
In sectors $0$ and $1$, respectively, use the following orders for the unweighted coefficient vectors in the six-dimensional cores:
$$
(c_{1,A,0},c_{2,A,0},c_{1,B,-1},c_{2,B,-1},c_{1,B,2},c_{2,B,2})^\top
$$
and
$$
(c_{1,B,0},c_{2,B,0},c_{1,A,1},c_{2,A,1},c_{1,A,-2},c_{2,A,-2})^\top.
$$

\begin{lemma}\label{lem:second-mode-self-block}
In either case of~\eqref{eq:delta-c-candidate-neighborhood} or \eqref{eq:delta-m-candidate-neighborhood}, the blocks for $(B,2)$ in sector $0$ and $(A,-2)$ in sector $1$ satisfy
$$
E_j^\varepsilon=
\begin{pmatrix}
-\dfrac{1}{4}&\dfrac{1}{4}\\
-\dfrac{1}{2\mu_{\mathrm m}}&-\dfrac{1}{2\mu}
\end{pmatrix}
+ \mathcal{O}(\varepsilon^2|\ln\varepsilon|), \qquad  j=0,1.
$$
As a result, $E_j^\varepsilon$ is invertible for sufficiently small $\varepsilon$.
\end{lemma}

\begin{proof}
The free-space self-interaction entries give $\widehat{\mathcal{S}}_{H_0,\omega_s}^{pp,\varepsilon}\varphi_{p,\pm2} =-\frac14\varphi_{p,\pm2} +\mathcal{O}(\varepsilon^2|\ln\varepsilon|)$, and $\mathcal{K}_{H_0,\omega_s}^{pp,\varepsilon,*}\varphi_{p,\pm2} =\mathcal{O}(\varepsilon^2|\ln\varepsilon|)$. On the other hand, $\Xi_{pp}(0)=3$, $J_2(z_D)=z_D^2/8+\mathcal{O}(\varepsilon^4)$, and $J_2'(z_D)=z_D/4+\mathcal{O}(\varepsilon^3)$. Equations~\eqref{eq:exact-res-S}--\eqref{eq:exact-res-K} imply
\begin{align*}
(\widehat{\mathcal{S}}_{\Lambda_0^*,\omega_s}^{pp,\varepsilon})_{\pm2,\pm2} &= \frac{6\pi}{|Y|\Delta_s}J_2(z_D)^2 = \mathcal{O} \left(\frac{\varepsilon^4}{|\Delta_s|}\right) = \mathcal{O}(\varepsilon^2),\\
(\mathcal{K}_{\Lambda_0^*,\omega_s}^{pp,\varepsilon,*})_{\pm2,\pm2} &= \frac{6\pi\varepsilon k_D}{|Y|\Delta_s} J_2'(z_D)J_2(z_D) = \mathcal{O} \left(\frac{\varepsilon^4}{|\Delta_s|}\right) = \mathcal{O}(\varepsilon^2).
\end{align*}
The corresponding smooth-remainder entries are also $\mathcal{O}(\varepsilon^2)$. Therefore $(\widehat{\mathcal{S}}_{\omega_s}^{pp,\varepsilon})_{\pm2,\pm2} =-1/4+\mathcal{O}(\varepsilon^2|\ln\varepsilon|)$, $(\mathcal{K}_{\omega_s}^{pp,\varepsilon,*})_{\pm2,\pm2} =\mathcal{O}(\varepsilon^2|\ln\varepsilon|)$ and
\begin{align*}
E_j^\varepsilon &=
\begin{pmatrix}
(\widehat{\mathcal{S}}_{\omega_{\mathrm m}}^{pp,\varepsilon})_{\pm2,\pm2} & -(\widehat{\mathcal{S}}_{\omega_{\mathrm c}}^{pp,\varepsilon})_{\pm2,\pm2}\\
-\mu_{\mathrm m}^{-1} \left(\dfrac12+(\mathcal{K}_{\omega_{\mathrm m}}^{pp,\varepsilon,*})_{\pm2,\pm2}\right) & \mu^{-1} \left(-\dfrac12+(\mathcal{K}_{\omega_{\mathrm c}}^{pp,\varepsilon,*})_{\pm2,\pm2}\right)
\end{pmatrix}
+\mathcal{O}(\varepsilon^2)\\
&=
\begin{pmatrix}
-\dfrac14&\dfrac14\\
-\dfrac1{2\mu_{\mathrm m}}&-\dfrac1{2\mu}
\end{pmatrix}
+ \mathcal{O}(\varepsilon^2|\ln\varepsilon|).
\end{align*}
The leading determinant is $(\mu_{\mathrm m}+\mu)/(8\mu\mu_{\mathrm m})$, which proves the assertion.
\end{proof}

\subsection{Schur equations for the two candidate frequencies}\label{sec:sec4-3}

Proposition~\ref{prop:six-dimensional-reduction} reduces the kernel equations in sectors $0$ and $1$ to the cores in $\mathfrak H_{\mathrm{core}}^{(j)}$, $j=0,1$. In sector $0$, the Schur matrix on $\mathfrak H_{\mathrm{core}}^{(0)}$ has the following $6\times6$ expansion:
$$
\mathbf C_0^\varepsilon:=
\begin{pmatrix}
\mathsf B_{A,0;A,0}^\varepsilon&\mathsf B_{A,0;B,-1}^\varepsilon&\mathsf B_{A,0;B,2}^\varepsilon\\
\mathsf B_{B,-1;A,0}^\varepsilon&\mathsf B_{B,-1;B,-1}^\varepsilon&\mathsf B_{B,-1;B,2}^\varepsilon\\
\mathsf B_{B,2;A,0}^\varepsilon&\mathsf B_{B,2;B,-1}^\varepsilon&E_0^\varepsilon
\end{pmatrix}
+\mathcal{O}(\varepsilon^2).
$$
Each entry is a $2\times2$ block. Sector $1$ has the same structure by Proposition~\ref{prop:conjugate-isospectral}. Eliminating the last two coordinates associated with $E_j^\varepsilon$ gives the four-dimensional matrix
\begin{equation}\label{eq:four-dim-matrix}
\begin{pmatrix}
Z_j^\varepsilon&R_j^\varepsilon\\
L_j^\varepsilon&D_j^\varepsilon
\end{pmatrix}
,  \qquad j=0,1.
\end{equation}
For $j=0$,
\begin{align*}
Z_j^\varepsilon &= \mathsf B_{A,0;A,0}^\varepsilon-\mathsf B_{A,0;B,2}^\varepsilon(E_j^\varepsilon)^{-1}\mathsf B_{B,2;A,0}^\varepsilon, \\
R_j^\varepsilon &= \mathsf B_{A,0;B,-1}^\varepsilon-\mathsf B_{A,0;B,2}^\varepsilon(E_j^\varepsilon)^{-1}\mathsf B_{B,2;B,-1}^\varepsilon, \\
L_j^\varepsilon &= \mathsf B_{B,-1;A,0}^\varepsilon-\mathsf B_{B,-1;B,2}^\varepsilon(E_j^\varepsilon)^{-1}\mathsf B_{B,2;A,0}^\varepsilon, \\
D_j^\varepsilon &= \mathsf B_{B,-1;B,-1}^\varepsilon - \mathsf B_{B,-1;B,2}^\varepsilon(E_j^\varepsilon)^{-1}\mathsf B_{B,2;B,-1}^\varepsilon.
\end{align*}
For $j=1$, we replace modes $(A,0), (B,-1)$ and $(B,2)$ with $(B,0), (A,-1), (A,-2)$. Since the eigenvalues of sector 0 and 1 coincide, we focus on the calcuation for $j=0$. In the subsequent discussions we compute \eqref{eq:four-dim-matrix} in the two candidate neighborhoods \eqref{eq:delta-c-candidate-neighborhood} and \eqref{eq:delta-m-candidate-neighborhood} separately.

\subsubsection{The $\Delta_{\mathrm c}$ candidate}\label{sec:delta-c-candidate}

Let $z_{\mathrm c}:=4\pi Q^2/(3|Y|)$ and $\Delta_{\mathrm c}(\omega)=\varepsilon^2z$. Choose $r_{\mathrm c}>0$ and take $|z-z_{\mathrm c}|<r_{\mathrm c}$. Set $a(z):=2\pi Q^2/(3|Y|z)$.

Lemma~\ref{lem:class0-low-cross-blocks} gives
\begin{align*}
\mathsf B_{A,0;B,2}^\varepsilon &=
\begin{pmatrix}
\mathcal{O}(\varepsilon^2)&-\dfrac{a(z)}2+\mathcal{O}(\varepsilon^2)\\
\mathcal{O}(\varepsilon^2)&\mathcal{O}(\varepsilon^2)
\end{pmatrix},
& \mathsf B_{B,-1;B,2}^\varepsilon&=\mathcal{O}(\varepsilon),\\
\mathsf B_{B,2;A,0}^\varepsilon &=
\begin{pmatrix}
\mathcal{O}(\varepsilon^2)&-\dfrac{a(z)}2+\mathcal{O}(\varepsilon^2)\\
\mathcal{O}(\varepsilon^2)&\dfrac{a(z)}\mu+\mathcal{O}(\varepsilon^2)
\end{pmatrix}, 
& \mathsf B_{B,2;B,-1}^\varepsilon&=\mathcal{O}(\varepsilon).
\end{align*}

Lemma~\ref{lem:second-mode-self-block} also gives
$$
E_0^\varepsilon=E_*+\mathcal{O}(\varepsilon^2|\ln\varepsilon|), \qquad  E_*:=
\begin{pmatrix}
-\dfrac14&\dfrac14\\
-\dfrac1{2\mu_{\mathrm m}}&-\dfrac1{2\mu}
\end{pmatrix}.
$$
Its inverse is
$$
E_*^{-1}=
\begin{pmatrix}
-\dfrac{4\mu_{\mathrm m}}{\mu_{\mathrm m}+\mu}&-\dfrac{2\mu\mu_{\mathrm m}}{\mu_{\mathrm m}+\mu}\\
\dfrac{4\mu}{\mu_{\mathrm m}+\mu}&-\dfrac{2\mu\mu_{\mathrm m}}{\mu_{\mathrm m}+\mu}
\end{pmatrix}.
$$

By Lemma~\ref{lem:second-mode-self-block}, the block $E_0^\varepsilon$ associated with $(B,2)$ is uniformly invertible. Eliminating its two density coefficients therefore produces the following correction on the four coefficients of $(A,0)$ and $(B,-1)$:
$$
\begin{pmatrix}
\mathsf B_{A,0;B,2}^{\varepsilon}\\
\mathsf B_{B,-1;B,2}^{\varepsilon}
\end{pmatrix}
(E_0^\varepsilon)^{-1}
\begin{pmatrix}
\mathsf B_{B,2;A,0}^{\varepsilon} & \mathsf B_{B,2;B,-1}^{\varepsilon}
\end{pmatrix}.
$$
Under $\Delta_{\mathrm c}=\varepsilon^2z$, insert the entries from Lemma~\ref{lem:class0-low-cross-blocks} and the inverse of $E_0^\varepsilon$ to obtain
\begin{align*}
&
\begin{pmatrix}
\mathsf B_{A,0;B,2}^{\varepsilon}\\
\mathsf B_{B,-1;B,2}^{\varepsilon}
\end{pmatrix}
(E_0^\varepsilon)^{-1}
\begin{pmatrix}
\mathsf B_{B,2;A,0}^{\varepsilon} & \mathsf B_{B,2;B,-1}^{\varepsilon}
\end{pmatrix}
\\
&\qquad =
\begin{pmatrix}
0&a(z)^2&0&0\\
0&0&0&0\\
0&0&0&0\\
0&0&0&0
\end{pmatrix}
+
\begin{pmatrix}
\mathcal{O}(\varepsilon^2|\ln\varepsilon|) & \mathcal{O}(\varepsilon^2|\ln\varepsilon|) & \mathcal{O}(\varepsilon)&\mathcal{O}(\varepsilon) \\
\mathcal{O}(\varepsilon^2|\ln\varepsilon|) & \mathcal{O}(\varepsilon^2|\ln\varepsilon|) & \mathcal{O}(\varepsilon)&\mathcal{O}(\varepsilon) \\
\mathcal{O}(\varepsilon)&\mathcal{O}(\varepsilon) & \mathcal{O}(\varepsilon^2)&\mathcal{O}(\varepsilon^2) \\
\mathcal{O}(\varepsilon)&\mathcal{O}(\varepsilon) & \mathcal{O}(\varepsilon^2)&\mathcal{O}(\varepsilon^2)
\end{pmatrix}
.
\end{align*}

Thus $Z_0^\varepsilon$ after this elimination is
\begin{align*}
Z_0^\varepsilon &=
\begin{pmatrix}
(\widehat{\mathcal{S}}_{\omega_{\mathrm m}}^{pp,\varepsilon})_{0,0}&-(\widehat{\mathcal{S}}_{\omega_{\mathrm c}}^{pp,\varepsilon})_{0,0}-a(z)^2\\
-\mu_{\mathrm m}^{-1} \left(\dfrac12+(\mathcal{K}_{\omega_{\mathrm m}}^{*,pp,\varepsilon})_{0,0}\right)& \mu^{-1} \left(-\dfrac12+(\mathcal{K}_{\omega_{\mathrm c}}^{*,pp,\varepsilon})_{0,0}\right)
\end{pmatrix}
+\mathcal{O}(\varepsilon^2|\ln\varepsilon|)\\
&=
\begin{pmatrix}
\mathcal{O}(|\ln\varepsilon|)&-\dfrac{6\pi}{|Y|\varepsilon^2z}+\mathcal{O}(|\ln\varepsilon|+1)\\
-\mu_{\mathrm m}^{-1}+\mathcal{O}(\varepsilon^2|\ln\varepsilon|)& -\dfrac{2a(z)}\mu+\mathcal{O}(\varepsilon^2|\ln\varepsilon|)
\end{pmatrix}.
\end{align*}
Consequently,
$$
\det Z_0^\varepsilon =-\frac{6\pi}{\mu_{\mathrm m}|Y|\varepsilon^2z}+\mathcal{O}(|\ln\varepsilon|),
$$
and
\begin{align*}
\bigl[(Z_0^\varepsilon)^{-1}\bigr]_{11}&=\mathcal{O}(\varepsilon^2),& \bigl[(Z_0^\varepsilon)^{-1}\bigr]_{12}&=-\mu_{\mathrm m}+\mathcal{O}(\varepsilon^2|\ln\varepsilon|),\\
\bigl[(Z_0^\varepsilon)^{-1}\bigr]_{21}&=-\frac{|Y|\varepsilon^2z}{6\pi}+\mathcal{O}(\varepsilon^4|\ln\varepsilon|),& \bigl[(Z_0^\varepsilon)^{-1}\bigr]_{22}&=\mathcal{O}(\varepsilon^2|\ln\varepsilon|).
\end{align*}

The same elimination changes the blocks $R_0^\varepsilon$ and $L_0^\varepsilon$ by $\mathcal{O}(\varepsilon)$ and block $D_0^\varepsilon$ by $\mathcal{O}(\varepsilon^2)$. Lemmas~\ref{lem:class0-low-cross-blocks} and~\ref{lem:low-self-fourier} therefore give
\begin{align*}
R_0^\varepsilon&=
\begin{pmatrix}
-\rho_{\mathrm m}^\varepsilon+\mathcal{O}(\varepsilon)&\rho_{\mathrm c}^\varepsilon+\mathcal{O}(\varepsilon)\\
\mathcal{O}(\varepsilon)&\mathcal{O}(\varepsilon)
\end{pmatrix},\\
L_0^\varepsilon&=
\begin{pmatrix}
-\rho_{\mathrm m}^\varepsilon+\mathcal{O}(\varepsilon)&\rho_{\mathrm c}^\varepsilon+\mathcal{O}(\varepsilon)\\
\mu_{\mathrm m}^{-1}\rho_{\mathrm m}^\varepsilon+\mathcal{O}(\varepsilon)&-\mu^{-1}\rho_{\mathrm c}^\varepsilon+\mathcal{O}(\varepsilon)
\end{pmatrix},\\
D_0^\varepsilon&=
\begin{pmatrix}
-\dfrac12&\dfrac12-a(z)\\
-\dfrac1{2\mu_{\mathrm m}}&\dfrac{1}{\mu}\left(-\dfrac12+a(z)\right)
\end{pmatrix}
+\mathcal{O}(\varepsilon^2|\ln\varepsilon|).
\end{align*}
Here $\rho_{\mathrm c}^\varepsilon=\mathcal{O}(\varepsilon^{-1})$ and $\rho_{\mathrm m}^\varepsilon=\mathcal{O}(\varepsilon)$. Multiplying these matrices with $(Z_0^\varepsilon)^{-1}$ gives
\begin{align*}
L_0^\varepsilon(Z_0^\varepsilon)^{-1}R_0^\varepsilon &=
\begin{pmatrix}
0&(\rho_{\mathrm c}^\varepsilon)^2\bigl[(Z_0^\varepsilon)^{-1}\bigr]_{21}\\
0&-\mu^{-1}(\rho_{\mathrm c}^\varepsilon)^2\bigl[(Z_0^\varepsilon)^{-1}\bigr]_{21}
\end{pmatrix}
+\mathcal{O}(\varepsilon^2|\ln\varepsilon|)\\
&=
\begin{pmatrix}
0&-a(z)\\
0&\mu^{-1}a(z)
\end{pmatrix}
+\mathcal{O}(\varepsilon^2|\ln\varepsilon|).
\end{align*}
Another Schur elimination of \eqref{eq:four-dim-matrix} yields
$$
D_0^\varepsilon-L_0^\varepsilon(Z_0^\varepsilon)^{-1}R_0^\varepsilon =
\begin{pmatrix}
-\dfrac12&\dfrac12\\
-\dfrac1{2\mu_{\mathrm m}}&-\dfrac1{2\mu}
\end{pmatrix}
+\mathcal{O}(\varepsilon^2|\ln\varepsilon|),
$$
and hence
$$
\det \left(D_0^\varepsilon-L_0^\varepsilon(Z_0^\varepsilon)^{-1}R_0^\varepsilon\right) = \frac{\mu_{\mathrm m}+\mu}{4\mu\mu_{\mathrm m}} +\mathcal{O}(\varepsilon^2|\ln\varepsilon|).
$$
Since $\mu_{\mathrm m}\mu(\mu_{\mathrm m}+\mu)\ne0$, choose $r_{\mathrm c}$ small enough that this determinant is nonzero for sufficiently small $\varepsilon$ throughout $|z-z_{\mathrm c}|<r_{\mathrm c}$. For this choice of $r_{\mathrm c}$, the Schur matrix $D_0^\varepsilon-L_0^\varepsilon(Z_0^\varepsilon)^{-1}R_0^\varepsilon$ and hence the sector-$0$ restriction are invertible. Proposition~\ref{prop:conjugate-isospectral} gives the same conclusion in sector $1$. Thus candidate $\Delta_c$ does not yield a common characteristic frequency of sectors $0$ and $1$.

\subsubsection{The $\Delta_{\mathrm m}$ candidate}

Let
\begin{equation}\label{eq:z-star-exact}
z_*:=\frac{4\pi Q^2}{3|Y|}\frac{\mu_{\mathrm m}-\mu}{\mu_{\mathrm m}+\mu}.
\end{equation}
Choose $r_{\mathrm m}>0$ so that $\overline{B_{r_{\mathrm m}}(z_*)}$ avoids zero. In~\eqref{eq:delta-m-candidate-neighborhood}, write $\Delta_{\mathrm m}(\omega)=\varepsilon^2z$ with $z\in\overline{B_{r_{\mathrm m}}(z_*)}$. Lemma~\ref{lem:class0-low-cross-blocks} now gives
\begin{align*}
\mathsf B_{A,0;B,2}^\varepsilon &=
\begin{pmatrix}
\dfrac{a(z)}2+\mathcal{O}(\varepsilon^2)&\mathcal{O}(\varepsilon^2)\\
\mathcal{O}(\varepsilon^2)&\mathcal{O}(\varepsilon^2)
\end{pmatrix}, 
& \mathsf B_{B,-1;B,2}^\varepsilon&=\mathcal{O}(\varepsilon),\\
\mathsf B_{B,2;A,0}^\varepsilon &=
\begin{pmatrix}
\dfrac{a(z)}2+\mathcal{O}(\varepsilon^2)&\mathcal{O}(\varepsilon^2)\\
-\dfrac{a(z)}{\mu_{\mathrm m}}+\mathcal{O}(\varepsilon^2)&\mathcal{O}(\varepsilon^2)
\end{pmatrix}, 
& \mathsf B_{B,2;B,-1}^\varepsilon&=\mathcal{O}(\varepsilon).
\end{align*}
Combining these expressions with $E_*^{-1}$ and Lemma~\ref{lem:second-mode-self-block}, we find that eliminating the two $(B,2)$ coefficients gives
\begin{align*}
&
\begin{pmatrix}
\mathsf B_{A,0;B,2}^{\varepsilon}\\
\mathsf B_{B,-1;B,2}^{\varepsilon}
\end{pmatrix}
(E_0^\varepsilon)^{-1}
\begin{pmatrix}
\mathsf B_{B,2;A,0}^{\varepsilon} & \mathsf B_{B,2;B,-1}^{\varepsilon}
\end{pmatrix}
\\
&\qquad =
\begin{pmatrix}
a(z)^2\dfrac{\mu-\mu_{\mathrm m}} {\mu+\mu_{\mathrm m}}&0&0&0\\
0&0&0&0\\
0&0&0&0\\
0&0&0&0
\end{pmatrix}
+
\begin{pmatrix}
\mathcal{O}(\varepsilon^2|\ln\varepsilon|) &\mathcal{O}(\varepsilon^2|\ln\varepsilon|) &\mathcal{O}(\varepsilon)&\mathcal{O}(\varepsilon)\\
\mathcal{O}(\varepsilon^2|\ln\varepsilon|) &\mathcal{O}(\varepsilon^2|\ln\varepsilon|) &\mathcal{O}(\varepsilon)&\mathcal{O}(\varepsilon)\\
\mathcal{O}(\varepsilon)&\mathcal{O}(\varepsilon) &\mathcal{O}(\varepsilon^2)&\mathcal{O}(\varepsilon^2)\\
\mathcal{O}(\varepsilon)&\mathcal{O}(\varepsilon) &\mathcal{O}(\varepsilon^2)&\mathcal{O}(\varepsilon^2)
\end{pmatrix}.
\end{align*}
Thus only the $(1,1)$ entry acquires a nonzero $\mathcal{O}(1)$ leading correction. Denote the upper-left $2\times 2$ block by $\widetilde Z_j^\varepsilon$.

\begin{lemma}\label{lem:lem-5-9}
For fixed nonzero $\lambda$,
$$
(\mathcal{K}_{H_0,\lambda}^{pp,\varepsilon,*})_{0,0} =\frac12-\frac{\lambda^2\varepsilon^2}{2} \bigl(\ln\varepsilon+\ln\lambda+\gamma_0\bigr) +\mathcal{O}(\varepsilon^4|\ln\varepsilon|).
$$
\end{lemma}

\begin{proof}
For $r>\varepsilon$, Graf's addition theorem gives
$$
H_0^{(1)}(\lambda|\bm x-\bm y|) = \sum_{n\in\mathbb Z} H_n^{(1)}(\lambda r) J_n(\lambda\varepsilon)e^{\ii n(\theta-t)},
$$
where $\bm x=r(\cos\theta,\sin\theta)$ and $\bm y=\varepsilon(\cos t,\sin t)$. Taking the exterior normal derivative of the single-layer potential and applying the jump relation yields
$$
(\mathcal{K}_{H_0,\lambda}^{pp,\varepsilon,*})_{0,0} = -\frac{\ii\pi\lambda\varepsilon}{2} J_0(\lambda\varepsilon) H_0^{(1)\prime}(\lambda\varepsilon) -\frac12.
$$
Set $x=\lambda\varepsilon$. Using $H_0^{(1)\prime}(x)=-H_1^{(1)}(x)$, $J_0(x)=1-x^2/4+\mathcal{O}(x^4)$ and $Y_1(x)=-2/(\pi x)+x[\ln(x/2)+\gamma-1/2]/\pi+\mathcal{O}(x^3|\ln x|)$ gives
$$
-\frac{\ii\pi x}{2}J_0(x)H_0^{(1)\prime}(x)-\frac12 =\frac12-\frac{x^2}{2} \left(\ln x-\ln2+\gamma-\frac{\ii\pi}{2}\right) +\mathcal{O}(x^4|\ln x|),
$$
which proves the formula.
\end{proof}

Lemma~\ref{lem:low-self-fourier} and Lemma~\ref{lem:lem-5-9} give
\begin{align*}
\widetilde Z_j^\varepsilon &=
\begin{pmatrix}
(\widehat{\mathcal{S}}_{\omega_{\mathrm m}}^{pp,\varepsilon})_{0,0}-a(z)^2\dfrac{\mu-\mu_{\mathrm m}}{\mu+\mu_{\mathrm m}}&-(\widehat{\mathcal{S}}_{\omega_{\mathrm c}}^{pp,\varepsilon})_{0,0}\\ 
-\mu_{\mathrm m}^{-1} \left(\dfrac12+(\mathcal{K}_{\omega_{\mathrm m}}^{*,pp,\varepsilon})_{0,0}\right)& \mu^{-1} \left(-\dfrac12+(\mathcal{K}_{\omega_{\mathrm c}}^{*,pp,\varepsilon})_{0,0}\right)
\end{pmatrix}
+
\begin{pmatrix}
\mathcal{O}(\varepsilon^2|\ln\varepsilon|)&\mathcal{O}(\varepsilon^2)\\
\mathcal{O}(\varepsilon^2)&\mathcal{O}(\varepsilon^4)
\end{pmatrix}
\\
&=
\begin{pmatrix}
\dfrac{6\pi}{|Y|\varepsilon^2z}+\mathcal{O}(|\ln\varepsilon|+1)&-\ln\varepsilon+\mathcal{O}(1)\\ 
-\mu_{\mathrm m}^{-1}(1-2a(z))+\mathcal{O}(\varepsilon^2|\ln\varepsilon|)& -\dfrac{\omega_{\mathrm c}(\omega)^2}{2\mu}\varepsilon^2\ln\varepsilon+\mathcal{O}(\varepsilon^2)
\end{pmatrix}.
\end{align*}
Therefore
\begin{align*}
\det\widetilde Z_j^\varepsilon &=-\frac{3\pi\omega_{\mathrm c}(\omega)^2}{\mu|Y|z}\ln\varepsilon -\frac{1-2a(z)}{\mu_{\mathrm m}}\ln\varepsilon+\mathcal{O}(1)\\
&=-\frac{2a(z)}{\mu_{\mathrm m}}\ln\varepsilon -\frac{1-2a(z)}{\mu_{\mathrm m}}\ln\varepsilon+\mathcal{O}(1) =-\frac1{\mu_{\mathrm m}}\ln\varepsilon+\mathcal{O}(1).
\end{align*}
The block $\widetilde Z_j^\varepsilon$ is hence invertible, and its inverse is
\begin{equation}\label{eq:Ztilde-inverse-orders}
(\widetilde Z_j^\varepsilon)^{-1} =
\begin{pmatrix}
\dfrac{k_D^2}{2}\varepsilon^2+\mathcal{O} \left(\dfrac{\varepsilon^2}{|\ln\varepsilon|}\right)& -\mu_{\mathrm m}+\mathcal{O} \left(\dfrac1{|\ln\varepsilon|}\right)\\
\mathcal{O} \left(\dfrac1{|\ln\varepsilon|}\right)& \mathcal{O} \left(\dfrac{\varepsilon^{-2}}{|\ln\varepsilon|}\right)
\end{pmatrix}, \qquad  j=0,1.
\end{equation}

We next compute $R_0^\varepsilon$ and $L_0^\varepsilon$. Besides the terms in Lemma~\ref{lem:class0-low-cross-blocks}, we also need to retain the leading order term of $(\mathcal{K}_{\Lambda_0^*,\omega_{\mathrm m}}^{AB,\varepsilon,*})_{0,-1}$. Equation~\eqref{eq:exact-res-K}, $J_0'=-J_1$, $J_{-1}=-J_1$, and $\Xi_{AB}(-1)=-3\ii$ give
$$
(\mathcal{K}_{\Lambda_0^*,\omega_{\mathrm m}}^{AB,\varepsilon,*})_{0,-1} =\frac{6\pi\varepsilon k_D}{|Y|\Delta_{\mathrm m}}J_1(z_D)^2 =\frac{3\pi k_D^3\varepsilon^3}{2|Y|\Delta_{\mathrm m}}+\mathcal{O}(\varepsilon^3) =\frac{k_D^2\varepsilon^2}{2}\rho_{\mathrm m}^\varepsilon+\mathcal{O}(\varepsilon^3).
$$
Then, after one Schur elimination, the $R_0^\varepsilon$ and $L_0^\varepsilon$ in \eqref{eq:four-dim-matrix} are
$$
R_0^\varepsilon=
\begin{pmatrix}
-\rho_{\mathrm m}^\varepsilon+\mathcal{O}(\varepsilon)&\mathcal{O}(\varepsilon)\\
-\dfrac{k_D^2\varepsilon^2}{2\mu_{\mathrm m}}\rho_{\mathrm m}^\varepsilon+\mathcal{O}(\varepsilon^2)&\mathcal{O}(\varepsilon^2)
\end{pmatrix},  \qquad L_0^\varepsilon=
\begin{pmatrix}
-\rho_{\mathrm m}^\varepsilon+\mathcal{O}(\varepsilon)&\mathcal{O}(\varepsilon)\\
\mu_{\mathrm m}^{-1}\rho_{\mathrm m}^\varepsilon+\mathcal{O}(\varepsilon)&\mathcal{O}(\varepsilon)
\end{pmatrix},
$$
where $\rho_{\mathrm m}^\varepsilon=\mathcal{O}(\varepsilon^{-1})$. By~\eqref{eq:Ztilde-inverse-orders}, the first entry of the product with $R_0^\varepsilon$ satisfies
\begin{align*}
\bigl[(\widetilde Z_0^\varepsilon)^{-1}R_0^\varepsilon\bigr]_{11} &=\left(\frac{k_D^2}{2}\varepsilon^2+\mathcal{O} \left(\frac{\varepsilon^2}{|\ln\varepsilon|}\right)\right) \bigl(-\rho_{\mathrm m}^\varepsilon+\mathcal{O}(\varepsilon)\bigr)\\
&\quad + \left(-\mu_{\mathrm m}+\mathcal{O} \left(\frac1{|\ln\varepsilon|}\right)\right) \left(-\frac{k_D^2\varepsilon^2}{2\mu_{\mathrm m}}\rho_{\mathrm m}^\varepsilon+\mathcal{O}(\varepsilon^2)\right)\\
&=\mathcal{O} \left(\frac{\varepsilon}{|\ln\varepsilon|}+\varepsilon^2\right).
\end{align*}
The remaining three entries of $(\widetilde Z_0^\varepsilon)^{-1}R_0^\varepsilon$ are calculated similarly. Combining these calculations gives
$$
(\widetilde Z_0^\varepsilon)^{-1}R_0^\varepsilon =
\begin{pmatrix}
\mathcal{O} \left(\dfrac{\varepsilon}{|\ln\varepsilon|}+\varepsilon^2\right)&\mathcal{O}(\varepsilon^2)\\
\mathcal{O} \left(\dfrac{\varepsilon^{-1}}{|\ln\varepsilon|}\right)&\mathcal{O} \left(\dfrac1{|\ln\varepsilon|}\right)
\end{pmatrix}.
$$
Multiplication on the left by $L_0^\varepsilon$ gives
\begin{equation}\label{eq:constant-elimination-exact-order}
\left\|L_0^\varepsilon(\widetilde Z_0^\varepsilon)^{-1}R_0^\varepsilon\right\| \le \frac{C}{|\ln\varepsilon|}.
\end{equation}
By the complex intertwining identity~\eqref{eq:J-complex-intertwining}, the estimate~\eqref{eq:constant-elimination-exact-order} also holds for $\|L_1^\varepsilon(\widetilde{Z}_1^\varepsilon)^{-1}R_1^\varepsilon\|$ throughout the conjugation-invariant frequency neighborhood. On its real part, this also follows from Proposition~\ref{prop:conjugate-isospectral}.

Finally, Lemma~\ref{lem:low-self-fourier} and one Schur elimination give
$$
D_j^\varepsilon=
\begin{pmatrix}
-\dfrac12+a(z)&\dfrac12\\
-\mu_{\mathrm m}^{-1} \left(\dfrac12+a(z)\right)&-\dfrac1{2\mu}
\end{pmatrix}
+\mathcal{O}(\varepsilon^2|\ln\varepsilon|), \qquad  j=0,1.
$$
As in Section~\ref{sec:delta-c-candidate}, a second Schur elimination gives $F_j^\varepsilon(z):=\det\bigl(D_j^\varepsilon-L_j^\varepsilon(\widetilde Z_j^\varepsilon)^{-1}R_j^\varepsilon\bigr)$, which satisfies
\begin{equation}\label{eq:exact-scalar-characteristic-expansion}
F_j^\varepsilon(z)=F_*(z)+\mathcal{O} \left(\frac1{|\ln\varepsilon|}\right), \qquad  j=0,1,
\end{equation}
uniformly on $\overline{B_{r_{\mathrm m}}(z_*)}$, where
\begin{equation}\label{eq:F-star-z}
F_*(z)= \frac{\mu_{\mathrm m}+\mu}{4\mu\mu_{\mathrm m}} +\frac{\mu-\mu_{\mathrm m}}{2\mu\mu_{\mathrm m}} \frac{2\pi Q^2}{3|Y|z}.
\end{equation}

\subsection{Calculation of eigenvalues}

We now convert the effective equation near $z_*$ into an existence and multiplicity statement for the full Bloch problem. Define $\omega_{0,m}:=\sqrt{|\bm\alpha^*|^2/(\epsilon\mu_m)}$, which satisfies $\omega_m(\omega_{0,m}) = |\bm\alpha^*|$.

\begin{theorem}\label{thm:exact-dirac-frequency}
Assume the fixed positive material parameters of Section~\ref{sec:sec2-1}, $\kappa_g=0$, and~\eqref{eq:boundary-material-nonresonance}. For sufficiently small $\varepsilon$, each of the restrictions of $\mathcal A^\varepsilon(\bm\alpha^*,\omega)$ to sectors $0$ and $1$ has the same characteristic frequency $\omega_*^\varepsilon$. Furthermore,
\begin{align}
\Delta_{\mathrm m}(\omega_*^\varepsilon) &=z_*\varepsilon^2+ \mathcal{O} \left(\frac{\varepsilon^2}{|\ln\varepsilon|}\right), \label{eq:exact-Delta-m-root}\\
\omega_*^\varepsilon &=\omega_{0,m} +\frac{2\pi Q^2} {3|Y||\bm\alpha^*|\omega_{\mathrm m}'(\omega_{0,m})} \frac{\mu_{\mathrm m}-\mu}{\mu_{\mathrm m}+\mu}\varepsilon^2 +\mathcal{O} \left(\frac{\varepsilon^2}{|\ln\varepsilon|}\right). \label{eq:exact-omega-m-root}
\end{align}
The kernel of each restriction of $\mathcal A^\varepsilon(\bm\alpha^*,\omega_*^\varepsilon)$ is one-dimensional. After decreasing $r_{\mathrm m}$ if necessary, the restriction of $\mathcal A^\varepsilon(\bm\alpha^*,\omega)$ to sector $2$ is invertible throughout the same frequency neighborhood. Consequently,
\begin{equation}\label{eq:exact-two-dimensional-kernel}
\dim\ker\mathcal A^\varepsilon (\bm\alpha^*,\omega_*^\varepsilon)=2.
\end{equation}
\end{theorem}

\begin{proof}
By~\eqref{eq:F-star-z} and~\eqref{eq:z-star-exact}, $F_*(z_*)=0$ and
$$
F_*'(z_*)= -\frac{\mu-\mu_{\mathrm m}}{2\mu\mu_{\mathrm m}}  \frac{2\pi Q^2}{3|Y|z_*^2}\ne0.
$$
Choose $r_{\mathrm m}$ sufficiently small that $z_*$ is the only zero of $F_*$ in $\overline{B_{r_{\mathrm m}}(z_*)}$ and set
$$
c_*:=\min_{|z-z_*|=r_{\mathrm m}}|F_*(z)|>0.
$$
By Proposition~\ref{prop:six-dimensional-reduction} and Lemma~\ref{lem:second-mode-self-block}, the operators $Q_j\mathbb A^\varepsilon Q_j$ and $E_j^\varepsilon$ are uniformly invertible for $j=0,1$. Moreover, the calculation in Section~\ref{sec:sec4-3} gives
$$
\det\widetilde Z_j^\varepsilon = -\frac{1}{\mu_{\mathrm m}}\ln\varepsilon+\mathcal{O}(1),
$$
so $\widetilde Z_j^\varepsilon$ is also invertible for sufficiently small $\varepsilon$. These conclusions hold uniformly on a slightly larger complex neighborhood of $\overline{B_{r_{\mathrm m}}(z_*)}$, after reducing $r_{\mathrm m}$ if necessary. Since the operator families are analytic there, their inverses and the resulting functions $F_j^\varepsilon(z)$ are analytic on that neighborhood. Equation~\eqref{eq:exact-scalar-characteristic-expansion} gives, for small $\varepsilon$,
$$
\sup_{|z-z_*|=r_{\mathrm m}} |F_j^\varepsilon(z)-F_*(z)|<c_*.
$$
Rouch\'e's theorem implies that $F_j^\varepsilon$ has exactly one zero $z_j^\varepsilon$ in $B_{r_{\mathrm m}}(z_*)$, counted with algebraic multiplicity. This zero is therefore simple. Since
$$
F_*(z)=\frac{\mu_{\mathrm m}+\mu}{4\mu\mu_{\mathrm m}}  \frac{z-z_*}{z},
$$
the estimate in~\eqref{eq:exact-scalar-characteristic-expansion} also gives
$$
z_j^\varepsilon=z_*+\mathcal{O}(|\ln\varepsilon|^{-1}),\qquad  j=0,1.
$$
Using the relation $\Delta_m(\omega)=\varepsilon^2z$, we obtain the corresponding characteristic frequency $\omega_j^\varepsilon:=\sqrt{(|\bm\alpha^*|^2+\varepsilon^2z_j^\varepsilon)/(\epsilon\mu_m)}$. The Schur eliminations show that $\ker\mathcal A_j^\varepsilon(\omega_j^\varepsilon)\ne\{0\}$. The invertibility of single-layer potential and Proposition~\ref{prop:pde-bie-equivalence} give a nonzero Bloch eigenfield $u$ for~\eqref{eq:system}. Green's identity gives
$$
\frac1{\mu_{\mathrm m}} \int_{Y\setminus\overline{D^\varepsilon}}|\nabla u|^2 d\bm x +\frac1\mu\int_{D^\varepsilon}|\nabla u|^2 d\bm x =\epsilon(\omega_j^\varepsilon)^2 \int_Y|u|^2 d\bm x.
$$
It follows that $(\omega_j^\varepsilon)^2$ is real. Since $\omega_j^\varepsilon\to\omega_{0,m}>0$, the frequencies are real and positive for sufficiently small $\varepsilon$, and hence $z_j^\varepsilon\in\mathbb R$.

From Proposition~\ref{prop:conjugate-isospectral}, $\omega_0^\varepsilon$ is also a characteristic frequency of sector $1$. The uniqueness of the zero of $F_1^\varepsilon$ in the same disk implies $z_0^\varepsilon=z_1^\varepsilon$. Denote the common zero by $z^\varepsilon$ and set $\omega_*^\varepsilon:=\sqrt{(|\bm\alpha^*|^2+\varepsilon^2z^\varepsilon)/(\epsilon\mu_m)}$. This proves~\eqref{eq:exact-Delta-m-root}. Since
$$
\Delta_{\mathrm m}'(\omega_{0,m}) =2|\bm\alpha^*|\omega_{\mathrm m}'(\omega_{0,m})\ne0,
$$
the Taylor expansion,
$$
\Delta_{\mathrm m}(\omega) =2|\bm\alpha^*|\omega_{\mathrm m}'(\omega_{0,m}) (\omega-\omega_{0,m}) +\mathcal{O}((\omega-\omega_{0,m})^2)
$$
gives~\eqref{eq:exact-omega-m-root}. The final $2\times2$ Schur matrix, after repeated Schur eliminations, has a determinant with a simple zero. Its rank at that zero is one, thus each of sectors $0$ and $1$ has a one-dimensional kernel.

For sector $2$, the Fourier indices on $A$ and $B$ are respectively $n=3k-1$ and $n=3k+1$, so the lowest modes are $(A,-1)$ and $(B,1)$. Continue to write $a(z)=2\pi Q^2/(3|Y|z)$. By Lemma~\ref{lem:low-self-fourier}, both diagonal two-density blocks equal $\mathsf D_0(z)+\mathcal{O}(\varepsilon^2|\ln\varepsilon|)$, where
$$
\mathsf D_0(z)=
\begin{pmatrix}
-\dfrac12+a(z)&\dfrac12\\
-\mu_{\mathrm m}^{-1} \left(\dfrac12+a(z)\right)&-\dfrac1{2\mu}
\end{pmatrix}.
$$
The phase sums $\Xi_{AB}(2)=\Xi_{BA}(-2)=-3$ and \eqref{eq:exact-res-S}--\eqref{eq:exact-res-K} give
\begin{align*}
(\widehat{\mathcal{S}}_{\omega_{\mathrm m}}^{AB,\varepsilon})_{-1,1} &=(\widehat{\mathcal{S}}_{\omega_{\mathrm m}}^{BA,\varepsilon})_{1,-1} =-a(z)+\mathcal{O}(\varepsilon^2),\\
(\mathcal{K}_{\omega_{\mathrm m}}^{AB,\varepsilon,*})_{-1,1} &=(\mathcal{K}_{\omega_{\mathrm m}}^{BA,\varepsilon,*})_{1,-1} =-a(z)+\mathcal{O}(\varepsilon^2).
\end{align*}
The corresponding entries for frequency $\omega_c$ are all $\mathcal{O}(\varepsilon^2)$. Consequently the two off-diagonal blocks satisfy
$$
\mathsf B_{A,-1;B,1}^\varepsilon =\mathsf C(z)+\mathcal{O}(\varepsilon^2),\qquad  \mathsf B_{B,1;A,-1}^\varepsilon =\mathsf C(z)+\mathcal{O}(\varepsilon^2),\qquad  \mathsf C(z):=
\begin{pmatrix}
-a(z)&0\\
\mu_{\mathrm m}^{-1}a(z)&0
\end{pmatrix}.
$$
In the coefficient order $(c_{1,A,-1},c_{2,A,-1},c_{1,B,1},c_{2,B,1})^\top$, the leading matrix is therefore
$$
\begin{pmatrix}
\mathsf D_0(z)&\mathsf C(z)\\
\mathsf C(z)&\mathsf D_0(z)
\end{pmatrix}.
$$
Write the two pairs of density coefficients as $X_A,X_B\in\mathbb C^2$ and set $X_\pm=X_A\pm X_B$. Adding and subtracting the two equations reduces this matrix to $\mathsf D_+(z)X_+=0$ and $\mathsf D_-(z)X_-=0$, where
$$
\mathsf D_+(z)=
\begin{pmatrix}
-\dfrac12&\dfrac12\\
-\dfrac1{2\mu_{\mathrm m}}&-\dfrac1{2\mu}
\end{pmatrix},  \qquad \mathsf D_-(z)=
\begin{pmatrix}
-\dfrac12+2a(z)&\dfrac12\\
-\mu_{\mathrm m}^{-1} \left(\dfrac12+2a(z)\right)&-\dfrac1{2\mu}
\end{pmatrix}.
$$
Their determinants are
$$
\det\mathsf D_+ =\frac{\mu_{\mathrm m}+\mu}{4\mu\mu_{\mathrm m}},\qquad  \det\mathsf D_-(z) =\frac{4a(z)(\mu-\mu_{\mathrm m})+\mu_{\mathrm m}+\mu} {4\mu\mu_{\mathrm m}}.
$$
By~\eqref{eq:z-star-exact}, $a(z_*)=(\mu_{\mathrm m}+\mu)/(2(\mu_{\mathrm m}-\mu))$, and hence
$$
\det\mathsf D_-(z_*)= -\frac{\mu_{\mathrm m}+\mu}{4\mu\mu_{\mathrm m}}\ne0,\qquad  \det\mathsf D_+= \frac{\mu_{\mathrm m}+\mu}{4\mu\mu_{\mathrm m}}\ne0.
$$
After decreasing $r_{\mathrm m}$, both matrices are uniformly invertible for $|z-z_*|\le r_{\mathrm m}$.

It remains to justify elimination of modes with $|n|\geq 2$ in this sector. With Proposition~\ref{prop:six-dimensional-reduction} and similar Schur eliminations, the sector-2 matrix reduces to
$$
\begin{pmatrix}
\mathsf D_0(z)&\mathsf C(z)\\
\mathsf C(z)&\mathsf D_0(z)
\end{pmatrix}
+\mathcal{O}(\varepsilon^2|\ln\varepsilon|),
$$
and is invertible for sufficiently small $\varepsilon$ throughout the disk. This proves the asserted invertibility of $\mathcal A_2^\varepsilon(\omega)$. Combining it with the one-dimensional kernels in sectors $0$ and $1$ and~\eqref{eq:kernel-direct-sum} proves \eqref{eq:exact-two-dimensional-kernel}.
\end{proof}

\subsection{Existence of Dirac cone}

The double frequency established above does not yet determine the dispersion away from $\bm\alpha^*$. We therefore analyze the first-order Bloch perturbation of the full boundary operator on its two-dimensional kernel. According to Theorem \ref{thm:exact-dirac-frequency}, we set
$$
\operatorname{span}\{C^{(0)},C^{(1)}\}:=\ker\mathbb A^\varepsilon(\bm\alpha^*,\omega_*^\varepsilon),
$$
where $C^{(j)}$ belongs to sector $j$, $j=0,1$. For brevity we write $\mathbb A_*^\varepsilon:=\mathbb A^\varepsilon(\bm\alpha^*,\omega_*^\varepsilon)$.

We now prove that the operator $\mathbb A_*^\varepsilon$ is Fredholm of index zero. By Proposition~\ref{prop:six-dimensional-reduction}, for $j=0,1$ the complementary operator
$$
Q_j\mathbb A_*^\varepsilon Q_j: Q_j\mathfrak H\longrightarrow Q_j\mathfrak H
$$
is boundedly invertible. Consequently, the Schur reduction in Proposition \ref{prop:six-dimensional-reduction} shows that $\mathbb A_*^\varepsilon|_{\mathfrak H^{(j)}}$ is equivalent, through boundedly invertible block operations, to the direct sum of $Q_j\mathbb A_*^\varepsilon Q_j$ and the six-dimensional Schur matrix
$$
P_j\mathbb A_*^\varepsilon P_j - P_j\mathbb A_*^\varepsilon Q_j (Q_j\mathbb A_*^\varepsilon Q_j)^{-1} Q_j\mathbb A_*^\varepsilon P_j.
$$
The first operator is invertible, while the second is a square finite-dimensional matrix. Both are Fredholm of index zero. Hence $\mathbb A_*^\varepsilon|_{\mathfrak H^{(j)}}$ is Fredholm of index zero for $j=0,1$. By Theorem~\ref{thm:exact-dirac-frequency}, the restriction to sector $2$ is invertible. The direct-sum decomposition into the three symmetry sectors therefore gives $\operatorname{ind}\mathbb A_*^\varepsilon=0$. It follows that $\dim\ker(\mathbb A_*^\varepsilon)^*=2$, with one vector $\widetilde C^{(j)}\in\ker(\mathbb A_*^\varepsilon)^*$ in each of sectors $0$ and $1$.

\begin{lemma}
For $j=0,1$,
\begin{equation}\label{eq:lem-proof-1}
\left\langle\widetilde C^{(j)}, \partial_\omega\mathbb A_*^\varepsilon C^{(j)} \right\rangle_{\mathfrak H}\ne0.
\end{equation}
For $j\ne\ell$ the corresponding pairing is zero. Thus $\widetilde C^{(j)}$ may be rescaled so that
\begin{equation}\label{eq:boundary-frequency-normalization}
\left\langle\widetilde C^{(j)}, \partial_\omega\mathbb A_*^\varepsilon C^{(\ell)} \right\rangle_{\mathfrak H}=\delta_{j\ell}, \qquad  j,\ell=0,1.
\end{equation}
\end{lemma}
\begin{proof}
From the definition of $\mathcal R$ in \eqref{eq:coefficient-rotation-action} and Theorem \ref{thm:mod3-invariant-subspaces}, $\mathcal R$ preserves sectors and commutes with $\mathbb A^\varepsilon$, so the derivative $\partial_\omega\mathbb A_*^\varepsilon$ preserves the orthogonal sector decomposition, which proves the vanishing of the pairings for $j\ne\ell$. If there exists $j$ such that
$$
\left\langle\widetilde C^{(j)}, \partial_\omega\mathbb A_*^\varepsilon C^{(j)} \right\rangle_{\mathfrak H}=0,
$$
then, since $\partial_\omega\mathbb A_*^\varepsilon\mathfrak H^{(j)}\subset\mathfrak H^{(j)}$, it follows that for $\ell\neq j$,
$$
\left\langle\widetilde C^{(\ell)}, \partial_\omega\mathbb A_*^\varepsilon C^{(j)} \right\rangle_{\mathfrak H}=0,
$$
hence $\partial_\omega\mathbb A_*^\varepsilon C^{(j)}\in\bigl(\ker(\mathbb A_*^\varepsilon)^*\bigr)^\perp=\operatorname{Ran}\mathbb A_*^\varepsilon$. There then exists $W^{(j)}\in\mathfrak H^{(j)}$ satisfying $\mathbb A_*^\varepsilon W^{(j)}=-\partial_\omega\mathbb A_*^\varepsilon C^{(j)}$. Consequently,
$$
\mathbb A^\varepsilon(\bm\alpha^*,\omega_*^\varepsilon+t)(C^{(j)}+tW^{(j)})=\mathcal{O}(t^2).
$$
The associated scalar characteristic function therefore vanishes to order at least two. It differs from $F_j^\varepsilon$, defined before~\eqref{eq:exact-scalar-characteristic-expansion}, only by a nonvanishing analytic factor. Since $\Delta_{\mathrm m}'(\omega_*^\varepsilon)\ne0$, this contradicts the simple zero proved in Theorem~\ref{thm:exact-dirac-frequency}. This proves \eqref{eq:lem-proof-1}, and rescaling proves \eqref{eq:boundary-frequency-normalization}.
\end{proof}

We next consider the first derivative with respect to the Bloch vector. Let $\bm\alpha(t):=\bm\alpha^*+t\delta\bm\alpha$. Since $R^\top\bm\alpha^*=\bm\alpha^*+\bm q_R$ and $R^\top=R^{-1}$,
\begin{equation}\label{eq:rotation-alpha-variance}
R^\top\bm\alpha(t) =\bm\alpha^*+\bm q_R+tR^{-1}\delta\bm\alpha.
\end{equation}
The rotated Green-function identity, together with the coefficient rotation~\eqref{eq:coefficient-rotation-action}, gives the corresponding operator covariance identity. Differentiating that identity at $t=0$ gives the following lemma.
\begin{lemma}\label{lem:lem4-5}
For $j,\ell=0,1$,
\begin{equation}\label{eq:boundary-bloch-rotation-pairing}
\left\langle\widetilde C^{(j)}, (\delta\bm\alpha\cdot\nabla_{\bm\alpha}) \mathbb A_*^\varepsilon C^{(\ell)} \right\rangle_{\mathfrak H} =e^{2\pi \ii(\ell-j)/3} \left\langle\widetilde C^{(j)}, (R^{-1}\delta\bm\alpha\cdot\nabla_{\bm\alpha}) \mathbb A_*^\varepsilon C^{(\ell)} \right\rangle_{\mathfrak H}.
\end{equation}
\end{lemma}
\begin{proof}
For any Bloch vector $\bm\alpha$, we have $G_\lambda^{\bm\alpha}(R\bm x)=G_\lambda^{R^\top\bm\alpha}(\bm x)$. Recall that $R\bm x_p(\theta)=\bm x_p(\theta+\vartheta)+\bm\ell_p$, where $\bm\ell_A=\bm 0$, $\bm\ell_B=-\bm e_2$. Consequently,
$$
G_\lambda^{R^\top\bm\alpha} \bigl(\bm x_p(\theta)-\bm x_q(t)\bigr)= G_\lambda^{\bm\alpha} \bigl(R\bm x_p(\theta)-R\bm x_q(t)\bigr) =e^{\ii\bm\alpha\cdot(\bm\ell_p-\bm\ell_q)} G_\lambda^{\bm\alpha} \bigl( \bm x_p(\theta+\vartheta) -\bm x_q(t+\vartheta) \bigr).
$$

Thus every Fourier matrix entry of $\mathbb A^\varepsilon (R^\top\bm\alpha,\omega)$ satisfies
$$
\bigl[ \mathbb A^\varepsilon (R^\top\bm\alpha,\omega) \bigr]_{\sigma,p,m;\tau,q,n} = e^{\ii\bm\alpha\cdot(\bm\ell_p-\bm\ell_q)} e^{2\pi \ii(m-n)/3} \bigl[ \mathbb A^\varepsilon(\bm\alpha,\omega) \bigr]_{\sigma,p,m;\tau,q,n}.
$$
Here $\sigma,\tau\in\{1,2\}$ are the two-density indices. Extend the coefficient rotation to $\mathcal R_{\bm\alpha}$ by $(\mathcal R_{\bm\alpha}C)_{\sigma,p,n}:=e^{\ii\bm\alpha\cdot\bm\ell_p}e^{2\pi \ii n/3} C_{\sigma,p,n}$. At $\bm\alpha=\bm\alpha^*$, it equals $\mathcal R$. Then
\begin{align*}
\bigl[ \mathbb A^\varepsilon (R^\top\bm\alpha,\omega) \mathcal R_{\bm\alpha} \bigr]_{\sigma,p,m;\tau,q,n} &= e^{\ii\bm\alpha\cdot\bm\ell_q} e^{2\pi \ii n/3} \bigl[ \mathbb A^\varepsilon (R^\top\bm\alpha,\omega) \bigr]_{\sigma,p,m;\tau,q,n} \\
&= e^{\ii\bm\alpha\cdot\bm\ell_p} e^{2\pi \ii m/3} \bigl[ \mathbb A^\varepsilon(\bm\alpha,\omega) \bigr]_{\sigma,p,m;\tau,q,n},
\end{align*}
which implies
$$
\mathbb A^\varepsilon (R^\top\bm\alpha,\omega) \mathcal R_{\bm\alpha} = \mathcal R_{\bm\alpha} \mathbb A^\varepsilon(\bm\alpha,\omega).
$$
From \eqref{eq:rotation-alpha-variance} and periodicity of $\mathbb{A}^\varepsilon$ in the Bloch vector modulo $\Lambda^*$, we have
\begin{equation}\label{eq:nearby-bloch-rotation-covariance}
\mathbb A^\varepsilon \bigl( \bm\alpha^*+tR^{-1}\delta\bm\alpha, \omega_*^\varepsilon \bigr) \mathcal R_{\bm\alpha(t)} = \mathcal R_{\bm\alpha(t)} \mathbb A^\varepsilon \bigl( \bm\alpha^*+t\delta\bm\alpha, \omega_*^\varepsilon \bigr).
\end{equation}
For brevity, in this calculation, let
$$
\dot{\mathcal R} := \left. \frac{d}{dt}\mathcal R_{\bm\alpha(t)} \right|_{t=0}.
$$
Its Fourier coefficient action is explicitly
$$
(\dot{\mathcal R}C)_{\sigma,p,n} = \left. \frac{d}{dt} \left( e^{\ii(\bm\alpha^*+t\delta\bm\alpha) \cdot\bm\ell_p} e^{2\pi \ii n/3} C_{\sigma,p,n} \right) \right|_{t=0} = \ii(\delta\bm\alpha\cdot\bm\ell_p) e^{2\pi \ii(n+s_p)/3}C_{\sigma,p,n}.
$$
In particular, $(\dot{\mathcal R}C)_{\sigma,A,n}=0$ and $(\dot{\mathcal R}C)_{\sigma,B,n}=-\ii(\delta\bm\alpha\cdot\bm e_2)e^{2\pi \ii(n+1)/3}C_{\sigma,B,n}$.

Differentiate \eqref{eq:nearby-bloch-rotation-covariance} with respect to $t$ at $t=0$. The product rule yields
\begin{equation}\label{eq:differentiated-bloch-rotation-covariance}
(R^{-1}\delta\bm\alpha\cdot\nabla_{\bm\alpha}) \mathbb A_*^\varepsilon \mathcal R+ \mathbb A_*^\varepsilon\dot{\mathcal R}= \dot{\mathcal R}\mathbb A_*^\varepsilon+ \mathcal R (\delta\bm\alpha\cdot\nabla_{\bm\alpha}) \mathbb A_*^\varepsilon.
\end{equation}

Apply \eqref{eq:differentiated-bloch-rotation-covariance} to $C^{(\ell)}$, and take the $\mathfrak H$ inner product with $\widetilde C^{(j)}$ on both sides. Since $\mathbb A_*^\varepsilon C^{(\ell)}=0$ and $(\mathbb A_*^\varepsilon)^* \widetilde C^{(j)}=0$, the equation reduces to
$$
\left\langle \widetilde C^{(j)}, (R^{-1}\delta\bm\alpha\cdot\nabla_{\bm\alpha}) \mathbb A_*^\varepsilon \mathcal R C^{(\ell)} \right\rangle_{\mathfrak H}= \left\langle \widetilde C^{(j)}, \mathcal R (\delta\bm\alpha\cdot\nabla_{\bm\alpha}) \mathbb A_*^\varepsilon C^{(\ell)} \right\rangle_{\mathfrak H}.
$$

Since $\mathcal R C^{(\ell)}=e^{2\pi \ii\ell/3} C^{(\ell)}$ and $\mathcal R^*\widetilde C^{(j)}=e^{-2\pi \ii j/3}\widetilde C^{(j)}$, we prove the desired result.
\end{proof}

With the help of Lemma \ref{lem:lem4-5}, the calculation of $\bm\alpha$-derivative can be significantly simplified as the following lemma demonstrates.

\begin{lemma}
For $j=\ell$, the pairing \eqref{eq:boundary-bloch-rotation-pairing} satisfies
\begin{equation}\label{eq:boundary-bloch-diagonal-zero}
\left\langle\widetilde C^{(j)}, (\delta\bm\alpha\cdot\nabla_{\bm\alpha}) \mathbb A^\varepsilon(\bm\alpha^*,\omega_*^\varepsilon) C^{(j)} \right\rangle_{\mathfrak H}=0,\qquad  j=0,1.
\end{equation}
Denote $\delta\bm\alpha:=(\delta\alpha_1,\delta\alpha_2)$. For $(j,\ell)=(0,1)$ and $(j,\ell)=(1,0)$, the inner product is a multiple of $\delta\alpha_2-\ii\delta\alpha_1$ and $\delta\alpha_2+\ii\delta\alpha_1$, respectively.
\end{lemma}
\begin{proof}
For brevity, write $L_{j\ell}(\delta\bm\alpha):=\langle\widetilde C^{(j)},(\delta\bm\alpha\cdot\nabla_{\bm\alpha})\mathbb A^\varepsilon(\bm\alpha^*,\omega_*^\varepsilon) C^{(\ell)}\rangle_{\mathfrak H}$. Lemma \ref{lem:lem4-5} implies
\begin{equation}\label{eq:lem4-6-2}
L_{j\ell}(\delta\bm\alpha) = e^{2\pi \ii(\ell-j)/3}L_{j\ell}(R^{-1}\delta\bm\alpha).
\end{equation}
Let $L_-(\delta\bm\alpha):=\delta\alpha_2-\ii\delta\alpha_1$, $L_+(\delta\bm\alpha):=\delta\alpha_2+\ii\delta\alpha_1$. Direct calculation shows that
\begin{equation}\label{eq:lem4-6-1}
L_-(R^{-1}\delta\bm\alpha) = e^{-2\pi \ii/3}L_-(\delta\bm\alpha),  \qquad L_+(R^{-1}\delta\bm\alpha) = e^{2\pi \ii/3} L_+(\delta\bm\alpha).
\end{equation}
Notice that any complex-valued linear functional $L$ on $\mathbb R^2$ can be expressed uniquely as
$$
L(\delta\bm\alpha) = aL_-(\delta\bm\alpha)+bL_+(\delta\bm\alpha), \qquad  a,b\in\mathbb C.
$$
For the case of $(j,\ell)=(0,1)$, if we write $L_{01}(\delta\bm\alpha) = a_{01}L_-(\delta\bm\alpha)+b_{01}L_+(\delta\bm\alpha)$, then from \eqref{eq:lem4-6-2} and \eqref{eq:lem4-6-1},
$$
e^{2\pi \ii/3} L_{01}(R^{-1}\delta\bm\alpha) = e^{2\pi \ii/3}(a_{01}e^{-2\pi \ii/3}L_-(\delta\bm\alpha)+be^{2\pi \ii/3} L_+(\delta\bm\alpha)) = aL_-(\delta\bm\alpha) + be^{4\pi \ii/3} L_+(\delta\bm\alpha).
$$
Since $e^{2\pi \ii/3}\neq 1$, it immediately follows that $b=0$. The result for the case $(j,\ell)=(1,0)$ can be proved similarly. For $j=\ell$, with \eqref{eq:lem4-6-2} and \eqref{eq:lem4-6-1} we have
$$
a_{jj}L_-+b_{jj}L_+ = a_{jj}e^{-2\pi \ii/3}L_-+b_{jj}e^{2\pi \ii/3} L_+,
$$
which means $a_{jj}=b_{jj}=0$.
\end{proof}

We now calculate the coefficients of $\delta\alpha_2-\ii\delta\alpha_1$ and $\delta\alpha_2+\ii\delta\alpha_1$ for $(j,\ell)=(0,1)$ and $(j,\ell)=(1,0)$ cases. By~\eqref{eq:exact-Delta-m-root} and $z_*\ne0$,
$$
\frac{|z_*|}{2}\varepsilon^2 \le|\Delta_{\mathrm m}(\omega_*^\varepsilon)| \le\frac{3|z_*|}{2}\varepsilon^2
$$
for small $\varepsilon$. Near $\bm\alpha^*$, the part $G_{\Lambda_0^*,\omega_{\mathrm m}}^{\bm\alpha}$ of the exterior Green function is
$$
G_{\Lambda_0^*,\omega_{\mathrm m}}^{\bm\alpha}(\bm x) =\frac1{|Y|}\sum_{\ell=1}^3 \frac{e^{\ii(\bm\alpha+\bm q_\ell)\cdot\bm x}} {\omega_{\mathrm m}(\omega)^2 -(\bm\alpha+\bm q_\ell)\cdot(\bm\alpha+\bm q_\ell)}, \qquad  \bm\alpha^*+\bm q_\ell=\bm\xi_\ell.
$$
We remove the $t$-dependence and directly write $\delta\bm\alpha:=\bm\alpha-\bm\alpha^*$ and $\delta\omega:=\omega-\omega_*^\varepsilon$. If $|\delta\bm\alpha|+|\delta\omega|\le c|\Delta_{\mathrm m}(\omega_*^\varepsilon)|$ with sufficiently small $c>0$, the denominators of $G_{\Lambda_0^*,\omega_{\mathrm m}}^{\bm\alpha}$ are nonzero and
$$
\frac1{\omega_{\mathrm m}(\omega)^2 -|\bm\alpha+\bm q_\ell|^2}=\frac1{\Delta_{\mathrm m}(\omega_*^\varepsilon)} -\frac{ \partial_\omega(\omega_{\mathrm m}^2)(\omega_*^\varepsilon) \delta\omega-2\bm\xi_\ell\cdot\delta\bm\alpha} {\Delta_{\mathrm m}(\omega_*^\varepsilon)^2} +\mathcal{O} \left( \frac{(|\delta\bm\alpha|+|\delta\omega|)^2} {|\Delta_{\mathrm m}(\omega_*^\varepsilon)|^3} \right).
$$
Furthermore,
$$
e^{\ii(\bm\alpha+\bm q_\ell)\cdot\bm x} =e^{\ii\bm\xi_\ell\cdot\bm x} \left(1+\ii\delta\bm\alpha\cdot\bm x +\mathcal{O}(|\delta\bm\alpha|^2|\bm x|^2)\right).
$$
The phases in~\eqref{eq:AB-resonant-phases}, the explicit $\bm\xi_\ell$, and $\theta_1=-\pi/6$, $\theta_2=-5\pi/6$, $\theta_3=\pi/2$ give
\begin{align*}
\sum_{\ell=1}^3(\bm\xi_\ell\cdot\delta\bm\alpha) e^{\ii\bm\xi_\ell\cdot\bm d_{BA}}e^{2\ii\theta_\ell} &=\frac{2\pi}{a_{\rm lat}} (\ii\delta\alpha_1-\delta\alpha_2),\\
\sum_{\ell=1}^3(\bm\xi_\ell\cdot\delta\bm\alpha) e^{\ii\bm\xi_\ell\cdot\bm d_{AB}}e^{-2\ii\theta_\ell} &=-\frac{2\pi}{a_{\rm lat}} (\ii\delta\alpha_1+\delta\alpha_2).
\end{align*}

\begin{lemma}\label{lem:boundary-dirac-coupling}
Normalize the right vectors $C^{(j)}$ by $c_{1,B,-1}^{(0)}=1$ and $c_{1,A,1}^{(1)}=1$, and the left vectors $\widetilde C^{(j)}$, $j=0,1$, by~\eqref{eq:boundary-frequency-normalization}. Then
\begin{align*}
\left\langle\widetilde C^{(0)}, \partial_{\alpha_1}\mathbb A_*^\varepsilon C^{(1)}\right\rangle_{\mathfrak H} &=-\ii \frac{4\pi} {3a_{\rm lat} \partial_\omega(\omega_{\mathrm m}^2)(\omega_*^\varepsilon)} \left(1+\mathcal{O} \left(|\ln\varepsilon|^{-1}\right)\right),\\
\left\langle\widetilde C^{(0)}, \partial_{\alpha_2}\mathbb A_*^\varepsilon C^{(1)}\right\rangle_{\mathfrak H} &=\frac{4\pi} {3a_{\rm lat} \partial_\omega(\omega_{\mathrm m}^2)(\omega_*^\varepsilon)} \left(1+\mathcal{O} \left(|\ln\varepsilon|^{-1}\right)\right).
\end{align*}
\end{lemma}

\begin{proof}
To determine the normalization factor, we first calculate the frequency derivatives, then the Bloch derivative, and combine the two calculations. By \eqref{eq:exact-res-S}, at $(\bm\alpha^*,\omega_*^\varepsilon)$ we have
\begin{align}
\partial_\omega (\widehat{\mathcal{S}}_{\Lambda_0^*,\omega_{\mathrm m}}^{BB,\varepsilon})_{-1,-1} &=-\frac{6\pi J_1(z_D)^2} {|Y|\Delta_{\mathrm m}(\omega_*^\varepsilon)^2} \partial_\omega(\omega_{\mathrm m}^2)(\omega_*^\varepsilon), \notag\\
(\delta\bm\alpha\cdot\nabla_{\bm\alpha}) (\widehat{\mathcal{S}}_{\Lambda_0^*,\omega_{\mathrm m}}^{BA,\varepsilon})_{-1,1} &=\frac{4\pi J_1(z_D)^2} {|Y|\Delta_{\mathrm m}(\omega_*^\varepsilon)^2} \sum_{\ell=1}^{3}(\bm\xi_\ell\cdot\delta\bm\alpha) e^{\ii\bm\xi_\ell\cdot\bm d_{BA}}e^{2\ii\theta_\ell} +\mathcal{O}(|\delta\bm\alpha|). \label{eq:partial-omega-s}
\end{align}
Similarly, \eqref{eq:exact-res-K} gives
\begin{align}
\partial_\omega (\mathcal{K}_{\Lambda_0^*,\omega_{\mathrm m}}^{BB,\varepsilon,*})_{-1,-1} &=-\frac{6\pi\varepsilon k_DJ_1'(z_D)J_1(z_D)} {|Y|\Delta_{\mathrm m}(\omega_*^\varepsilon)^2} \partial_\omega(\omega_{\mathrm m}^2)(\omega_*^\varepsilon), \notag\\
(\delta\bm\alpha\cdot\nabla_{\bm\alpha}) (\mathcal{K}_{\Lambda_0^*,\omega_{\mathrm m}}^{BA,\varepsilon,*})_{-1,1} &=\frac{4\pi\varepsilon k_DJ_1'(z_D)J_1(z_D)} {|Y|\Delta_{\mathrm m}(\omega_*^\varepsilon)^2} \sum_{\ell=1}^{3}(\bm\xi_\ell\cdot\delta\bm\alpha) e^{\ii\bm\xi_\ell\cdot\bm d_{BA}}e^{2\ii\theta_\ell} +\mathcal{O}(|\delta\bm\alpha|). \label{eq:partial-omega-k}
\end{align}
Substituting $\sum_{\ell=1}^{3}(\bm\xi_\ell\cdot\delta\bm\alpha) e^{\ii\bm\xi_\ell\cdot\bm d_{BA}}e^{2\ii\theta_\ell} =(2\pi/a_{\rm lat})(\ii\delta\alpha_1-\delta\alpha_2)$ into the single-layer and normal-derivative entries, and dividing each by its corresponding $\omega$-derivative, gives
\begin{align*}
\frac{ (\delta\bm\alpha\cdot\nabla_{\bm\alpha}) (\widehat{\mathcal{S}}_{\Lambda_0^*,\omega_{\mathrm m}}^{BA,\varepsilon})_{-1,1}} {\partial_\omega (\widehat{\mathcal{S}}_{\Lambda_0^*,\omega_{\mathrm m}}^{BB,\varepsilon})_{-1,-1}} &=-\frac{4\pi(\ii\delta\alpha_1-\delta\alpha_2)} {3a_{\rm lat}  \partial_\omega(\omega_{\mathrm m}^2)(\omega_*^\varepsilon)} +\mathcal{O}(\varepsilon^2|\delta\bm\alpha|),\\
\frac{ (\delta\bm\alpha\cdot\nabla_{\bm\alpha}) (\mathcal{K}_{\Lambda_0^*,\omega_{\mathrm m}}^{*,BA,\varepsilon})_{-1,1}} {\partial_\omega (\mathcal{K}_{\Lambda_0^*,\omega_{\mathrm m}}^{*,BB,\varepsilon})_{-1,-1}} &=-\frac{4\pi(\ii\delta\alpha_1-\delta\alpha_2)} {3a_{\rm lat}  \partial_\omega(\omega_{\mathrm m}^2)(\omega_*^\varepsilon)} +\mathcal{O}(\varepsilon^2|\delta\bm\alpha|).
\end{align*}
As a result,
$$
(\delta\bm\alpha\cdot\nabla_{\bm\alpha})\mathsf B_{B,-1;A,1}^\varepsilon = \frac{8\pi^2 J_1(z_D)^2}{a_{\mathrm{lat}}|Y|\Delta_m^2}(\ii\delta\alpha_1-\delta\alpha_2)
\begin{pmatrix}
1 & 0 \\
-\frac{\varepsilon k_D J_1'(z_D)}{\mu_mJ_1(z_D)} & 0
\end{pmatrix}
+ \mathcal{O}(|\delta\bm\alpha|).
$$

It remains to estimate the derivatives of the terms produced by Schur elimination. Notice that
$$
\partial_\omega\left[ L_j^\varepsilon(\widetilde Z_j^\varepsilon)^{-1}R_j^\varepsilon \right]= (\partial_\omega L_j^\varepsilon)(\widetilde Z_j^\varepsilon)^{-1}R_j^\varepsilon-L_j^\varepsilon(\widetilde Z_j^\varepsilon)^{-1} (\partial_\omega\widetilde Z_j^\varepsilon) (\widetilde Z_j^\varepsilon)^{-1}R_j^\varepsilon+ L_j^\varepsilon(\widetilde Z_j^\varepsilon)^{-1} (\partial_\omega R_j^\varepsilon).
$$
Recall that the estimate \eqref{eq:constant-elimination-exact-order} holds on the larger complex $z$-disk specified before Theorem~\ref{thm:exact-dirac-frequency}. Cauchy's estimate on a smaller disk, followed by $\partial_\omega z=\varepsilon^{-2}\partial_\omega(\omega_{\mathrm m}^2)$, therefore gives
\begin{equation}\label{eq:freq-derivative-bound}
\left\| \partial_\omega \left[ L_j^\varepsilon(\widetilde Z_j^\varepsilon)^{-1}R_j^\varepsilon \right] \right\| \le \frac{C\varepsilon^{-2}}{|\ln\varepsilon|}, \qquad  j=0,1.
\end{equation}

We now compute the Bloch derivative of the coupling from the mode $(A,1)$ in sector $1$ to the mode $(B,-1)$ in sector $0$. By Theorem~\ref{thm:mod3-invariant-subspaces}, $\mathbb A_*^\varepsilon$ preserves each symmetry sector. However, its Bloch derivative $\partial_{\alpha_k}\mathbb A_*^\varepsilon$ need not preserve this decomposition and may couple different sectors. We therefore cannot obtain the off-diagonal Bloch derivative by differentiating the sectorwise reduced matrices separately.

We first eliminate sector $2$ and show that this elimination does not contribute to the first-order coupling between sectors $0$ and $1$. Recall the projections $\Pi_0$ and $\Pi_1$ defined before Proposition~\ref{prop:six-dimensional-reduction}, and set $\Pi_2:=I-\Pi_0-\Pi_1$, the orthogonal projection onto $\mathfrak H^{(2)}$. By Theorem~\ref{thm:exact-dirac-frequency}, $\Pi_2\mathbb A_*^\varepsilon\Pi_2$ is invertible. Consequently, for $\bm\alpha$ sufficiently close to $\bm\alpha^*$, the restriction $\Pi_2\mathbb A^\varepsilon\Pi_2$ remains invertible at $\omega=\omega_*^\varepsilon$. The block mapping sector $1$ into sector $0$ in the Schur complement obtained by eliminating sector $2$ is
$$
\Pi_0\mathbb A^\varepsilon\Pi_1 -\Pi_0\mathbb A^\varepsilon\Pi_2 (\Pi_2\mathbb A^\varepsilon\Pi_2)^{-1} \Pi_2\mathbb A^\varepsilon\Pi_1,
$$
where all operators are evaluated at $(\bm\alpha,\omega_*^\varepsilon)$. The first term represents the direct coupling between sectors $1$ and $0$, while the second accounts for the indirect coupling through sector $2$. We now show that the Bloch derivative of the second term vanishes at $\bm\alpha=\bm\alpha^*$. The remaining coefficients within sectors $0$ and $1$ will then be eliminated using Proposition~\ref{prop:six-dimensional-reduction}, Lemma~\ref{lem:second-mode-self-block}, and~\eqref{eq:Ztilde-inverse-orders}.

Differentiating the second term with respect to $\alpha_k$ gives
\begin{align*}
&\Pi_0(\partial_{\alpha_k}\mathbb A^\varepsilon)\Pi_2 (\Pi_2\mathbb A^\varepsilon\Pi_2)^{-1} \Pi_2\mathbb A^\varepsilon\Pi_1 \\
+&\Pi_0\mathbb A^\varepsilon\Pi_2 \partial_{\alpha_k} \left[(\Pi_2\mathbb A^\varepsilon\Pi_2)^{-1}\right] \Pi_2\mathbb A^\varepsilon\Pi_1\\
+&\Pi_0\mathbb A^\varepsilon\Pi_2 (\Pi_2\mathbb A^\varepsilon\Pi_2)^{-1} \Pi_2(\partial_{\alpha_k}\mathbb A^\varepsilon)\Pi_1,
\end{align*}
which is identically 0 at $\bm\alpha=\bm\alpha^*$ since $\Pi_i\mathbb A_*^\varepsilon\Pi_j=0$ for $i\neq j$. The same argument applies to the other retained blocks. These sector-2 terms do not contribute at first order.

After Schur elimination of the constant modes, the mixed Bloch derivative between the retained dipole modes $(A,1)$ and $(B,-1)$ is
\begin{align}
&
\begin{pmatrix}
-L_0^\varepsilon(\widetilde{Z_0^\varepsilon})^{-1} & I_2
\end{pmatrix}
\begin{pmatrix}
\partial_{\alpha_k}\mathsf B_{A,0;B,0}^\varepsilon & \partial_{\alpha_k}\mathsf B_{A,0;A,1}^\varepsilon \\
\partial_{\alpha_k}\mathsf B_{B,-1;B,0}^\varepsilon & \partial_{\alpha_k}\mathsf B_{B,-1;A,1}^\varepsilon
\end{pmatrix}
\begin{pmatrix}
-(\widetilde{Z_1^\varepsilon})^{-1}R_1^\varepsilon \\
I_2
\end{pmatrix}
\notag\\
=&\partial_{\alpha_k}\mathsf B_{B,-1;A,1}^{\varepsilon} -\bigl(\partial_{\alpha_k}\mathsf B_{B,-1;B,0}^{\varepsilon}\bigr) (\widetilde Z_1^\varepsilon)^{-1}R_1^\varepsilon- L_0^\varepsilon(\widetilde Z_0^\varepsilon)^{-1} \bigl(\partial_{\alpha_k}\mathsf B_{A,0;A,1}^{\varepsilon}\bigr) \notag\\
&\quad + L_0^\varepsilon(\widetilde Z_0^\varepsilon)^{-1} \bigl(\partial_{\alpha_k}\mathsf B_{A,0;B,0}^{\varepsilon}\bigr) (\widetilde Z_1^\varepsilon)^{-1}R_1^\varepsilon +\mathcal{O}(\varepsilon^{-1}). \label{eq:boundary-mixed-schur-derivative}
\end{align}
The remainder accounts for the differentiated terms involving the $|n|=2$ modes and $Q_0\mathfrak H\oplus Q_1\mathfrak H$.

To check the bounds in \eqref{eq:boundary-mixed-schur-derivative}, \eqref{eq:Ztilde-inverse-orders} and the $L_j^\varepsilon,R_j^\varepsilon$ entries used in \eqref{eq:constant-elimination-exact-order} give
$$
(\widetilde Z_j^\varepsilon)^{-1}R_j^\varepsilon=
\begin{pmatrix}
\mathcal{O}(\varepsilon/|\ln\varepsilon|)&\mathcal{O}(\varepsilon^2)\\
\mathcal{O}(\varepsilon^{-1}/|\ln\varepsilon|)&\mathcal{O}(|\ln\varepsilon|^{-1})
\end{pmatrix},\quad  L_j^\varepsilon(\widetilde Z_j^\varepsilon)^{-1} =
\begin{pmatrix}
\mathcal{O}(\varepsilon)&\mathcal{O}(\varepsilon^{-1})\\
\mathcal{O}(\varepsilon)&\mathcal{O}(\varepsilon^{-1})
\end{pmatrix}.
$$
For a constant target mode, \eqref{eq:exact-res-S} and \eqref{eq:exact-res-K} also give the exact identity
$$
\begin{pmatrix}
(\widehat{\mathcal{S}}_{\Lambda_0^*,\omega_{\mathrm m}}^{pq,\varepsilon})_{0,n}\\
-\mu_{\mathrm m}^{-1} (\mathcal{K}_{\Lambda_0^*,\omega_{\mathrm m}}^{*,pq,\varepsilon})_{0,n}
\end{pmatrix}
=
\begin{pmatrix}
1\\
-\dfrac{\varepsilon k_DJ_0'(z_D)}{\mu_{\mathrm m}J_0(z_D)}
\end{pmatrix}
(\widehat{\mathcal{S}}_{\Lambda_0^*,\omega_{\mathrm m}}^{pq,\varepsilon})_{0,n}.
$$
Since $J_0'=-J_1$, \eqref{eq:Ztilde-inverse-orders} implies
\begin{equation}\label{eq:boundary-constant-derivative-cancellation}
L_j^\varepsilon(\widetilde Z_j^\varepsilon)^{-1}
\begin{pmatrix}
1\\
-\dfrac{\varepsilon k_DJ_0'(z_D)}{\mu_{\mathrm m}J_0(z_D)}
\end{pmatrix}
=\mathcal{O} \left(\frac{\varepsilon}{|\ln\varepsilon|}\right).
\end{equation}

Differentiating the denominator $\omega_m^2-|\bm\alpha+\bm q_\ell|^2$ in $G_{\Lambda_0^*,\omega_m}$ gives a term of order $\mathcal{O}(\varepsilon^{-4})$. Thus its single-layer contributions between a dipole and a constant, and between two constants, have orders $\mathcal{O}(\varepsilon^{-3})$ and $\mathcal{O}(\varepsilon^{-4})$, respectively, by \eqref{eq:exact-res-S}. By \eqref{eq:partial-omega-s}, \eqref{eq:partial-omega-k} and \eqref{eq:boundary-constant-derivative-cancellation}, both the second and third terms of \eqref{eq:boundary-mixed-schur-derivative} are of order $\mathcal{O}(\varepsilon^{-2}|\ln\varepsilon|^{-1})$. For the fourth term, a similar calculation gives $\mathcal{O}(\varepsilon^{-2}|\ln\varepsilon|^{-2})$.

Now we study the order of the modes with $|n|=2$ and higher Fourier modes in $Q_0\mathfrak H\oplus Q_1\mathfrak H$. We use \eqref{eq:exact-res-S}--\eqref{eq:exact-res-K} and \eqref{eq:bessel-estimate} to control those remaining terms. Straightforward but long calculations show that the contributions involving $|n|=2$ are bounded by $\mathcal{O}(\varepsilon^{-1})$, while those involving $Q_j\mathfrak H$ are bounded by $\mathcal{O}(1)$. This proves the stated remainder. Consequently, the off-diagonal Bloch-derivative block differs from $\partial_{\alpha_k}\mathsf B_{B,-1;A,1}^{\varepsilon}$ by $\mathcal{O}(\varepsilon^{-2}|\ln\varepsilon|^{-1})$. The reverse block is treated in the same way. Together with the frequency-derivative estimate \eqref{eq:freq-derivative-bound}, we control the derivatives of the full Schur equation on the two retained dipole modes $(A,1)$ and $(B,-1)$. To conclude, \eqref{eq:boundary-mixed-schur-derivative} equals
\begin{equation}\label{eq:boundary-mixed-schur-derivative-result}
M_{01}^\varepsilon(\delta\bm\alpha):=\frac{8\pi^2 J_1(z_D)^2}{a_{\mathrm{lat}}|Y|\Delta_m^2}(\ii\delta\alpha_1-\delta\alpha_2)
\begin{pmatrix}
1 & 0 \\
-\frac{\varepsilon k_D J_1'(z_D)}{\mu_mJ_1(z_D)} & 0
\end{pmatrix}
+ \mathcal{O}(|\delta\bm\alpha|)+\mathcal{O}\left(\frac{\varepsilon^{-2}|\delta\bm\alpha|}{|\ln\varepsilon|}\right).
\end{equation}

Finally, we calculate $\widetilde C^{(0)}$ and $C^{(1)}$. Because the preceding Schur reductions reconstruct the eliminated coefficients, it suffices to calculate the reduced nullvectors of $\mathcal D_j^\varepsilon := D_j^\varepsilon-L_j^\varepsilon(\widetilde{Z}_j^\varepsilon)^{-1}R_j^\varepsilon$. Using $J_1(z_D)=Q\varepsilon/3+\mathcal{O}(\varepsilon^3)$ and \eqref{eq:exact-Delta-m-root} gives $\Delta_{\mathrm m}(\omega_*^\varepsilon) =z_*\varepsilon^2+\mathcal{O}(\varepsilon^2|\ln\varepsilon|^{-1})$. The right nullvector $r_j^\varepsilon$ and a left nullrow $\ell_j^\varepsilon$ of $\mathcal D_j^\varepsilon$ have, in the coefficient order of \eqref{eq:four-dim-matrix}, the leading values
$$
r_j^\varepsilon =
\begin{pmatrix}
1\\
1-2a(z_*)
\end{pmatrix}
+\mathcal{O}(|\ln\varepsilon|^{-1}), \qquad  \ell_j^\varepsilon =
\begin{pmatrix}
1&\mu
\end{pmatrix}
+\mathcal{O}(|\ln\varepsilon|^{-1}),
$$
by \eqref{eq:constant-elimination-exact-order} and the displayed formula for $D_j^\varepsilon$ immediately preceding \eqref{eq:exact-scalar-characteristic-expansion}. This also shows that the first-density normalization in the lemma is possible. Contracting the frequency derivative with these vectors, before normalizing the left row, gives the normalization factor
$$
\ell_0^\varepsilon(\partial_\omega\mathcal D_0^\varepsilon)r_0^\varepsilon = -\frac{6\pi J_1(z_D)^2} {|Y|\Delta_{\mathrm m}(\omega_*^\varepsilon)^2} \partial_\omega(\omega_{\mathrm m}^2)(\omega_*^\varepsilon) \left(1-\frac{\mu}{\mu_{\mathrm m}}\right) \left(1+\mathcal{O}(|\ln\varepsilon|^{-1})\right).
$$
Multiplying~\eqref{eq:boundary-mixed-schur-derivative-result} on the left by $\ell_0^\varepsilon$ and on the right by $r_1^\varepsilon$ gives
$$
\ell_0^\varepsilon M_{01}^\varepsilon(\delta\bm\alpha)r_1^\varepsilon= \frac{8\pi^2J_1(z_D)^2}{a_{\rm lat}|Y|\Delta_{\mathrm m}^2}\left(1-\frac{\mu}{\mu_{\mathrm m}}\right)\left(1+\mathcal{O}(|\ln\varepsilon|^{-1})\right)(\ii\delta\alpha_1-\delta\alpha_2).
$$
Then
\begin{align*}
\left\langle\widetilde C^{(0)}, (\delta\bm\alpha\cdot\nabla_{\bm\alpha}) \mathbb A_*^\varepsilon C^{(1)}\right\rangle_{\mathfrak H} &= \frac{\ell_0^\varepsilon M_{01}^\varepsilon(\delta\bm\alpha)r_1^\varepsilon}{\ell_0^\varepsilon(\partial_\omega\mathcal D_0^\varepsilon)r_0^\varepsilon} \\
&= \frac{4\pi} {3a_{\rm lat} \partial_\omega(\omega_{\mathrm m}^2)(\omega_*^\varepsilon)} \left(1+\mathcal{O} \left(|\ln\varepsilon|^{-1}\right)\right) (\delta\alpha_2-\ii\delta\alpha_1).
\end{align*}
\end{proof}

The preceding lemma determines the first-order Bloch coupling from the mode $(A,1)$ in sector $1$ to the mode $(B,-1)$ in sector $0$. To obtain the complete two-dimensional reduced system, we also need the coupling in the opposite direction, from sector $0$ to sector $1$. Rather than repeating the full calculation, we first determine its dependence on the Bloch increment $\delta\bm\alpha$ and then use the inversion--conjugation symmetry to relate its coefficient to the one already computed.

For the reverse coupling, the same calculation as in Lemma~\ref{lem:boundary-dirac-coupling} applies, with the roles of the two boundary components exchanged. The corresponding resonant phase sum is
$$
\sum_{\ell=1}^{3} (\bm\xi_\ell\cdot\delta\bm\alpha) e^{\ii\bm\xi_\ell\cdot\bm d_{AB}}e^{-2\ii\theta_\ell} = -\frac{2\pi}{a_{\rm lat}} (\ii\delta\alpha_1+\delta\alpha_2).
$$
Comparing this identity with the phase sum used for the coupling from $(A,1)$ to $(B,-1)$, we see that the factor $\delta\alpha_2-\ii\delta\alpha_1$ is replaced by $\delta\alpha_2+\ii\delta\alpha_1$. Applying the same Schur reduction and estimates as in the preceding lemma therefore gives
$$
\left\langle\widetilde C^{(1)}, (\delta\bm\alpha\cdot\nabla_{\bm\alpha}) \mathbb A_*^\varepsilon C^{(0)} \right\rangle_{\mathfrak H}= \left\langle\widetilde C^{(1)}, \partial_{\alpha_2}\mathbb A_*^\varepsilon C^{(0)} \right\rangle_{\mathfrak H} (\delta\alpha_2+\ii\delta\alpha_1).
$$
Thus the dependence on $\delta\bm\alpha$ is known. It remains to relate the coefficient on the right-hand side to the corresponding coefficient for the coupling from sector $1$ to sector $0$.

For this purpose, recall the inversion--conjugation symmetry introduced in~\eqref{eq:J-definition}. By~\eqref{eq:J-mapping}, this symmetry exchanges the boundary components $A$ and $B$ and changes the Fourier index $n$ into $-n$. In particular, it maps sector $0$ onto sector $1$. Since the Fourier--Sobolev weights defining $\mathcal U_X$ and $\mathcal U_Y$ depend only on $|n|$, their induced actions on the coefficient space $\mathfrak H$ coincide and satisfy
$$
\bigl(\mathcal U_X\mathcal J_X \mathcal U_X^{-1}C\bigr)_{\sigma,p,n}= \bigl(\mathcal U_Y\mathcal J_Y \mathcal U_Y^{-1}C\bigr)_{\sigma,p,n} =(-1)^n\overline{C_{\sigma,\overline p,-n}}.
$$
These transformations are antilinear isometries of $\mathfrak H$.

The intertwining relation~\eqref{eq:J2-commutation} expresses the invariance of the nonmagnetic boundary operator under this symmetry. The Green-function identity used in its proof holds for real Bloch vectors in a neighborhood of $\bm\alpha^*$, so the same intertwining relation holds there and can be differentiated with respect to $\alpha_k$. At $(\bm\alpha^*,\omega_*^\varepsilon)$, this yields
$$
(\partial_{\alpha_k}\mathbb A_*^\varepsilon) \mathcal U_X\mathcal J_X\mathcal U_X^{-1}= \mathcal U_Y\mathcal J_Y\mathcal U_Y^{-1} (\partial_{\alpha_k}\mathbb A_*^\varepsilon).
$$
Moreover, because the symmetry exchanges sectors $0$ and $1$, the right and left nullvectors can be chosen consistently with the normalizations used in Lemma~\ref{lem:boundary-dirac-coupling} so that
$$
C^{(1)} =-\mathcal U_X\mathcal J_X \mathcal U_X^{-1}C^{(0)},\qquad  \widetilde C^{(1)} =-\mathcal U_Y\mathcal J_Y \mathcal U_Y^{-1}\widetilde C^{(0)}.
$$
The minus signs agree with the first-density normalizations of the modes $(B,-1)$ and $(A,1)$.

Using these relations, the differentiated intertwining identity, and the antiunitary property of the coefficient-space transformations, we obtain
$$
\left\langle\widetilde C^{(1)}, \partial_{\alpha_2}\mathbb A_*^\varepsilon C^{(0)}\right\rangle_{\mathfrak H}= \overline{ \left\langle\widetilde C^{(0)}, \partial_{\alpha_2}\mathbb A_*^\varepsilon C^{(1)}\right\rangle_{\mathfrak H} }.
$$
Consequently, the coefficients of the two off-diagonal Bloch couplings are complex conjugates. In particular, the reverse coupling is completely determined by the result of Lemma~\ref{lem:boundary-dirac-coupling}, without a separate calculation of its leading coefficient.

Finally, we fix the remaining phase freedom so that the two couplings have the same positive real coefficient. After the first-density normalization used for the calculation in the preceding lemma, we may rephase the right and left nullvectors in sector $0$ by a common factor $e^{-\ii\chi}$ and those in sector $1$ by $e^{\ii\chi}$. These opposite phase changes preserve both the inversion--conjugation relations and the frequency normalization \eqref{eq:boundary-frequency-normalization}. Choosing $\chi$ appropriately, we obtain
$$
\left\langle\widetilde C^{(0)}, \partial_{\alpha_2}\mathbb A_*^\varepsilon C^{(1)}\right\rangle_{\mathfrak H} = \left\langle\widetilde C^{(1)}, \partial_{\alpha_2}\mathbb A_*^\varepsilon C^{(0)}\right\rangle_{\mathfrak H} =:v_D^\varepsilon>0.
$$
By Lemma~\ref{lem:boundary-dirac-coupling}, the magnitude of this coefficient is
$$
v_D^\varepsilon = \frac{4\pi} {3a_{\rm lat} |\partial_\omega(\omega_{\mathrm m}^2) (\omega_*^\varepsilon)|} \left( 1+\mathcal{O} \left(|\ln\varepsilon|^{-1}\right) \right).
$$
Since $\partial_\omega(\omega_{\mathrm m}^2) =2\omega_{\mathrm m}\omega_{\mathrm m}'$, we conclude that
\begin{equation}\label{eq:newK-dirac-velocity}
v_D^\varepsilon = \frac{2\pi} {3a_{\rm lat} |\omega_{\mathrm m}(\omega_*^\varepsilon) \omega_{\mathrm m}'(\omega_*^\varepsilon)|} \left( 1+\mathcal{O} \left(|\ln\varepsilon|^{-1}\right) \right)>0.
\end{equation}
Thus the two off-diagonal Bloch couplings have the respective forms $v_D^\varepsilon(\delta\alpha_2-\ii\delta\alpha_1)$ and $v_D^\varepsilon(\delta\alpha_2+\ii\delta\alpha_1)$. Their conjugate structure is the key ingredient in the two-dimensional reduced system used below to establish the Dirac conical dispersion.

\begin{theorem}
Under the assumptions of Theorem~\ref{thm:exact-dirac-frequency}, for each sufficiently small fixed $\varepsilon>0$ there exists $r_D^\varepsilon>0$ such that, when $|\bm\alpha-\bm\alpha^*|<r_D^\varepsilon$, there are two characteristic-frequency branches
\begin{equation}\label{eq:dirac-branches-exact}
\omega_\pm^\varepsilon(\bm\alpha) = \omega_*^\varepsilon \pm v_D^\varepsilon|\bm\alpha-\bm\alpha^*| +\mathcal{O}_\varepsilon \left(|\bm\alpha-\bm\alpha^*|^2\right),
\end{equation}
where $v_D^\varepsilon>0$ is given by \eqref{eq:newK-dirac-velocity}, and we may take $r_D^\varepsilon\le C\varepsilon^2$. By \eqref{eq:normalized-kernel-equivalence} and Proposition~\ref{prop:pde-bie-equivalence}, \eqref{eq:dirac-branches-exact} also gives the two characteristic-frequency branches of the original PDE Bloch problem. Hence $(\bm\alpha^*,\omega_*^\varepsilon)$ is a local Dirac point.
\end{theorem}

\begin{proof}
By Theorem~\ref{thm:exact-dirac-frequency}, $\ker\mathbb A_*^\varepsilon =\operatorname{span}\{C^{(0)},C^{(1)}\}$. Since $\mathbb A_*^\varepsilon$ is Fredholm of index zero, its range is closed and $\operatorname{Ran}\mathbb A_*^\varepsilon =(\ker(\mathbb A_*^\varepsilon)^*)^\perp$. Moreover, $\operatorname{ind}\mathbb A_*^\varepsilon=0$ and $\dim\ker\mathbb A_*^\varepsilon=2$ imply $\dim\ker(\mathbb A_*^\varepsilon)^*=2$. Therefore
$$
\left.\mathbb A_*^\varepsilon\right|_{(\ker\mathbb A_*^\varepsilon)^\perp} : (\ker\mathbb A_*^\varepsilon)^\perp \longrightarrow \operatorname{Ran}\mathbb A_*^\varepsilon
$$
is bijective, and its inverse $(\mathbb A_*^\varepsilon)^{-1}$ is bounded by the bounded inverse theorem.

Recall $\delta\bm\alpha=\bm\alpha-\bm\alpha^*$, $\delta\omega=\omega-\omega_*^\varepsilon$. Write $C=a_0C^{(0)}+a_1C^{(1)}+C^\perp$, where $C^\perp\in(\ker\mathbb A_*^\varepsilon)^\perp$. Let $P_Y$ denote the orthogonal projection onto $\operatorname{Ran}\mathbb A_*^\varepsilon$, and set $\delta\mathbb A^\varepsilon:=\mathbb A^\varepsilon(\bm\alpha,\omega) -\mathbb A_*^\varepsilon$. Since $\mathbb A_*^\varepsilon C^{(j)}=0$, $j=0,1$, applying $P_Y$ to $\mathbb A^\varepsilon(\bm\alpha,\omega)C=0$ gives
$$
\mathbb A_*^\varepsilon C^\perp +P_Y\delta\mathbb A^\varepsilon C^\perp =-P_Y\delta\mathbb A^\varepsilon(a_0C^{(0)}+a_1C^{(1)}).
$$
Applying $(\mathbb A_*^\varepsilon)^{-1}$ yields
$$
(I+(\mathbb A_*^\varepsilon)^{-1}P_Y\delta\mathbb A^\varepsilon)C^\perp=-(\mathbb A_*^\varepsilon)^{-1}P_Y\delta\mathbb A^\varepsilon(a_0C^{(0)}+a_1C^{(1)}).
$$
For fixed $\varepsilon$, the smooth dependence of $\mathbb A^\varepsilon$ on its parameters implies
$$
\|\delta\mathbb A^\varepsilon\|_{\mathcal L(\mathfrak H)} \le C_\varepsilon (|\delta\bm\alpha|+|\delta\omega|).
$$
Consequently, when $|\delta\bm\alpha|+|\delta\omega|$ is sufficiently small, we have $\|(\mathbb A_*^\varepsilon)^{-1}P_Y\delta\mathbb A^\varepsilon\|<1$. The operator $I+(\mathbb A_*^\varepsilon)^{-1}P_Y\delta\mathbb A^\varepsilon$ is therefore invertible by a Neumann series, which uniquely determines $C^\perp$ in terms of $a_0,a_1$ and gives
$$
\|C^\perp\|_{\mathfrak H} \le C_\varepsilon (|\delta\bm\alpha|+|\delta\omega|) (|a_0|+|a_1|).
$$

Expand $\mathbb A^\varepsilon(\bm\alpha,\omega)$ at $(\bm\alpha^*,\omega_*^\varepsilon)$ and pair with $\widetilde C^{(j)}\in\ker((\mathbb A_*^\varepsilon)^*)$, $j=0,1$. Using $\mathbb A_*^\varepsilon C^{(\ell)}=0$, $(\mathbb A_*^\varepsilon)^*\widetilde C^{(j)}=0$, and the estimate for $C^\perp$, we obtain
\begin{align*}
0={}& \sum_{\ell=0}^{1}a_\ell \left[ \delta\omega \left\langle\widetilde C^{(j)}, \partial_\omega\mathbb A_*^\varepsilon C^{(\ell)} \right\rangle_{\mathfrak H} + \left\langle\widetilde C^{(j)}, (\delta\bm\alpha\cdot\nabla_{\bm\alpha}) \mathbb A_*^\varepsilon C^{(\ell)} \right\rangle_{\mathfrak H} \right]\\
&+ \mathcal{O}_\varepsilon \left( (|\delta\bm\alpha|+|\delta\omega|)^2 \right)(|a_0|+|a_1|), \qquad  j=0,1.
\end{align*}
By \eqref{eq:boundary-frequency-normalization},
$$
\left\langle\widetilde C^{(j)}, \partial_\omega\mathbb A_*^\varepsilon C^{(\ell)} \right\rangle_{\mathfrak H} =\delta_{j\ell}.
$$
Also, \eqref{eq:boundary-bloch-diagonal-zero} gives
$$
\left\langle\widetilde C^{(j)}, (\delta\bm\alpha\cdot\nabla_{\bm\alpha}) \mathbb A_*^\varepsilon C^{(j)} \right\rangle_{\mathfrak H}=0, \qquad  j=0,1,
$$
and the phase choice above gives
\begin{align*}
\left\langle\widetilde C^{(0)}, (\delta\bm\alpha\cdot\nabla_{\bm\alpha}) \mathbb A_*^\varepsilon C^{(1)} \right\rangle_{\mathfrak H} &=v_D^\varepsilon(\delta\alpha_2-\ii\delta\alpha_1),\\
\left\langle\widetilde C^{(1)}, (\delta\bm\alpha\cdot\nabla_{\bm\alpha}) \mathbb A_*^\varepsilon C^{(0)} \right\rangle_{\mathfrak H} &=v_D^\varepsilon(\delta\alpha_2+\ii\delta\alpha_1).
\end{align*}
Consequently,
$$
\left[
\begin{pmatrix}
\delta\omega& v_D^\varepsilon(\delta\alpha_2-\ii\delta\alpha_1)\\
v_D^\varepsilon(\delta\alpha_2+\ii\delta\alpha_1)& \delta\omega
\end{pmatrix}
+ \mathcal{O}_\varepsilon \left( (|\delta\bm\alpha|+|\delta\omega|)^2 \right) \right]
\begin{pmatrix}
a_0\\
a_1
\end{pmatrix}
=0.
$$
A nonzero vector $(a_0,a_1)^\top$ exists if and only if the determinant of this two-dimensional matrix vanishes. Thus
$$
(\delta\omega)^2 -(v_D^\varepsilon)^2|\delta\bm\alpha|^2 + \mathcal{O}_\varepsilon \left( (|\delta\bm\alpha|+|\delta\omega|)^3 \right)=0.
$$

For $\delta\bm\alpha\ne\bm0$, let $r=|\delta\bm\alpha|$, $\bm e=\delta\bm\alpha/r$, and $\delta\omega=r\zeta$. For bounded $|\delta\omega|/|\delta\bm\alpha|$, dividing the preceding equation by $|\delta\bm\alpha|^2$ gives
$$
\frac{(\delta\omega)^2}{|\delta\bm\alpha|^2}-(v_D^\varepsilon)^2+\mathcal{O}_\varepsilon(|\delta\bm\alpha|)=0.
$$
At $r=0$, the two roots are $\zeta=\pm v_D^\varepsilon$, and $2\zeta\ne0$. By the implicit function theorem, there is a sufficiently small range of $r$, independent of $\delta\bm\alpha/r$, in which $\zeta_\pm=\pm v_D^\varepsilon+\mathcal{O}_\varepsilon(r)$. Therefore
$$
\omega_\pm^\varepsilon(\bm\alpha) = \omega_*^\varepsilon \pm v_D^\varepsilon|\bm\alpha-\bm\alpha^*| + \mathcal{O}_\varepsilon \left(|\bm\alpha-\bm\alpha^*|^2\right).
$$

By \eqref{eq:exact-Delta-m-root}, $\Delta_{\mathrm m}(\omega_*^\varepsilon) =z_*\varepsilon^2+\mathcal{O}(\varepsilon^2|\ln\varepsilon|^{-1})$, with $z_*\ne0$. For $\bm\alpha=\bm\alpha^*+\delta\bm\alpha$ and $\omega=\omega_*^\varepsilon+\delta\omega$,
$$
\omega_{\mathrm m}(\omega)^2-|\bm\alpha+\bm q_\ell|^2 =\Delta_{\mathrm m}(\omega_*^\varepsilon) +\partial_\omega(\omega_{\mathrm m}^2)(\omega_*^\varepsilon) \delta\omega -2\bm\xi_\ell\cdot\delta\bm\alpha+\mathcal{O} \left(|\delta\bm\alpha|^2+|\delta\omega|^2\right).
$$
Hence there exists $c>0$ such that these three quantities are nonzero whenever $|\delta\bm\alpha|+|\delta\omega|\le c\varepsilon^2$. The remaining denominators stay separated by \eqref{eq:boundary-material-nonresonance}, and the single-layer inverses at the vertex persist in a sufficiently small parameter neighborhood for fixed $\varepsilon$. Since \eqref{eq:dirac-branches-exact} gives $|\delta\omega|\le C_\varepsilon|\delta\bm\alpha|$ in a sufficiently small neighborhood, we may reduce $r_D^\varepsilon$ so that both branches remain in this admissibility neighborhood and $r_D^\varepsilon\le C\varepsilon^2$.
\end{proof}

\section{Splitting of the Dirac point and local gaps when $\kappa_g\neq 0$}\label{sec:5}

We now introduce a small gyrotropic bias and determine how it changes the nonmagnetic Dirac structure obtained in Section~4. The two bulk configurations have identical geometry and material parameters except for the sign of $\kappa_g$, which is taken to be $+\kappa$ or $-\kappa$. Rotational symmetry still preserves the three sectors, but the nonmagnetic relation between sectors $0$ and $1$ is replaced by a correspondence between opposite magnetizations. We first calculate the biased dipole determinants to identify the leading displacement of the two candidate frequencies. As in the nonmagnetic problem, these truncated equations must then be justified at the level of the full boundary operator. For this purpose, we eliminate the interior density using the explicit Fourier--Bessel solution in each disk and separate the three resonant terms of the exterior Green function from the remaining boundary equation. Inverting this remaining equation gives an exact $3\times3$ system for the resonant Bloch amplitudes. Its sector decomposition at $\bm\alpha^*$ yields the actual split frequencies and excludes the third sector from the relevant frequency window. Inversion--conjugation then proves that reversing the magnetization exchanges the two frequencies exactly. Finally, for a fixed sufficiently small nonzero bias, we choose a closed interval strictly between the split frequencies and use continuity and a uniform inverse bound to extend its spectral exclusion to nearby Bloch vectors. The resulting gap is common to both magnetization directions and is local near $\bm\alpha^*$.

\subsection{Candidate frequencies}

As in the nonmagnetic case, we fix $\bm\alpha=\bm\alpha^*$. Following Proposition~\ref{prop:leading-candidate-determinant}, we begin by calculating the candidate frequencies. Using the self entries $(\widehat{\mathcal{S}}_{\omega_s}^{pp,\varepsilon})_{n,n}$, $(\mathcal{K}_{\omega_s}^{pp,\varepsilon,*})_{n,n}$, and $(\mathcal{T}_{\omega_s}^{pp,\varepsilon})_{n,n}$, write $D_{0,\pm}^{\mathrm{dip},\varepsilon}$ and $D_{1,\pm}^{\mathrm{dip},\varepsilon}$ for the raw two-density self blocks on $\varphi_{B,-1}$ and $\varphi_{A,1}$, respectively, with coefficient order $(c_{1,p,n},c_{2,p,n})^\top$:
$$
D_{0,\pm}^{\mathrm{dip},\varepsilon} =
\begin{pmatrix}
(\widehat{\mathcal{S}}_{\omega_{\mathrm m}}^{pp,\varepsilon})_{-1,-1} & -(\widehat{\mathcal{S}}_{\omega_{\mathrm c}}^{pp,\varepsilon})_{-1,-1} \\
-\mu_{\mathrm m}^{-1} \left(\dfrac12+(\mathcal{K}_{\omega_{\mathrm m}}^{pp,\varepsilon,*})_{-1,-1}\right) & \displaystyle \frac{\mu}{\mu^2-\kappa^2} \left(-\dfrac12+(\mathcal{K}_{\omega_{\mathrm c}}^{pp,\varepsilon,*})_{-1,-1}\right) \pm \ii\frac{\kappa}{\mu^2-\kappa^2} (\mathcal{T}_{\omega_{\mathrm c}}^{pp,\varepsilon})_{-1,-1}
\end{pmatrix},
$$
and
$$
D_{1,\pm}^{\mathrm{dip},\varepsilon} =
\begin{pmatrix}
(\widehat{\mathcal{S}}_{\omega_{\mathrm m}}^{pp,\varepsilon})_{1,1} & -(\widehat{\mathcal{S}}_{\omega_{\mathrm c}}^{pp,\varepsilon})_{1,1} \\
-\mu_{\mathrm m}^{-1} \left(\dfrac12+(\mathcal{K}_{\omega_{\mathrm m}}^{pp,\varepsilon,*})_{1,1}\right) & \displaystyle \frac{\mu}{\mu^2-\kappa^2} \left(-\dfrac12+(\mathcal{K}_{\omega_{\mathrm c}}^{pp,\varepsilon,*})_{1,1}\right) \pm \ii\frac{\kappa}{\mu^2-\kappa^2} (\mathcal{T}_{\omega_{\mathrm c}}^{pp,\varepsilon})_{1,1}
\end{pmatrix}.
$$

Section~\ref{sec:4} places the Dirac frequency on the $\Delta_{\mathrm m}$ branch, so $\Delta_{\mathrm c}$ branch will not be considered here. Write $\Delta_{\mathrm m}=\varepsilon^2z$, with $z$ near $z_*$, and assume $|\Delta_{\mathrm c}|\ge c_{\mathrm c}>0$. Set
\begin{equation}\label{eq:gyro-main-block-function}
f_{j,\pm}^\varepsilon(\omega,\kappa):=\det D_{j,\pm}^{\mathrm{dip},\varepsilon}, \qquad  j=0,1.
\end{equation}
The candidate frequencies are the roots of $f_{j,\pm}^\varepsilon(\omega,\kappa)$.

\begin{proposition}[Principal characteristic equations for opposite magnetizations] \label{prop:gyro-main-equations}
Let $\Delta_{\mathrm m}=\varepsilon^2z$, with $z$ in a fixed neighborhood of $z_*$, and assume $|\Delta_{\mathrm c}|\ge c_{\mathrm c}>0$, $|\kappa|<\mu$. The two dipole characteristic functions in~\eqref{eq:gyro-main-block-function} satisfy
$$
f_{0,\pm}^\varepsilon(\omega,\kappa)=\left( -\frac12+ \frac{2\pi Q^2\varepsilon^2} {3|Y|\Delta_{\mathrm c}} \right) \Bigg[ \frac{1}{\mu\pm\kappa} \left( -\frac12+ \frac{2\pi Q^2\varepsilon^2} {3|Y|\Delta_{\mathrm m}} \right)- \mu_{\mathrm m}^{-1} \left( \frac12+ \frac{2\pi Q^2\varepsilon^2} {3|Y|\Delta_{\mathrm m}} \right) \Bigg]+ \mathcal{O}(\varepsilon^2|\ln\varepsilon|),
$$
and
$$
f_{1,\pm}^\varepsilon(\omega,\kappa) =\left( -\frac12+ \frac{2\pi Q^2\varepsilon^2} {3|Y|\Delta_{\mathrm c}} \right) \Bigg[ \frac{1}{\mu\mp\kappa} \left( -\frac12+ \frac{2\pi Q^2\varepsilon^2} {3|Y|\Delta_{\mathrm m}} \right)- \mu_{\mathrm m}^{-1} \left( \frac12+ \frac{2\pi Q^2\varepsilon^2} {3|Y|\Delta_{\mathrm m}} \right) \Bigg] + \mathcal{O}(\varepsilon^2|\ln\varepsilon|).
$$
Consequently, opposite magnetizations satisfy
$$
f_{0,+}^\varepsilon(\omega,\kappa) = f_{1,-}^\varepsilon(\omega,\kappa),  \qquad f_{1,+}^\varepsilon(\omega,\kappa) = f_{0,-}^\varepsilon(\omega,\kappa).
$$
\end{proposition}

\begin{proof}
Section~\ref{sec:nine-fourier-elements} gives
$$
(\widehat{\mathcal{S}}_{\omega_s}^{pp,\varepsilon})_{\pm1,\pm1} =-\frac12+\frac{2\pi Q^2\varepsilon^2}{3|Y|\Delta_s} +\mathcal{O}(\varepsilon^2|\ln\varepsilon|),\qquad  (\mathcal{K}_{\omega_s}^{pp,\varepsilon,*})_{\pm1,\pm1} =\frac{2\pi Q^2\varepsilon^2}{3|Y|\Delta_s} +\mathcal{O}(\varepsilon^2|\ln\varepsilon|),
$$
where $s\in\{\mathrm m,\mathrm c\}$, and
\begin{align*}
(\mathcal{T}_{\omega_{\mathrm c}}^{pp,\varepsilon})_{-1,-1} &=\ii\left(-\frac12+\frac{2\pi Q^2\varepsilon^2}{3|Y|\Delta_{\mathrm c}}\right) +\mathcal{O}(\varepsilon^2|\ln\varepsilon|),\\
(\mathcal{T}_{\omega_{\mathrm c}}^{pp,\varepsilon})_{1,1} &=-\ii\left(-\frac12+\frac{2\pi Q^2\varepsilon^2}{3|Y|\Delta_{\mathrm c}}\right) +\mathcal{O}(\varepsilon^2|\ln\varepsilon|).
\end{align*}
Hence
$$
\frac{\mu}{\mu^2-\kappa^2} \left(-\frac12+\frac{2\pi Q^2\varepsilon^2}{3|Y|\Delta_{\mathrm c}}\right) \pm \ii\frac{\kappa}{\mu^2-\kappa^2}(\mathcal{T}_{\omega_{\mathrm c}}^{pp,\varepsilon})_{-1,-1} = \frac1{\mu\pm\kappa} \left(-\frac12+\frac{2\pi Q^2\varepsilon^2}{3|Y|\Delta_{\mathrm c}}\right) +\mathcal{O}(\varepsilon^2|\ln\varepsilon|).
$$
Furthermore, the exact self-diagonal single-layer and normal entries are even in $n$, while $(\mathcal{T}_{\omega_s}^{pp,\varepsilon})_{-n,-n}=-(\mathcal{T}_{\omega_s}^{pp,\varepsilon})_{n,n}$. Since $\omega_{\mathrm c}$ is even in the bias, these identities give the exact exchange under opposite magnetizations. Substitution into $\det D_{j,\pm}^{\mathrm{dip},\varepsilon}$ gives the two asymptotic formulas.
\end{proof}

\begin{proposition}
Let $\omega_{0,m}>0$ satisfy
$$
\omega_{\mathrm m}(\omega_{0,m})=|\bm\alpha^*|, \qquad  \omega_{\mathrm m}'(\omega_{0,m})\neq0,
$$
and impose the admissibility assumptions of Proposition~\ref{prop:gyro-main-equations}. Assume $\mu\ne\mu_{\mathrm m}$ and restrict $|\kappa|<\kappa_0$, with $\kappa_0<\min\{\mu/2,|\mu_{\mathrm m}-\mu|/2\}$ small enough that the following scaled roots stay in the prescribed neighborhood of $z_*$. In particular their numerators and material denominators are separated from zero. The sector-$0$ and sector-$1$ principal candidate roots satisfy
$$
\Delta_{\mathrm m,0}^{\pm} = \frac{4\pi Q^2}{3|Y|} \frac{\mu_{\mathrm m}-\mu\mp\kappa} {\mu_{\mathrm m}+\mu\pm\kappa} \varepsilon^2 + \mathcal{O}(\varepsilon^4|\ln\varepsilon|),
$$
and
$$
\Delta_{\mathrm m,1}^{\pm} = \frac{4\pi Q^2}{3|Y|} \frac{\mu_{\mathrm m}-\mu\pm\kappa} {\mu_{\mathrm m}+\mu\mp\kappa} \varepsilon^2 + \mathcal{O}(\varepsilon^4|\ln\varepsilon|).
$$
\end{proposition}
\begin{proof}
For sector $0$, the second factor in Proposition~\ref{prop:gyro-main-equations} gives
$$
(\mu_{\mathrm m}-\mu\mp\kappa) \frac{2\pi Q^2\varepsilon^2}{3|Y|\Delta_{\mathrm m}} -\frac{\mu_{\mathrm m}+\mu\pm\kappa}{2} =\mathcal{O}(\varepsilon^2|\ln\varepsilon|),
$$
whose leading zero has a nonzero derivative with respect to $z$, uniformly under the stated restrictions. On a small complex disk around that zero the remainder is analytic and $\mathcal{O}(\varepsilon^2|\ln\varepsilon|)$. Rouch\'e's theorem gives one simple root and a displacement of this order in $z$, hence $\mathcal{O}(\varepsilon^4|\ln\varepsilon|)$ in $\Delta_{\mathrm m}$. Sector $1$ replaces $\mu\pm\kappa$ by $\mu\mp\kappa$, which yields the desired result.
\end{proof}

\subsection{Calculation of resonant frequencies}

In this subsection we calculate the resonant frequencies based on the estimated candidates in previous discussion. To justify the magnetic candidates without discarding their coupling to other modes, we reduce the full boundary equation to three resonant Bloch amplitudes. A regular separated solution of the Helmholtz equation in a disk, with boundary Fourier mode $n$, is
$$
\frac{J_{|n|}(\omega_c r)}{J_{|n|}(\omega_c\varepsilon)}\frac{e^{\ii n\theta}}{\sqrt{2\pi}}.
$$
Let $h_{p,n}:=\langle\varphi_{p,n},\mathcal{S}_{\omega_m}^{\bm\alpha,\varepsilon}\phi\rangle_{p,\varepsilon}$. Then the solution of \eqref{eq:system} in $D_p^\varepsilon$ is
$$
u|_{D_p^\varepsilon}(r,\theta) = \sum\limits_{n\in\mathbb{Z}}h_{p,n}\frac{J_{|n|}(\omega_c r)}{J_{|n|}(\omega_c\varepsilon)}\frac{e^{\ii n\theta}}{\sqrt{2\pi}}.
$$
Hence the normal derivative on $\partial D_p^\varepsilon$ is
$$
\partial_{\bm\nu}u|_- = \sum\limits_{n\in\mathbb{Z}}h_{p,n}\omega_c\frac{J_{|n|}'(\omega_c\varepsilon)}{J_{|n|}(\omega_c\varepsilon)}\varphi_{p,n}.
$$
Eliminate $\psi$ from \eqref{eq:bie-continuity}:
$$
\psi=(\mathcal{S}_{\omega_{\mathrm c}}^{\bm\alpha,\varepsilon})^{-1} \mathcal{S}_{\omega_{\mathrm m}}^{\bm\alpha,\varepsilon}\phi.
$$
We then have
\begin{equation}\label{eq:normal-derivative}
\left\langle\varphi_{p,n},\left(-\frac{1}{2}I + (\mathcal{K}_{\omega_c}^{-\bm\alpha,\varepsilon})^*\right)\psi\right\rangle_{p,\varepsilon} = \omega_c\frac{J_{|n|}'(\omega_c\varepsilon)}{J_{|n|}(\omega_c\varepsilon)}\langle\varphi_{p,n},\mathcal{S}_{\omega_m}^{\bm\alpha,\varepsilon}\phi\rangle_{p,\varepsilon}.
\end{equation}
Similarly, for tangential derivatives,
$$
\partial_{\bm\tau}\varphi_{p,n} = -\frac{\ii n}{\varepsilon}\varphi_{p,n},
$$
so
\begin{equation}\label{eq:tangential-derivative}
\langle\varphi_{p,n},\mathcal{T}_{\omega_c}^{\bm\alpha,\varepsilon}\psi\rangle_{p,\varepsilon} = -\frac{\ii n}{\varepsilon}\langle\varphi_{p,n},\mathcal{S}_{\omega_m}^{\bm\alpha,\varepsilon}\phi\rangle_{p,\varepsilon}.
\end{equation}

Putting \eqref{eq:normal-derivative} and \eqref{eq:tangential-derivative} into \eqref{eq:bie-flux}, we have for $p=A,B$ and $n\in\mathbb Z$,
\begin{equation}\label{eq:one-density-boundary-operator}
\left\langle\varphi_{p,n}, \left(\frac12 I+(\mathcal{K}_{\omega_{\mathrm m}}^{-\bm\alpha,\varepsilon})^*\right) \phi\right\rangle_{p,\varepsilon}- \frac{\mu_{\mathrm m}}{\mu^2-\kappa_g^2} \left(\mu\omega_{\mathrm c} \frac{J_{|n|}'(\omega_{\mathrm c}\varepsilon)} {J_{|n|}(\omega_{\mathrm c}\varepsilon)} +\frac{\kappa_g n}{\varepsilon}\right) \left\langle\varphi_{p,n}, \mathcal{S}_{\omega_{\mathrm m}}^{\bm\alpha,\varepsilon}\phi \right\rangle_{p,\varepsilon}=0.
\end{equation}

\begin{lemma}\label{lem:boundary-three-pole-reduction}
Fix positive material parameters and a sufficiently small complex disk of $\kappa_g$ about zero. Assume $|\Delta_m(\omega)|\leq C_1\varepsilon^2$ and $|\bm\alpha-\bm\alpha^*|\leq C_2\varepsilon^2$ for fixed constants $C_1,C_2>0$. For all sufficiently small $\varepsilon$, the following equations uniquely determine three densities $\phi^{(\ell)}\in H^{-1/2}(\partial D^\varepsilon)$, $\ell=1,2,3$, $p=A,B$, $n\in\mathbb Z$:
\begin{align}
&\left\langle\varphi_{p,n}, \left(\frac12 I+\mathcal{K}_{H_0,\omega_{\mathrm m}}^{\varepsilon,*} +\widetilde{\mathcal{K}}_{\omega_{\mathrm m}}^{\varepsilon,*}\right) \phi^{(\ell)} -\frac{\partial_{\bm\nu}e^{\ii(\bm\alpha+\bm q_\ell)\cdot\bm x}} {\sqrt{|Y|}}\right\rangle_{p,\varepsilon} \notag\\
&\quad =\frac{\mu_{\mathrm m}}{\mu^2-\kappa_g^2} \left(\mu\omega_{\mathrm c} \frac{J_{|n|}'(\omega_{\mathrm c}\varepsilon)} {J_{|n|}(\omega_{\mathrm c}\varepsilon)} +\frac{\kappa_g n}{\varepsilon}\right) \left\langle\varphi_{p,n}, (\mathcal{S}_{H_0,\omega_{\mathrm m}}^\varepsilon+ \widetilde{\mathcal{S}}_{\omega_{\mathrm m}}^\varepsilon) \phi^{(\ell)} -\frac{e^{\ii(\bm\alpha+\bm q_\ell)\cdot\bm x}}{\sqrt{|Y|}} \right\rangle_{p,\varepsilon}. \label{eq:boundary-auxiliary-density}
\end{align}
The $3\times3$ matrix $\mathsf M^\varepsilon=(\mathsf M^\varepsilon_{\ell,r})$, defined with $\phi^{(r)}$, $r=1,2,3$ is also analytic:
\begin{equation}\label{eq:boundary-three-pole-matrix}
\mathsf M^\varepsilon_{\ell,r}(\bm\alpha,\omega,\kappa_g):=\bigl(\omega_{\mathrm m}(\omega)^2 -(\bm\alpha+\bm q_\ell)\cdot(\bm\alpha+\bm q_\ell)\bigr) \delta_{\ell r} +\frac1{\sqrt{|Y|}}\int_{\partial D^\varepsilon} e^{-\ii(\bm\alpha+\bm q_\ell)\cdot\bm x} \phi^{(r)}(\bm x) d\sigma(\bm x),  \ell,r=1,2,3.
\end{equation}
The nullspace of $\mathsf M^\varepsilon$ corresponds bijectively to that of \eqref{eq:one-density-boundary-operator} with the same multiplicities. Furthermore, these are also the nullity and multiplicities of $\mathbb A^\varepsilon$.
\end{lemma}

\begin{proof}
Write $\phi^{(\ell)} = \sum_{q=A,B}\sum_{n\in\mathbb{Z}}c_{q,n}^{(\ell)}\varphi_{q,n}$ and let
\begin{align*}
&\left( \mathcal{K}_{H_0,\omega_{\mathrm m}}^{\varepsilon,*} +\widetilde{\mathcal{K}}_{\omega_{\mathrm m}}^{\varepsilon,*} \right)_{p,m;q,n}:= \left\langle\varphi_{p,m}, \left( \mathcal{K}_{H_0,\omega_{\mathrm m}}^{\varepsilon,*} +\widetilde{\mathcal{K}}_{\omega_{\mathrm m}}^{\varepsilon,*} \right)\varphi_{q,n} \right\rangle_{p,\varepsilon}, \\
&\left( \mathcal{S}_{H_0,\omega_{\mathrm m}}^\varepsilon +\widetilde{\mathcal{S}}_{\omega_{\mathrm m}}^\varepsilon \right)_{p,m;q,n}:= \left\langle\varphi_{p,m}, \left( \mathcal{S}_{H_0,\omega_{\mathrm m}}^\varepsilon +\widetilde{\mathcal{S}}_{\omega_{\mathrm m}}^\varepsilon \right)\varphi_{q,n} \right\rangle_{p,\varepsilon}.
\end{align*}
By rearranging terms in \eqref{eq:boundary-auxiliary-density} we have
$$
\sum\limits_{q=A,B}\sum\limits_{n\in\mathbb{Z}}\mathcal L_{p,m;q,n}^\varepsilon c_{q,n}^{(\ell)} = \left\langle\varphi_{p,m}, \frac{\partial_{\bm\nu} e^{\ii(\bm\alpha+\bm q_\ell)\cdot\bm x}}{\sqrt{|Y|}} \right\rangle_{p,\varepsilon}-b_m^\varepsilon \left\langle\varphi_{p,m}, \frac{e^{\ii(\bm\alpha+\bm q_\ell)\cdot\bm x}}{\sqrt{|Y|}} \right\rangle_{p,\varepsilon},
$$
where
\begin{align*}
\mathcal L^\varepsilon_{p,m;q,n}:=&\frac12\delta_{pq}\delta_{mn}+ \left( \mathcal{K}_{H_0,\omega_{\mathrm m}}^{\varepsilon,*} +\widetilde{\mathcal{K}}_{\omega_{\mathrm m}}^{\varepsilon,*} \right)_{p,m;q,n}-b_m^\varepsilon \left( \mathcal{S}_{H_0,\omega_{\mathrm m}}^\varepsilon +\widetilde{\mathcal{S}}_{\omega_{\mathrm m}}^\varepsilon \right)_{p,m;q,n}, \\
b_m^\varepsilon :=& \frac{\mu_{\mathrm m}}{\mu^2-\kappa_g^2} \left( \mu\omega_{\mathrm c} \frac{J_{|m|}'(\omega_{\mathrm c}\varepsilon)} {J_{|m|}(\omega_{\mathrm c}\varepsilon)} +\frac{\kappa_gm}{\varepsilon} \right).
\end{align*}
The Bessel power series and recurrence used in \eqref{eq:bessel-estimate} give
\begin{align}
\frac{\mu_{\mathrm m}\mu}{\mu^2-\kappa_g^2} \omega_{\mathrm c} \frac{J_0'(\omega_{\mathrm c}\varepsilon)} {J_0(\omega_{\mathrm c}\varepsilon)} &=-\frac{\omega_{\mathrm m}(\omega)^2\varepsilon}{2} +\mathcal{O}(\varepsilon^3), \notag\\
\frac{\mu_{\mathrm m}}{\mu^2-\kappa_g^2} \left(\mu\omega_{\mathrm c} \frac{J_{|n|}'(\omega_{\mathrm c}\varepsilon)} {J_{|n|}(\omega_{\mathrm c}\varepsilon)} +\frac{\kappa_g n}{\varepsilon}\right) &=\frac{\mu_{\mathrm m}}{\mu-\operatorname{sgn}(n)\kappa_g} \frac{|n|}{\varepsilon} +\mathcal{O} \left(\frac{\varepsilon}{1+|n|}\right),\qquad  n\ne0. \label{eq:boundary-local-flux-asymptotics}
\end{align}
Equations~\eqref{eq:self-K-H0-fourier}, \eqref{eq:self-S-H0-zero}, and~\eqref{eq:self-S-H0-nonzero}, together with Lemmas~\ref{lem:self-interaction-smooth} and~\ref{lem:cross-component}, give the following diagonal, mixed, and off-diagonal estimates:
\begin{align*}
\mathcal L_{p,0;p,0}^\varepsilon &= 1+\mathcal{O}(\varepsilon^2(1+|\ln\varepsilon|)), \\
\mathcal L_{p,n;p,n}^\varepsilon &= \frac{1}{2}+\frac{\mu_m}{2(\mu-\operatorname{sgn}(n)\kappa_g)} + \mathcal{O}(\varepsilon^2(1+|\ln\varepsilon|)),  n\ne0, \\
\mathcal L_{p,0;q,n}^\varepsilon &= \mathcal{O}(\varepsilon^3), \qquad \mathcal L_{p,n;q,0}^\varepsilon = \mathcal{O}(\varepsilon), \qquad n\ne0, \\
\mathcal L_{p,m;q,n}^\varepsilon &= \mathcal{O}(\varepsilon^2),  m,n\ne0,  (p,m)\ne(q,n).
\end{align*}
For each $\ell$, order the coefficients $c_{p,n}^{(\ell)}$ as follows: $(c_{A,0}^{(\ell)}, c_{B,0}^{(\ell)}; c_{A,1}^{(\ell)}, c_{A,-1}^{(\ell)}, c_{B,1}^{(\ell)}, c_{B,-1}^{(\ell)}, \cdots)$. The first two coefficients are regarded as one group, while the remaining ones are regarded as another group. In this decomposition, the operator $\mathcal L^\varepsilon$ has the block form
\begin{align*}
&
\begin{pmatrix}
\mathcal L_{0,0}^\varepsilon & \mathcal L_{0,\neq 0}^\varepsilon \\
\mathcal L_{\neq 0,0}^\varepsilon & \mathcal L_{\neq 0,\neq 0}^\varepsilon
\end{pmatrix}
\\
&=
\begin{pmatrix}
I+\mathcal{O} \left(\varepsilon^2(1+|\ln\varepsilon|)\right) &\mathcal{O}(\varepsilon^3)\\
\mathcal{O}(\varepsilon) &\displaystyle \operatorname{diag}_{p=A,B, n\ne0} \frac{1+\mu_{\mathrm m}/ (\mu-\operatorname{sgn}(n)\kappa_g)}{2} +\mathcal{O} \left(\varepsilon^2(1+|\ln\varepsilon|)\right)
\end{pmatrix}.
\end{align*}
These are operator-norm bounds in $h^{-1/2}$. It is immediate that this matrix is invertible with bounded inverse norm, hence \eqref{eq:boundary-auxiliary-density} uniquely determines $\phi^{(\ell)}$, $\ell=1,2,3$.

Now we calculate $\phi^{(\ell)}$. In \eqref{eq:boundary-auxiliary-density}, the right-hand side without terms involving $\phi^{(\ell)}$ is
\begin{equation}\label{eq:appendix-b-1}
f_{p,n}^{(\ell)}:=\left\langle\varphi_{p,n}, \frac{\partial_{\bm\nu}e^{\ii(\bm\alpha+\bm q_\ell)\cdot\bm x}} {\sqrt{|Y|}}\right\rangle_{p,\varepsilon}- \frac{\mu_{\mathrm m}}{\mu^2-\kappa_g^2} \left( \mu\omega_{\mathrm c} \frac{J_{|n|}'(\omega_{\mathrm c}\varepsilon)} {J_{|n|}(\omega_{\mathrm c}\varepsilon)} +\frac{\kappa_g n}{\varepsilon} \right) \left\langle\varphi_{p,n}, \frac{e^{\ii(\bm\alpha+\bm q_\ell)\cdot\bm x}} {\sqrt{|Y|}}\right\rangle_{p,\varepsilon}.
\end{equation}
When $n=0$, a direct calculation gives
\begin{align*}
&\left\langle\varphi_{p,0}, \frac{e^{\ii(\bm\alpha+\bm q_\ell)\cdot\bm x}} {\sqrt{|Y|}}\right\rangle_{p,\varepsilon}=\sqrt{\frac{2\pi}{|Y|}} e^{\ii(\bm\alpha+\bm q_\ell)\cdot\bm z_p} \left[ 1-\frac{\varepsilon^2}{4} (\bm\alpha+\bm q_\ell)\cdot(\bm\alpha+\bm q_\ell) +\mathcal{O}(\varepsilon^4) \right] \\
&\left\langle\varphi_{p,0}, \frac{\partial_{\bm\nu}e^{\ii(\bm\alpha+\bm q_\ell)\cdot\bm x}} {\sqrt{|Y|}}\right\rangle_{p,\varepsilon}= -\sqrt{\frac{2\pi}{|Y|}} e^{\ii(\bm\alpha+\bm q_\ell)\cdot\bm z_p} \frac{\varepsilon}{2} (\bm\alpha+\bm q_\ell)\cdot(\bm\alpha+\bm q_\ell) +\mathcal{O}(\varepsilon^3).
\end{align*}
Putting the two equations above and \eqref{eq:boundary-local-flux-asymptotics} into \eqref{eq:appendix-b-1} yields, when $n=0$,
$$
\sqrt{\frac{2\pi}{|Y|}}  e^{\ii(\bm\alpha+\bm q_\ell)\cdot\bm z_p} \frac{\varepsilon}{2} \bigl(\omega_{\mathrm m}(\omega)^2 -(\bm\alpha+\bm q_\ell)\cdot(\bm\alpha+\bm q_\ell)\bigr) +\mathcal{O}(\varepsilon^3)=\mathcal{O}(\varepsilon^3).
$$
When $n=\pm1$, we have
\begin{align*}
&\left\langle\varphi_{p,n}, \frac{e^{\ii(\bm\alpha+\bm q_\ell)\cdot\bm x}}{\sqrt{|Y|}} \right\rangle_{p,\varepsilon}= \frac{\ii\varepsilon}{2} \sqrt{\frac{2\pi}{|Y|}}  e^{\ii(\bm\alpha+\bm q_\ell)\cdot\bm z_p} ((\bm\alpha+\bm q_\ell)_1-\ii n(\bm\alpha+\bm q_\ell)_2) +\mathcal{O}(\varepsilon^3), \\
&\left\langle\varphi_{p,n}, \frac{\partial_\nu e^{\ii(\bm\alpha+\bm q_\ell)\cdot\bm x}}{\sqrt{|Y|}} \right\rangle_{p,\varepsilon}= \frac{\ii}{2} \sqrt{\frac{2\pi}{|Y|}}  e^{\ii(\bm\alpha+\bm q_\ell)\cdot\bm z_p} ((\bm\alpha+\bm q_\ell)_1-\ii n(\bm\alpha+\bm q_\ell)_2) +\mathcal{O}(\varepsilon^2).
\end{align*}
Then
$$
f_{p,n}^{(\ell)}=\frac{\ii}{2}\sqrt{\frac{2\pi}{|Y|}}  \frac{\mu-n\kappa_g-\mu_{\mathrm m}} {\mu-n\kappa_g}e^{\ii(\bm\alpha+\bm q_\ell)\cdot\bm z_p} \bigl((\bm\alpha+\bm q_\ell)_1 -\ii n(\bm\alpha+\bm q_\ell)_2\bigr)+\mathcal{O}(\varepsilon^2),
$$
which implies $f_{p,n}^{(\ell)}$ is bounded when $n=\pm 1$. Using the block inverse estimate above and the diagonal leading term in~\eqref{eq:boundary-auxiliary-density}, we obtain
\begin{align*}
\langle\varphi_{p,0},\phi^{(\ell)}\rangle_{p,\varepsilon} &=\mathcal{O}(\varepsilon^3),\\
\langle\varphi_{p,n},\phi^{(\ell)}\rangle_{p,\varepsilon} &=-\ii\sqrt{\frac{2\pi}{|Y|}}  \frac{\mu_{\mathrm m}-\mu+n\kappa_g} {\mu_{\mathrm m}+\mu-n\kappa_g} e^{\ii(\bm\alpha+\bm q_\ell)\cdot\bm z_p} \bigl((\bm\alpha+\bm q_\ell)_1 -\ii n(\bm\alpha+\bm q_\ell)_2\bigr)\\
&\qquad +\mathcal{O} \left(\varepsilon^2(1+|\ln\varepsilon|)\right), \qquad  n=\pm1.
\end{align*}

The Jacobi--Anger expansion gives, for every $n\ne0$,
$$
\left|\left\langle\varphi_{p,n}, \frac{e^{\ii\bm k\cdot\bm x}}{\sqrt{|Y|}} \right\rangle_{p,\varepsilon} \right|\le C\frac{(C\varepsilon)^{|n|}}{|n|!}, \qquad  \left|\left\langle\varphi_{p,n}, \frac{\partial_{\bm\nu}e^{\ii\bm k\cdot\bm x}} {\sqrt{|Y|}} \right\rangle_{p,\varepsilon} \right|\le C\frac{(C\varepsilon)^{|n|-1}}{(|n|-1)!}.
$$
Hence in $h^{-1/2}$ norm
$$
\left(\sum\limits_{p=A,B}\sum\limits_{|n|\geq 2} (1+n^2)^{-1/2}|\langle\varphi_{p,n},\phi^{(\ell)}\rangle_{p,\varepsilon}|^2\right)^{1/2}\leq C\varepsilon.
$$
The leading density coefficients occur at $n=\pm1$; the constant and higher Fourier coefficients have the smaller orders stated above.

Now we are ready to calculate \eqref{eq:boundary-three-pole-matrix}. The preceding estimates give
\begin{align*}
&\frac1{\sqrt{|Y|}}\int_{\partial D^\varepsilon} e^{-\ii(\bm\alpha+\bm q_\ell)\cdot\bm x} \phi^{(r)}(\bm x) d\sigma(\bm x)\\
=&-\frac{\pi\varepsilon^2}{|Y|}\left(\sum_{p=A,B} e^{\ii(\bm q_r-\bm q_\ell)\cdot\bm z_p}\right)\times \\
&\qquad \sum_{s=\pm1} \frac{\mu_{\mathrm m}-\mu+s\kappa_g} {\mu_{\mathrm m}+\mu-s\kappa_g} \bigl((\bm\alpha+\bm q_\ell)_1 +\ii s(\bm\alpha+\bm q_\ell)_2\bigr) \bigl((\bm\alpha+\bm q_r)_1 -\ii s(\bm\alpha+\bm q_r)_2\bigr)\\
&\qquad  +\mathcal{O} \left(\varepsilon^4(1+|\ln\varepsilon|)\right).
\end{align*}
Since $\bm\alpha+\bm q_\ell=\bm\xi_\ell+\delta\bm\alpha$, we have
$$
\omega_{\mathrm m}(\omega)^2 -(\bm\alpha+\bm q_\ell)\cdot(\bm\alpha+\bm q_\ell) =\varepsilon^2z-2\bm\xi_\ell\cdot\delta\bm\alpha -\delta\bm\alpha\cdot\delta\bm\alpha.
$$

Finally, reinstate the three finite-rank terms of \eqref{eq:operator-three-part}. For a density $\phi$ satisfying \eqref{eq:one-density-boundary-operator}, put
\begin{equation}\label{eq:boundary-amplitude-reconstruction}
a_\ell= \frac{\displaystyle |Y|^{-1/2}\int_{\partial D^\varepsilon} e^{-\ii(\bm\alpha+\bm q_\ell)\cdot\bm x}\phi(\bm x) d\sigma} {\displaystyle\omega_{\mathrm m}(\omega)^2 -(\bm\alpha+\bm q_\ell)\cdot(\bm\alpha+\bm q_\ell)},
\end{equation}
then from the uniqueness of the solution of \eqref{eq:boundary-auxiliary-density} we have
\begin{equation}\label{eq:phi-decomposition}
\phi=-\sum_{r=1}^3a_r\phi^{(r)}.
\end{equation}
Substitution into the first identity gives the homogeneous equation with matrix \eqref{eq:boundary-three-pole-matrix}. Conversely, a nullvector of $\mathsf M^\varepsilon$ reconstructs a solution by \eqref{eq:phi-decomposition}. A zero reconstructed density would force all $a_\ell=0$ because the denominators $\omega_{\mathrm m}(\omega)^2-(\bm\alpha+\bm q_\ell)\cdot(\bm\alpha+\bm q_\ell)$ are nonzero. The correspondence is therefore bijective.
\end{proof}

At $\bm\alpha=\bm\alpha^*$, use the following ordered orthonormal basis of the three amplitudes in \eqref{eq:boundary-amplitude-reconstruction}:
\begin{equation}\label{eq:boundary-amplitude-basis}
\frac1{\sqrt3}
\begin{pmatrix}
1\\
1\\
1
\end{pmatrix},\qquad  \frac1{\sqrt3}
\begin{pmatrix}
e^{-2\pi \ii/3}\\
e^{2\pi \ii/3}\\
1
\end{pmatrix},  \qquad \frac1{\sqrt3}
\begin{pmatrix}
1\\
e^{2\pi \ii/3}\\
e^{-2\pi \ii/3}
\end{pmatrix}.
\end{equation}
The matrix \eqref{eq:boundary-three-pole-matrix} is exactly diagonal under similarity transformation by \eqref{eq:boundary-amplitude-basis} at $\bm\alpha^*$. After division by $\varepsilon^2$, its three diagonal entries are
\begin{align}
&z-\frac{4\pi Q^2}{3|Y|} \frac{\mu_{\mathrm m}-\mu-\kappa_g} {\mu_{\mathrm m}+\mu+\kappa_g} +\mathcal{O} \left(\varepsilon^2(1+|\ln\varepsilon|)\right), \notag\\
&z-\frac{4\pi Q^2}{3|Y|} \frac{\mu_{\mathrm m}-\mu+\kappa_g} {\mu_{\mathrm m}+\mu-\kappa_g} +\mathcal{O} \left(\varepsilon^2(1+|\ln\varepsilon|)\right), \notag\\
&z-\frac{4\pi Q^2}{3|Y|} \left(\frac{\mu_{\mathrm m}-\mu-\kappa_g} {\mu_{\mathrm m}+\mu+\kappa_g} +\frac{\mu_{\mathrm m}-\mu+\kappa_g} {\mu_{\mathrm m}+\mu-\kappa_g}\right) +\mathcal{O} \left(\varepsilon^2(1+|\ln\varepsilon|)\right). \label{eq:boundary-response-sector-diagonal}
\end{align}
At $\kappa_g=0$, the leading expressions are $z-z_*$, $z-z_*$, and $z-2z_*$, with $z_*$ as in \eqref{eq:z-star-exact}.

Now we calculate the actual resonant frequencies based on $\Delta_{m,0}^\pm$ and $\Delta_{m,1}^\pm$.

\begin{theorem}\label{thm:actual-gyro-branches}
There is a fixed $\kappa_0>0$ such that, for real $|\kappa|<\kappa_0$, $\sigma\in\{+1,-1\}$, and sufficiently small $\varepsilon$, the operator on the left of \eqref{eq:one-density-boundary-operator}, with $\bm\alpha=\bm\alpha^*$ and $\kappa_g=\sigma\kappa$, has one simple characteristic frequency in each of sectors $0$ and $1$ in a neighborhood of $\omega_{0,m}$. Denote them by $\widehat\omega_{0,\sigma}^\varepsilon$ and $\widehat\omega_{1,\sigma}^\varepsilon$. Their $\Delta_{\mathrm m}$ values satisfy
\begin{align}
\widehat\Delta_{\mathrm m,0}^{\sigma} &= \frac{4\pi Q^2}{3|Y|} \frac{\mu_{\mathrm m}-\mu-\sigma\kappa} {\mu_{\mathrm m}+\mu+\sigma\kappa} \varepsilon^2 +\mathcal{O} \left(\varepsilon^3|\ln\varepsilon|^{1/2}\right), \label{eq:actual-gyro-delta0}\\
\widehat\Delta_{\mathrm m,1}^{\sigma} &= \frac{4\pi Q^2}{3|Y|} \frac{\mu_{\mathrm m}-\mu+\sigma\kappa} {\mu_{\mathrm m}+\mu-\sigma\kappa} \varepsilon^2 +\mathcal{O} \left(\varepsilon^3|\ln\varepsilon|^{1/2}\right). \label{eq:actual-gyro-delta1}
\end{align}
After reducing $r_{\mathrm m}$ and $\kappa_0$, sector $2$ has no characteristic frequency in the same window. These are also the transmission frequencies of~\eqref{eq:system}.

Opposite magnetizations satisfy the exact exchange relation
\begin{equation}\label{eq:actual-gyro-swap}
\widehat\omega_{0,+}^\varepsilon(\kappa) =\widehat\omega_{1,-}^\varepsilon(\kappa),  \qquad \widehat\omega_{1,+}^\varepsilon(\kappa) =\widehat\omega_{0,-}^\varepsilon(\kappa).
\end{equation}
\end{theorem}

\begin{proof}
Apply Lemma~\ref{lem:boundary-three-pole-reduction} with $\kappa_g=\sigma\kappa$, $\sigma=\pm 1$ and $\bm\alpha=\bm\alpha^*$. The matrix \eqref{eq:boundary-three-pole-matrix}, divided by $\varepsilon^2$, is diagonal in basis \eqref{eq:boundary-amplitude-basis}. By \eqref{eq:boundary-response-sector-diagonal}, its entries are
$$
k_j^\varepsilon(z,\kappa_g)=z-t_j(\kappa_g)+\rho_j^\varepsilon(z,\kappa_g),  \rho_j^\varepsilon(z,\kappa_g)=\mathcal{O}(\varepsilon^2(1+|\ln\varepsilon|)),
$$
where
$$
t_0(\kappa_g)=\frac{4\pi Q^2}{3|Y|} \frac{\mu_{\mathrm m}-\mu-\kappa_g} {\mu_{\mathrm m}+\mu+\kappa_g},\qquad  t_1(\kappa_g)=\frac{4\pi Q^2}{3|Y|} \frac{\mu_{\mathrm m}-\mu+\kappa_g} {\mu_{\mathrm m}+\mu-\kappa_g},\qquad  t_2=t_0+t_1.
$$
By Cauchy's estimates on slightly smaller fixed complex disks, the first derivatives of $\rho_j^\varepsilon$ with respect to $z$ and $\kappa_g$ have the same bound. At $\kappa_g=0$ these three principal roots are $z_*,z_*,2z_*$. Choose $0<r_{\mathrm m}<|z_*|/4$ and let $\kappa_g$ stay inside a neighborhood of $0$ small enough such that $t_0,t_1$ lie in $B_{r_{\mathrm m}/2}(z_*)$, while $t_2$ remains outside $\overline{B_{r_{\mathrm m}}(z_*)}$.

Rouch\'e's theorem on $B_{r_{\mathrm m}}(z_*)$ gives one zero, counted with multiplicity, of each of $k_0^\varepsilon,k_1^\varepsilon$ but none of $k_2^\varepsilon$. The first two zeros are simple, are analytic in $\kappa_g$, and satisfy
$$
z_j^\varepsilon(\kappa_g) =t_j(\kappa_g)+\mathcal{O}(\varepsilon^2(1+|\ln\varepsilon|)),\qquad  j=0,1.
$$
The exact boundary reconstruction excludes any other transmission frequency in this window. 

We finally prove the exact exchange relation \eqref{eq:actual-gyro-swap} between opposite magnetizations. Recall the antilinear operator $\mathcal J$ defined in~\eqref{eq:J-definition} by
$$
(\mathcal J f)(\bm x) =\overline{f(P\bm x)}, \qquad  P\bm x=2\bm c-\bm x.
$$
Its definition is independent of $\kappa_g$. As established in the proof of Proposition~\ref{prop:conjugate-isospectral}, for real frequencies, $\mathcal J$ commutes with the single-layer and normal-derivative boundary operators. Moreover, since $P$ exchanges the two disks and reverses both their outward normals and tangent vectors, the chain rule gives
$$
\mathcal J\partial_{\bm\nu} =\partial_{\bm\nu}\mathcal J,  \qquad \mathcal J\partial_{\bm\tau} =\partial_{\bm\tau}\mathcal J.
$$
However, $\mathcal J$ is antilinear, so that $\mathcal J(\ii f)=-\ii\mathcal J f$. Consequently, applying $\mathcal J$ to the gyrotropic flux condition in~\eqref{eq:system} changes the coefficient of the tangential derivative from $\ii\kappa_g$ to $-\ii\kappa_g$, while leaving the normal-derivative coefficient unchanged. Thus the flux condition with $\kappa_g=\kappa$ is transformed into that with $\kappa_g=-\kappa$.

Furthermore, $\mathcal J$ exchanges the boundary components and reverses the Fourier index. It follows that applying $\mathcal J$ to the equation indexed by $(p,n)$ with bias $\kappa_g=\kappa$ gives the equation indexed by $(\overline p,-n)$ with bias $\kappa_g=-\kappa$. Therefore, if $\phi$ is a nonzero solution of~\eqref{eq:one-density-boundary-operator} at $(\bm\alpha^*,\omega,\kappa)$, then $\mathcal J\phi$ is a nonzero solution of the same equation at $(\bm\alpha^*,\omega,-\kappa)$. Since $\mathcal J^2=I$, this correspondence is bijective. Furthermore, \eqref{eq:J-mapping} and \eqref{eq:J2-Vj-map} show that $\mathcal J$ exchanges sectors $0$ and $1$. Hence a real characteristic frequency in sector $0$ for $\kappa_g=\kappa$ is also a characteristic frequency in sector $1$ for $\kappa_g=-\kappa$, and conversely.

For real $\kappa$ in the prescribed range, the characteristic frequencies under consideration are real, as follows from the self-adjoint transmission problem with Hermitian positive-definite permeability tensor. Since the preceding argument has established that each of sectors $0$ and $1$ contains exactly one simple characteristic frequency for either magnetization, the above bijection identifies these frequencies uniquely. We conclude that
$$
\widehat\omega_{0,+}^\varepsilon(\kappa) =\widehat\omega_{1,-}^\varepsilon(\kappa),  \qquad \widehat\omega_{1,+}^\varepsilon(\kappa) =\widehat\omega_{0,-}^\varepsilon(\kappa),
$$
which proves~\eqref{eq:actual-gyro-swap}.
\end{proof}

From Theorem \ref{thm:actual-gyro-branches}, the corresponding characteristic frequencies are
\begin{align*}
\widehat\omega_{0,+}^\varepsilon(\kappa) &=\sqrt{\frac{k_D^2+\displaystyle\frac{4\pi Q^2}{3|Y|} \frac{\mu_{\mathrm m}-\mu-\kappa} {\mu_{\mathrm m}+\mu+\kappa}\varepsilon^2 }{\epsilon\mu_{\mathrm m}} }+\mathcal{O}(\varepsilon^3|\ln\varepsilon|^{1/2}), \\
\widehat\omega_{1,+}^\varepsilon(\kappa)&= \sqrt{\frac{k_D^2+\displaystyle\frac{4\pi Q^2}{3|Y|} \frac{\mu_{\mathrm m}-\mu+\kappa} {\mu_{\mathrm m}+\mu-\kappa}\varepsilon^2 }{\epsilon\mu_{\mathrm m}} }+\mathcal{O}(\varepsilon^3|\ln\varepsilon|^{1/2}).
\end{align*}

\subsection{A local gap near the Dirac point}

The preceding subsection separates the two characteristic frequencies at the vertex, but this alone does not exclude nearby Bloch eigenfrequencies between them. We now establish that exclusion on an interval strictly inside the vertex splitting. We suppress the dependence on $\kappa$ in the following notation. For $\sigma\in\{+1,-1\}$, define
$$
\bar\omega_\sigma^\varepsilon :=\frac{\widehat\omega_{0,\sigma}^\varepsilon+ \widehat\omega_{1,\sigma}^\varepsilon}{2}, \qquad  m_\sigma^\varepsilon :=\frac{\widehat\omega_{1,\sigma}^\varepsilon- \widehat\omega_{0,\sigma}^\varepsilon}{2}.
$$

\begin{theorem}
Fix $0<|\kappa|<\kappa_0$ as in Theorem~\ref{thm:actual-gyro-branches}. For any $c\in(0,1)$ and every sufficiently small fixed $\varepsilon$, there is $r_{\varepsilon,\kappa}>0$ such that
$$
I_\kappa^\varepsilon :=\left(\bar\omega_+^\varepsilon-c|m_+^\varepsilon|, \bar\omega_+^\varepsilon+c|m_+^\varepsilon|\right)
$$
contains no characteristic frequency of the operator on the left of \eqref{eq:one-density-boundary-operator}, with $\kappa_g=\sigma\kappa$, for either $\sigma=\pm1$, whenever $|\bm\alpha-\bm\alpha^*|<r_{\varepsilon,\kappa}$. The same interval contains no transmission Bloch eigenfrequency. The radius may be reduced so that $r_{\varepsilon,\kappa}\le C\varepsilon^2$.
\end{theorem}

\begin{proof}
Let $J_\kappa^\varepsilon=\overline{I_\kappa^\varepsilon}$. At $\bm\alpha^*$ this compact interval lies strictly between the two simple frequencies for each magnetization. By Theorem~\ref{thm:actual-gyro-branches}, the matrix \eqref{eq:boundary-three-pole-matrix} is uniformly invertible on $J_\kappa^\varepsilon$. Continuity of its inverse gives
\begin{equation}\label{eq:gap-uniform-inverse-bound}
\sup_{\substack{\sigma=\pm1\\
\omega\in J_\kappa^\varepsilon}} \left\|\mathsf M^\varepsilon(\bm\alpha^*, \omega,\sigma\kappa)^{-1}\right\|_{\mathbb C^3\to\mathbb C^3}<\infty. 
\end{equation}

On $J_\kappa^\varepsilon$, $\Delta_{\mathrm m}/\varepsilon^2$ stays away from zero. The analytic dependence proved in Lemma~\ref{lem:boundary-three-pole-reduction} therefore implies that the entries of $\mathsf M^\varepsilon$ converge to their values at $\bm\alpha^*$ uniformly for $\omega\in J_\kappa^\varepsilon$ and $\sigma=\pm1$ as $\bm\alpha\to\bm\alpha^*$. Choose a common neighborhood in which the norm of this matrix difference, multiplied by the bound in \eqref{eq:gap-uniform-inverse-bound}, is less than one. Set $\delta\mathsf M(\bm\alpha):=\mathsf M^\varepsilon(\bm\alpha,\omega,\sigma\kappa)-\mathsf M^\varepsilon(\bm\alpha^*, \omega,\sigma\kappa)$. Then
$$
\mathsf M^\varepsilon(\bm\alpha,\omega,\sigma\kappa) = \mathsf M^\varepsilon(\bm\alpha^*,\omega,\sigma\kappa)(I+\mathsf M^\varepsilon(\bm\alpha^*,\omega,\sigma\kappa)^{-1}\delta\mathsf M(\bm\alpha)).
$$
Choose $r_{\varepsilon,\kappa}>0$ so that $\|\mathsf M^\varepsilon(\bm\alpha^*,\omega,\sigma\kappa)^{-1}\delta\mathsf M(\bm\alpha)\|<1$. The inverse is then given by the convergent Neumann series
$$
(I+\mathsf M^\varepsilon(\bm\alpha^*,\omega,\sigma\kappa)^{-1}\delta\mathsf M(\bm\alpha))^{-1} = \sum\limits_{n=0}^\infty (-\mathsf M^\varepsilon(\bm\alpha^*,\omega,\sigma\kappa)^{-1}\delta\mathsf M(\bm\alpha))^n.
$$
Hence $\mathsf M^\varepsilon(\bm\alpha,\omega,\sigma\kappa)$ is invertible. The equivalence with the transmission problem then excludes all Bloch eigenfrequencies in the interval.
\end{proof}

\section*{Acknowledgment}
W. Wu was supported by NSFC grant No.~12301539.

\subsection*{Data availability}
No datasets were generated or analysed during the current study.

\begin{appendices}

\section{Small-radius layer-potential expansions and basic matrix elements}
In this appendix, we derive the small-radius expansions of the single-layer, normal-derivative, and tangential-derivative operators and compute the Fourier matrix elements used in the finite-dimensional reductions.

\subsection{Expansions on a single boundary component}

Fix $p\in\{A,B\}$. Throughout this subsection, we use $\mathcal{S}_{\lambda}^{pp,\varepsilon}[\varphi](\bm x)$, $\mathcal{K}_{\lambda}^{pp,\varepsilon,*}[\varphi](\bm x)$, $\mathcal{T}_{\lambda}^{pp,\varepsilon}[\varphi](\bm x)$ to denote the corresponding operators on component $p$, evaluated at $\bm x\in\partial D_p^\varepsilon$, with boundary density $\varphi$ supported on component $p$. For $\bm x,\bm y\in\partial D_p^\varepsilon$, write $\bm x=\bm z_p+\varepsilon\bm r(\theta)$, $\bm y=\bm z_p+\varepsilon\bm r(t)$. Then $\bm x-\bm y=\varepsilon(\bm r(\theta)-\bm r(t))$. We expand each term in \eqref{eq:operator-three-part} with this parametrization.

\begin{lemma}\label{lem:appendix-a-1}
On a single boundary component, $\mathcal{S}_{H_0,\lambda}^{pp,\varepsilon}$ satisfies
\begin{equation}\label{eq:self-S-H0-general}
\mathcal{S}_{H_0,\lambda}^{pp,\varepsilon}[\varphi](\theta) =\frac{\varepsilon}{2\pi}\int_0^{2\pi} \left(\ln\varepsilon+\ln\lambda+\gamma_0 +\ln|\bm r(\theta)-\bm r(t)|\right)\varphi(t) dt +\mathcal{O}(\varepsilon^3|\ln\varepsilon|).
\end{equation}
In particular,
\begin{equation}\label{eq:self-S-H0-zero}
\mathcal{S}_{H_0,\lambda}^{pp,\varepsilon}[1] =\varepsilon(\ln\varepsilon+\ln\lambda+\gamma_0) +\mathcal{O}(\varepsilon^3|\ln\varepsilon|),
\end{equation}
and, for $n\ne 0$,
\begin{equation}\label{eq:self-S-H0-nonzero}
\mathcal{S}_{H_0,\lambda}^{pp,\varepsilon}[e^{\ii n\theta}] =-\frac{\varepsilon}{2|n|}e^{\ii n\theta} +\mathcal{O}(\varepsilon^3|\ln\varepsilon|).
\end{equation}
\end{lemma}

\begin{proof}
Set $\bm d(\theta,t)=\bm r(\theta)-\bm r(t)$. Substitution of $\bm x=\varepsilon\bm d(\theta,t)$ into~\eqref{eq:hankel-local-expansion} gives
$$
H_0(\lambda;\varepsilon\bm d(\theta,t)) =\frac1{2\pi} \left(\ln\varepsilon+\ln\lambda+\gamma_0+\ln|\bm d(\theta,t)|\right) +\mathcal{O}(\varepsilon^2|\ln\varepsilon|).
$$
Since $d\sigma_{\bm y}=\varepsilon dt$,
$$
\mathcal{S}_{H_0,\lambda}^{pp,\varepsilon}[\varphi](\theta) =\varepsilon\int_0^{2\pi} H_0(\lambda;\varepsilon\bm d(\theta,t))\varphi(t) dt,
$$
which proves~\eqref{eq:self-S-H0-general}.

The circular logarithmic kernel has the $L^2$ Fourier expansion
$$
\ln|\bm r(\theta)-\bm r(t)| =-\sum_{k=1}^{\infty}\frac{\cos k(\theta-t)}{k}.
$$
Indeed, $ \cos k(\theta-t) =\frac12e^{\ii k\theta}e^{-\ii kt}+\frac12e^{-\ii k\theta}e^{\ii kt} $ and Fourier orthogonality give
\begin{equation}\label{eq:circle-log-fourier-integral}
\frac1{2\pi}\int_0^{2\pi} \ln|\bm r(\theta)-\bm r(t)|e^{\ii nt} dt =
\begin{cases}
0,&n=0,\\
-\dfrac1{2|n|}e^{\ii n\theta},&n\ne0.
\end{cases}
\end{equation}
Also, $\int_0^{2\pi}e^{\ii nt} dt=2\pi\delta_{n0}$. Taking $\varphi=e^{\ii n\theta}$ in~\eqref{eq:self-S-H0-general} proves~\eqref{eq:self-S-H0-zero} and~\eqref{eq:self-S-H0-nonzero}.
\end{proof}

\begin{lemma}
The self-interaction single layers induced by $G_{\Lambda_0^*}^{\lambda}$ and $\widetilde G^\lambda$ satisfy
\begin{align}
\mathcal{S}_{\Lambda_0^*,\lambda}^{pp,\varepsilon}[\varphi](\theta) &=\frac{\varepsilon}{|Y|\Delta_\lambda}\int_0^{2\pi} \left(3-\frac{Q^2\varepsilon^2}{3}|\bm r(\theta)-\bm r(t)|^2\right) \varphi(t) dt+\mathcal{O} \left(\frac{\varepsilon^4}{|\Delta_\lambda|}\right),\label{eq:self-S-res}\\
\widetilde{\mathcal{S}}_\lambda^{pp,\varepsilon}[\varphi](\theta) &=\varepsilon\widetilde G^\lambda(\bm0)\int_0^{2\pi}\varphi(t) dt+\varepsilon^2\int_0^{2\pi} \nabla_{\bm x}\widetilde G^\lambda(\bm0)\cdot (\bm r(\theta)-\bm r(t))\varphi(t) dt +\mathcal{O}(\varepsilon^3).\label{eq:self-S-smooth}
\end{align}
\end{lemma}

\begin{proof}
Again set $\bm d(\theta,t)=\bm r(\theta)-\bm r(t)$. Equation~\eqref{eq:resonant-green-expansion} gives
$$
G_{\Lambda_0^*}^{\lambda}(\varepsilon\bm d) =\frac1{|Y|\Delta_\lambda} \left(3-\frac{Q^2\varepsilon^2}{3}|\bm d|^2\right) +\mathcal{O} \left(\frac{\varepsilon^3}{|\Delta_\lambda|}\right).
$$
Parametrizing the circle yields
\begin{align*}
\mathcal{S}_{\Lambda_0^*,\lambda}^{pp,\varepsilon}[\varphi](\theta) &=\varepsilon\int_0^{2\pi} G_{\Lambda_0^*}^{\lambda}(\varepsilon\bm d(\theta,t))\varphi(t) dt\\
&=\frac{\varepsilon}{|Y|\Delta_\lambda}\int_0^{2\pi} \left(3-\frac{Q^2\varepsilon^2}{3}|\bm d(\theta,t)|^2\right)\varphi(t) dt +\mathcal{O} \left(\frac{\varepsilon^4}{|\Delta_\lambda|}\right).
\end{align*}

Since $\widetilde G^\lambda$ is twice continuously differentiable near $\bm0$,
$$
\widetilde G^\lambda(\varepsilon\bm d) =\widetilde G^\lambda(\bm0) +\varepsilon\nabla_{\bm x}\widetilde G^\lambda(\bm0)\cdot\bm d +\mathcal{O}(\varepsilon^2).
$$
Parametrization gives
\begin{align*}
\widetilde{\mathcal{S}}_\lambda^{pp,\varepsilon}[\varphi](\theta) &=\varepsilon\int_0^{2\pi} \widetilde G^\lambda(\varepsilon\bm d(\theta,t))\varphi(t) dt\\
&=\varepsilon\widetilde G^\lambda(\bm0)\int_0^{2\pi}\varphi(t) dt +\varepsilon^2\int_0^{2\pi}\nabla_{\bm x}\widetilde G^\lambda(\bm0)\cdot\bm d(\theta,t)\varphi(t) dt +\mathcal{O}(\varepsilon^3).
\end{align*}
\end{proof}

\begin{lemma}\label{lem:appendix-a-3}
The self-interaction operator $\mathcal{K}_{H_0,\lambda}^{pp,\varepsilon,*}$ satisfies
\begin{equation}\label{eq:self-K-H0-general}
\mathcal{K}_{H_0,\lambda}^{pp,\varepsilon,*}[\varphi](\theta) =\frac1{4\pi}\int_0^{2\pi}\varphi(t) dt+\mathcal{O}(\varepsilon^2|\ln\varepsilon|).
\end{equation}
Consequently,
\begin{equation}\label{eq:self-K-H0-fourier}
\mathcal{K}_{H_0,\lambda}^{pp,\varepsilon,*}[e^{\ii n\theta}] =
\begin{cases}
\dfrac12+\mathcal{O}(\varepsilon^2|\ln\varepsilon|),&n=0,\\
\mathcal{O}(\varepsilon^2|\ln\varepsilon|),&n\ne0.
\end{cases}
\end{equation}
\end{lemma}

\begin{proof}
Set $\bm d(\theta,t)=\bm r(\theta)-\bm r(t)$. By~\eqref{eq:hankel-gradient-expansion},
$$
\nabla_{\bm x}H_0(\lambda;\varepsilon\bm d) =\frac1{2\pi\varepsilon}\frac{\bm d}{|\bm d|^2} +\mathcal{O}(\varepsilon|\ln\varepsilon|).
$$
Hence
$$
\mathcal{K}_{H_0,\lambda}^{pp,\varepsilon,*}[\varphi](\theta) =\frac1{2\pi}\operatorname{p.v.}\int_0^{2\pi} \frac{\bm\nu(\theta)\cdot\bm d(\theta,t)}{|\bm d(\theta,t)|^2} \varphi(t) dt+\mathcal{O}(\varepsilon^2|\ln\varepsilon|).
$$
On the circle, $\bm\nu(\theta)=\bm r(\theta)$ and
$$
\bm\nu(\theta)\cdot\bm d(\theta,t)=1-\cos(\theta-t), \qquad  |\bm d(\theta,t)|^2=2(1-\cos(\theta-t)).
$$
The principal-value kernel reduces to the constant $1/2$. Therefore
$$
\mathcal{K}_{H_0,\lambda}^{pp,\varepsilon,*}[\varphi](\theta) =\frac1{4\pi}\int_0^{2\pi}\varphi(t) dt +\mathcal{O}(\varepsilon^2|\ln\varepsilon|).
$$
The desired results follow because $\frac1{4\pi}\int_0^{2\pi}e^{\ii nt} dt=\frac12\delta_{n0}$.
\end{proof}

\begin{lemma}
The Neumann--Poincar\'e operator induced by $G_{\Lambda_0^*}^{\lambda}$ satisfies
\begin{equation}\label{eq:self-K-res-general}
\mathcal{K}_{\Lambda_0^*,\lambda}^{pp,\varepsilon,*}[\varphi](\theta) =-\frac{2Q^2\varepsilon^2}{3|Y|\Delta_\lambda} \int_0^{2\pi}(1-\cos(\theta-t))\varphi(t) dt +\mathcal{O} \left(\frac{\varepsilon^3}{|\Delta_\lambda|}\right).
\end{equation}
In particular,
$$
\mathcal{K}_{\Lambda_0^*,\lambda}^{pp,\varepsilon,*}[e^{\ii n\theta}] =
\begin{cases}
-\dfrac{4\pi Q^2\varepsilon^2}{3|Y|\Delta_\lambda} +\mathcal{O}(\varepsilon^3/|\Delta_\lambda|),&n=0,\\
\dfrac{2\pi Q^2\varepsilon^2}{3|Y|\Delta_\lambda}e^{\ii n\theta} +\mathcal{O}(\varepsilon^3/|\Delta_\lambda|),&n=\pm1,\\
\mathcal{O}(\varepsilon^3/|\Delta_\lambda|),&|n|>1.
\end{cases}
$$
\end{lemma}

\begin{proof}
Set $\bm d(\theta,t)=\bm r(\theta)-\bm r(t)$. Equation~\eqref{eq:resonant-green-expansion} gives
$$
\nabla_{\bm x}G_{\Lambda_0^*}^{\lambda}(\varepsilon\bm d) =-\frac{2Q^2\varepsilon}{3|Y|\Delta_\lambda}\bm d +\mathcal{O}\left(\frac{\varepsilon^2}{|\Delta_\lambda|}\right).
$$
Thus
\begin{align*}
\mathcal{K}_{\Lambda_0^*,\lambda}^{pp,\varepsilon,*}[\varphi](\theta) &=\varepsilon\int_0^{2\pi} \bm\nu(\theta)\cdot \nabla_{\bm x}G_{\Lambda_0^*}^{\lambda}(\varepsilon\bm d(\theta,t)) \varphi(t) dt\\
&=-\frac{2Q^2\varepsilon^2}{3|Y|\Delta_\lambda} \int_0^{2\pi}(1-\cos(\theta-t))\varphi(t) dt +\mathcal{O} \left(\frac{\varepsilon^3}{|\Delta_\lambda|}\right).
\end{align*}
Fourier orthogonality gives
$$
\int_0^{2\pi}(1-\cos(\theta-t))e^{\ii nt} dt =
\begin{cases}
2\pi,&n=0,\\
-\pi e^{\ii\theta},&n=1,\\
-\pi e^{-\ii\theta},&n=-1,\\
0,&|n|>1.
\end{cases}
$$
Substituting these four Fourier integrals into~\eqref{eq:self-K-res-general} yields
$$
\mathcal{K}_{\Lambda_0^*,\lambda}^{pp,\varepsilon,*}[e^{\ii n\theta}] =
\begin{cases}
-\dfrac{4\pi Q^2\varepsilon^2}{3|Y|\Delta_\lambda} +\mathcal{O}(\varepsilon^3/|\Delta_\lambda|),&n=0,\\
\dfrac{2\pi Q^2\varepsilon^2}{3|Y|\Delta_\lambda}e^{\ii n\theta} +\mathcal{O}(\varepsilon^3/|\Delta_\lambda|),&n=\pm1,\\
\mathcal{O}(\varepsilon^3/|\Delta_\lambda|),&|n|>1.
\end{cases}
$$
\end{proof}

\begin{lemma}
The Neumann--Poincar\'e operator induced by $\widetilde G^\lambda$ satisfies
\begin{equation}\label{eq:self-K-smooth}
\widetilde{\mathcal{K}}_\lambda^{pp,\varepsilon,*}[\varphi](\theta) =\varepsilon\bm\nu(\theta)\cdot\nabla_{\bm x}\widetilde G^\lambda(\bm0) \int_0^{2\pi}\varphi(t) dt+\mathcal{O}(\varepsilon^2).
\end{equation}
For $n\ne0$, $\widetilde{\mathcal{K}}_\lambda^{pp,\varepsilon,*}[e^{\ii n\theta}]=\mathcal{O}(\varepsilon^2)$.
\end{lemma}

\begin{proof}
Let $\widetilde{\bm g}_\lambda=\nabla_{\bm x}\widetilde G^\lambda(\bm0)$. Smoothness gives
$$
\nabla_{\bm x}\widetilde G^\lambda (\varepsilon(\bm r(\theta)-\bm r(t))) =\widetilde{\bm g}_\lambda+\mathcal{O}(\varepsilon).
$$
Therefore
\begin{align*}
\widetilde{\mathcal{K}}_\lambda^{pp,\varepsilon,*}[\varphi](\theta) &=\varepsilon\int_0^{2\pi} \bm\nu(\theta)\cdot \nabla_{\bm x}\widetilde G^\lambda (\varepsilon(\bm r(\theta)-\bm r(t)))\varphi(t) dt\\
&=\varepsilon\bm\nu(\theta)\cdot\widetilde{\bm g}_\lambda \int_0^{2\pi}\varphi(t) dt+\mathcal{O}(\varepsilon^2).
\end{align*}
For $n\ne0$, $\int_0^{2\pi}e^{\ii nt} dt=0$, so $\widetilde{\mathcal{K}}_\lambda^{pp,\varepsilon,*}[e^{\ii n\theta}]=\mathcal{O}(\varepsilon^2)$.
\end{proof}

\begin{lemma}
The tangential operator $\mathcal{T}_{H_0,\lambda}^{pp,\varepsilon}$ satisfies
\begin{equation}\label{eq:self-T-H0-general}
\mathcal{T}_{H_0,\lambda}^{pp,\varepsilon}[\varphi](\theta) =-\frac1{4\pi}\operatorname{p.v.}\int_0^{2\pi} \cot\frac{\theta-t}{2} \varphi(t) dt +\mathcal{O}(\varepsilon^2|\ln\varepsilon|).
\end{equation}
For $n\ne0$,
$$
\mathcal{T}_{H_0,\lambda}^{pp,\varepsilon}[e^{\ii n\theta}] =\frac{\ii}{2}\operatorname{sgn}(n)e^{\ii n\theta}+\mathcal{O}(\varepsilon^2|\ln\varepsilon|),
$$
and $\mathcal{T}_{H_0,\lambda}^{pp,\varepsilon}[1]=\mathcal{O}(\varepsilon^2|\ln\varepsilon|)$.
\end{lemma}

\begin{proof}
Set $\bm d(\theta,t)=\bm r(\theta)-\bm r(t)$. By~\eqref{eq:hankel-gradient-expansion},
$$
\mathcal{T}_{H_0,\lambda}^{pp,\varepsilon}[\varphi](\theta) =\frac1{2\pi}\operatorname{p.v.}\int_0^{2\pi} \frac{\bm\tau(\theta)\cdot\bm d(\theta,t)}{|\bm d(\theta,t)|^2} \varphi(t) dt+\mathcal{O}(\varepsilon^2|\ln\varepsilon|).
$$
With the convention $\bm\tau(\theta)=(\sin\theta,-\cos\theta)^\top$,
$$
\bm\tau(\theta)\cdot\bm d(\theta,t)=-\sin(\theta-t), \qquad  |\bm d(\theta,t)|^2=2(1-\cos(\theta-t)).
$$
Thus
$$
\frac{\bm\tau(\theta)\cdot\bm d(\theta,t)}{|\bm d(\theta,t)|^2} =-\frac12\cot\frac{\theta-t}{2}.
$$
Consequently,
$$
\mathcal{T}_{H_0,\lambda}^{pp,\varepsilon}[\varphi](\theta) =-\frac1{4\pi}\operatorname{p.v.}\int_0^{2\pi} \cot\frac{\theta-t}{2} \varphi(t) dt +\mathcal{O}(\varepsilon^2|\ln\varepsilon|).
$$
Straightforward calculation shows that
$$
\mathcal{T}_{H_0,\lambda}^{pp,\varepsilon}[e^{\ii n\theta}] =\frac{\ii}{2}\operatorname{sgn}(n)e^{\ii n\theta}+\mathcal{O}(\varepsilon^2|\ln\varepsilon|),\quad  n\ne0,
$$
and $\mathcal{T}_{H_0,\lambda}^{pp,\varepsilon}[1]=\mathcal{O}(\varepsilon^2|\ln\varepsilon|)$.
\end{proof}

\begin{lemma}
The tangential operator induced by $G_{\Lambda_0^*}^{\lambda}$ satisfies
\begin{equation}\label{eq:self-T-res-general}
\mathcal{T}_{\Lambda_0^*,\lambda}^{pp,\varepsilon}[\varphi](\theta) =\frac{2Q^2\varepsilon^2}{3|Y|\Delta_\lambda} \int_0^{2\pi}\sin(\theta-t)\varphi(t) dt +\mathcal{O} \left(\frac{\varepsilon^3}{|\Delta_\lambda|}\right).
\end{equation}
Moreover,
\begin{equation}\label{eq:self-T-res-fourier}
\mathcal{T}_{\Lambda_0^*,\lambda}^{pp,\varepsilon}[e^{\ii n\theta}] =
\begin{cases}
-\dfrac{2\pi \ii Q^2\varepsilon^2}{3|Y|\Delta_\lambda}e^{\ii\theta} +\mathcal{O}(\varepsilon^3/|\Delta_\lambda|),&n=1,\\
\dfrac{2\pi \ii Q^2\varepsilon^2}{3|Y|\Delta_\lambda}e^{-\ii\theta} +\mathcal{O}(\varepsilon^3/|\Delta_\lambda|),&n=-1,\\
\mathcal{O}(\varepsilon^3/|\Delta_\lambda|),&n=0,\\
\mathcal{O}(\varepsilon^3/|\Delta_\lambda|),&|n|>1.
\end{cases}
\end{equation}
\end{lemma}
\begin{proof}
Set $\bm d(\theta,t)=\bm r(\theta)-\bm r(t)$. The resonant-kernel gradient is
$$
\nabla_{\bm x}G_{\Lambda_0^*}^{\lambda}(\varepsilon\bm d) =-\frac{2Q^2\varepsilon}{3|Y|\Delta_\lambda}\bm d +\mathcal{O} \left(\frac{\varepsilon^2}{|\Delta_\lambda|}\right).
$$
Therefore
\begin{align*}
\mathcal{T}_{\Lambda_0^*,\lambda}^{pp,\varepsilon}[\varphi](\theta) &=\varepsilon\int_0^{2\pi} \bm\tau(\theta)\cdot \nabla_{\bm x}G_{\Lambda_0^*}^{\lambda}(\varepsilon\bm d(\theta,t)) \varphi(t) dt\\
&=\frac{2Q^2\varepsilon^2}{3|Y|\Delta_\lambda} \int_0^{2\pi}\sin(\theta-t)\varphi(t) dt +\mathcal{O} \left(\frac{\varepsilon^3}{|\Delta_\lambda|}\right).
\end{align*}
Also,
$$
\int_0^{2\pi}\sin(\theta-t)e^{\ii nt} dt =
\begin{cases}
-\ii\pi e^{\ii\theta},&n=1,\\
\ii\pi e^{-\ii\theta},&n=-1,\\
0,&n=0\text{ or }|n|>1.
\end{cases}
$$
Substituting these Fourier integrals into~\eqref{eq:self-T-res-general} gives~\eqref{eq:self-T-res-fourier}.
\end{proof}

\begin{lemma}
For $\bm x\in\partial D_p^\varepsilon$, define the tangential operator associated with $\widetilde G^\lambda$ by
\begin{equation}\label{eq:T-smooth-definition}
\widetilde{\mathcal{T}}_\lambda^{pp,\varepsilon}[\varphi](\bm x) :=\operatorname{p.v.}\int_{\partial D_p^\varepsilon} \bm\tau(\bm x)\cdot\nabla_{\bm x}\widetilde G^\lambda(\bm x-\bm y) \varphi(\bm y) d\sigma(\bm y).
\end{equation}
In circular coordinates,
\begin{equation}\label{eq:self-T-smooth}
\widetilde{\mathcal{T}}_\lambda^{pp,\varepsilon}[\varphi](\theta) =\varepsilon\bm\tau(\theta)\cdot\nabla_{\bm x}\widetilde G^\lambda(\bm0) \int_0^{2\pi}\varphi(t) dt+\mathcal{O}(\varepsilon^2).
\end{equation}
For $n\ne0$, $\widetilde{\mathcal{T}}_\lambda^{pp,\varepsilon}[e^{\ii n\theta}]=\mathcal{O}(\varepsilon^2)$.
\end{lemma}

\begin{proof}
The gradient in~\eqref{eq:T-smooth-definition} acts on the first spatial argument of the Green function. Let $\widetilde{\bm g}_\lambda=\nabla_{\bm x}\widetilde G^\lambda(\bm0)$. Then
$$
\nabla_{\bm x}\widetilde G^\lambda (\varepsilon(\bm r(\theta)-\bm r(t))) =\widetilde{\bm g}_\lambda+\mathcal{O}(\varepsilon).
$$
After parametrization,
\begin{align*}
\widetilde{\mathcal{T}}_\lambda^{pp,\varepsilon}[\varphi](\theta) &=\varepsilon\int_0^{2\pi} \bm\tau(\theta)\cdot \nabla_{\bm x}\widetilde G^\lambda (\varepsilon(\bm r(\theta)-\bm r(t)))\varphi(t) dt\\
&=\varepsilon\bm\tau(\theta)\cdot\widetilde{\bm g}_\lambda \int_0^{2\pi}\varphi(t) dt+\mathcal{O}(\varepsilon^2).
\end{align*}
For $n\ne0$, $\int_0^{2\pi}e^{\ii nt} dt=0$, so $\widetilde{\mathcal{T}}_\lambda^{pp,\varepsilon}[e^{\ii n\theta}]=\mathcal{O}(\varepsilon^2)$.
\end{proof}

\subsection{Expansions between distinct boundary components}

Take $p,q\in\{A,B\}$ with $p\ne q$. Throughout this subsection, we use $\mathcal{S}_{\lambda}^{pq,\varepsilon}[\varphi](\bm x), \mathcal{K}_{\lambda}^{pq,\varepsilon,*}[\varphi](\bm x), \mathcal{T}_{\lambda}^{pq,\varepsilon}[\varphi](\bm x)$ to denote the corresponding operators from component $q$ to component $p$, evaluated at $\bm x\in\partial D_p^\varepsilon$, with boundary density $\varphi$ supported on component $q$. Assume $\bm x\in\partial D_p^\varepsilon$ and $\bm y\in\partial D_q^\varepsilon$ are on two different boundaries. Set $\bm x=\bm z_p+\varepsilon\bm r(\theta)$, $\bm y=\bm z_q+\varepsilon\bm r(t)$, $\bm d_{pq}:=\bm z_p-\bm z_q$. Since $\bm d_{pq}\ne\bm0$, there is $c_0>0$, independent of $\varepsilon$, such that, for sufficiently small $\varepsilon$,
\begin{equation}\label{eq:cross-positive-distance}
|\bm x-\bm y|\ge c_0,  \qquad \bm x\in\partial D_p^\varepsilon, \quad  \bm y\in\partial D_q^\varepsilon.
\end{equation}

\begin{lemma}
For $p\ne q$, the normal traces on component $p$ of the layer generated on component $q$ satisfy
\begin{equation}\label{eq:cross-normal-trace}
\left.\partial_{\bm\nu}\mathcal{S}_{\lambda}^{pq,\varepsilon}[\varphi]\right|_+ =\left.\partial_{\bm\nu}\mathcal{S}_{\lambda}^{pq,\varepsilon}[\varphi]\right|_- =\mathcal{K}_{\lambda}^{pq,\varepsilon,*}[\varphi].
\end{equation}
\end{lemma}

\begin{proof}
Equation~\eqref{eq:cross-positive-distance} gives $|\bm x-\bm y|\ge c_0$. Hence $\nabla_{\bm x}G_\lambda^{\bm\alpha^*}(\bm x-\bm y)$ is continuous on $\partial D_p^\varepsilon\times\partial D_q^\varepsilon$. For any $\bm x_0\in\partial D_p^\varepsilon$, approaching from either side gives
$$
\partial_{\bm\nu}\mathcal{S}_{\lambda}^{pq,\varepsilon}[\varphi](\bm x) \longrightarrow \int_{\partial D_q^\varepsilon} \bm\nu(\bm x_0)\cdot\nabla_{\bm x} G_\lambda^{\bm\alpha^*}(\bm x_0-\bm y) \varphi(\bm y) d\sigma(\bm y).
$$
The limits coincide, and the right-hand side is $\mathcal{K}_{\lambda}^{pq,\varepsilon,*}[\varphi](\bm x_0)$, proving~\eqref{eq:cross-normal-trace}.
\end{proof}

\begin{lemma}
For $p\ne q$, the three cross-interaction operators satisfy
\begin{align*}
\mathcal{S}_{\lambda}^{pq,\varepsilon}[\varphi](\theta) &=\varepsilon G_\lambda^{\bm\alpha^*}(\bm d_{pq}) \int_0^{2\pi}\varphi(t) dt+\varepsilon^2\int_0^{2\pi} \nabla_{\bm x}G_\lambda^{\bm\alpha^*}(\bm d_{pq})\cdot (\bm r(\theta)-\bm r(t))\varphi(t) dt +\mathcal{O}\left(\varepsilon^3(1+|\Delta_\lambda|^{-1})\right),\\
\mathcal{K}_{\lambda}^{pq,\varepsilon,*}[\varphi](\theta) &=\varepsilon\bm\nu(\theta)\cdot \nabla_{\bm x}G_\lambda^{\bm\alpha^*}(\bm d_{pq})\int_0^{2\pi}\varphi(t) dt+\mathcal{O}\left(\varepsilon^2(1+|\Delta_\lambda|^{-1})\right),\\
\mathcal{T}_{\lambda}^{pq,\varepsilon}[\varphi](\theta) &=\varepsilon\bm\tau(\theta)\cdot \nabla_{\bm x}G_\lambda^{\bm\alpha^*}(\bm d_{pq})\int_0^{2\pi}\varphi(t) dt+\mathcal{O}\left(\varepsilon^2(1+|\Delta_\lambda|^{-1})\right).
\end{align*}
The kernel remainders and the derivatives needed for the stated operator norms are uniform in $(\theta,t)\in[0,2\pi]^2$.
\end{lemma}
\begin{proof}
Set $\bm d(\theta,t)=\bm r(\theta)-\bm r(t)$. Since $|\bm d_{pq}|>0$, equation~\eqref{eq:cross-positive-distance} permits the uniform local expansion
$$
G_\lambda^{\bm\alpha^*}(\bm d_{pq}+\varepsilon\bm d) =G_\lambda^{\bm\alpha^*}(\bm d_{pq}) +\varepsilon\nabla_{\bm x}G_\lambda^{\bm\alpha^*}(\bm d_{pq})\cdot\bm d +\mathcal{O} \left(\varepsilon^2(1+|\Delta_\lambda|^{-1})\right),
$$
and
$$
\nabla_{\bm x}G_\lambda^{\bm\alpha^*}(\bm d_{pq}+\varepsilon\bm d) =\nabla_{\bm x}G_\lambda^{\bm\alpha^*}(\bm d_{pq})+\mathcal{O} \left(\varepsilon(1+|\Delta_\lambda|^{-1})\right).
$$
The single-layer operator is
\begin{align*}
\mathcal{S}_{\lambda}^{pq,\varepsilon}[\varphi](\theta) &=\varepsilon\int_0^{2\pi} G_\lambda^{\bm\alpha^*}(\bm d_{pq}+\varepsilon\bm d(\theta,t))\varphi(t) dt\\
&=\varepsilon G_\lambda^{\bm\alpha^*}(\bm d_{pq})\int_0^{2\pi}\varphi(t) dt +\varepsilon^2\int_0^{2\pi} \nabla_{\bm x}G_\lambda^{\bm\alpha^*}(\bm d_{pq})\cdot\bm d(\theta,t)\varphi(t) dt +\mathcal{O} \left(\varepsilon^3(1+|\Delta_\lambda|^{-1})\right).
\end{align*}
The normal and tangential operators are, respectively,
\begin{align*}
\mathcal{K}_{\lambda}^{pq,\varepsilon,*}[\varphi](\theta) &=\varepsilon\int_0^{2\pi} \bm\nu(\theta)\cdot\nabla_{\bm x}G_\lambda^{\bm\alpha^*} (\bm d_{pq}+\varepsilon\bm d(\theta,t))\varphi(t) dt\\
&=\varepsilon\bm\nu(\theta)\cdot\nabla_{\bm x}G_\lambda^{\bm\alpha^*}(\bm d_{pq}) \int_0^{2\pi}\varphi(t) dt+\mathcal{O} \left(\varepsilon^2(1+|\Delta_\lambda|^{-1})\right),\\
\mathcal{T}_{\lambda}^{pq,\varepsilon}[\varphi](\theta) &=\varepsilon\int_0^{2\pi} \bm\tau(\theta)\cdot\nabla_{\bm x}G_\lambda^{\bm\alpha^*} (\bm d_{pq}+\varepsilon\bm d(\theta,t))\varphi(t) dt\\
&=\varepsilon\bm\tau(\theta)\cdot\nabla_{\bm x}G_\lambda^{\bm\alpha^*}(\bm d_{pq}) \int_0^{2\pi}\varphi(t) dt+\mathcal{O} \left(\varepsilon^2(1+|\Delta_\lambda|^{-1})\right).
\end{align*}
\end{proof}

\begin{proposition}\label{prop:resonant-cross-leading}
For $p\ne q$, the cross-component normal and tangential operators induced by $G_{\Lambda_0^*}^{\lambda}$ satisfy
\begin{align}
\mathcal{K}_{\Lambda_0^*,\lambda}^{pq,\varepsilon,*}[\varphi](\theta) &=\frac{\ii\varepsilon}{|Y|\Delta_\lambda} \bm\nu(\theta)\cdot\bm P(\bm d_{pq}) \int_0^{2\pi}\varphi(t) dt +\mathcal{O} \left(\frac{\varepsilon^2}{|\Delta_\lambda|}\right),\label{eq:cross-K-res-leading}\\
\mathcal{T}_{\Lambda_0^*,\lambda}^{pq,\varepsilon}[\varphi](\theta) &=\frac{\ii\varepsilon}{|Y|\Delta_\lambda} \bm\tau(\theta)\cdot\bm P(\bm d_{pq}) \int_0^{2\pi}\varphi(t) dt +\mathcal{O} \left(\frac{\varepsilon^2}{|\Delta_\lambda|}\right),\label{eq:cross-T-res-leading}
\end{align}
where
$$
\bm P(\bm d_{pq}) :=\sum_{\bm q_j\in\Lambda_0^*} (\bm\alpha^*+\bm q_j)e^{\ii(\bm\alpha^*+\bm q_j)\cdot\bm d_{pq}}.
$$
\end{proposition}

\begin{proof}
Set $\bm d(\theta,t)=\bm r(\theta)-\bm r(t)$. The gradient of $G_{\Lambda_0^*}^{\lambda}$ is
\begin{align*}
\nabla_{\bm x}G_{\Lambda_0^*}^{\lambda} (\bm d_{pq}+\varepsilon\bm d) &=\frac{\ii}{|Y|\Delta_\lambda} \sum_{\bm q_j\in\Lambda_0^*} (\bm\alpha^*+\bm q_j)e^{\ii(\bm\alpha^*+\bm q_j)\cdot\bm d_{pq}} e^{\ii\varepsilon(\bm\alpha^*+\bm q_j)\cdot\bm d}\\
&=\frac{\ii}{|Y|\Delta_\lambda}\bm P(\bm d_{pq}) +\mathcal{O} \left(\frac{\varepsilon}{|\Delta_\lambda|}\right).
\end{align*}
Hence
\begin{align*}
\mathcal{K}_{\Lambda_0^*,\lambda}^{pq,\varepsilon,*}[\varphi](\theta) &=\varepsilon\int_0^{2\pi} \bm\nu(\theta)\cdot\nabla_{\bm x}G_{\Lambda_0^*}^{\lambda} (\bm d_{pq}+\varepsilon\bm d(\theta,t))\varphi(t) dt\\
&=\frac{\ii\varepsilon}{|Y|\Delta_\lambda} \bm\nu(\theta)\cdot\bm P(\bm d_{pq}) \int_0^{2\pi}\varphi(t) dt +\mathcal{O} \left(\frac{\varepsilon^2}{|\Delta_\lambda|}\right),
\end{align*}
and
\begin{align*}
\mathcal{T}_{\Lambda_0^*,\lambda}^{pq,\varepsilon}[\varphi](\theta) &=\varepsilon\int_0^{2\pi} \bm\tau(\theta)\cdot\nabla_{\bm x}G_{\Lambda_0^*}^{\lambda} (\bm d_{pq}+\varepsilon\bm d(\theta,t))\varphi(t) dt\\
&=\frac{\ii\varepsilon}{|Y|\Delta_\lambda} \bm\tau(\theta)\cdot\bm P(\bm d_{pq}) \int_0^{2\pi}\varphi(t) dt +\mathcal{O} \left(\frac{\varepsilon^2}{|\Delta_\lambda|}\right).
\end{align*}
\end{proof}

\begin{proposition}
When $\Delta_\lambda=\varepsilon^2z$ with $z\neq 0$, for every $\varphi\in\operatorname{span}\{e^{\ii t},e^{-\ii t}\}$, the $\mathcal{O}(\varepsilon/\Delta_\lambda)$ leading terms in Proposition~\ref{prop:resonant-cross-leading} vanish exactly, and
\begin{equation}\label{eq:cross-V1-bounds}
\mathcal{K}_{\Lambda_0^*,\lambda}^{pq,\varepsilon,*}[\varphi] =\mathcal{O}\left(\frac{\varepsilon^2}{|\Delta_\lambda|}\right), \qquad  \mathcal{T}_{\Lambda_0^*,\lambda}^{pq,\varepsilon}[\varphi] =\mathcal{O}\left(\frac{\varepsilon^2}{|\Delta_\lambda|}\right).
\end{equation}
In particular, these bounds are $\mathcal{O}(1)$ in the stated scaling.
\end{proposition}

\begin{proof}
Every $\varphi\in\operatorname{span}\{e^{\ii t},e^{-\ii t}\}$ has the form $\varphi(t)=c_1e^{\ii t}+c_{-1}e^{-\ii t}$, and hence $\int_0^{2\pi}\varphi(t) dt=0$. Substitution into~\eqref{eq:cross-K-res-leading} and~\eqref{eq:cross-T-res-leading} cancels the explicit $\varepsilon/\Delta_\lambda$ terms. The remaining terms satisfy~\eqref{eq:cross-V1-bounds}.
\end{proof}

\subsection{Representations of the nine basic Fourier matrix elements}\label{sec:nine-fourier-elements}

Equation~\eqref{eq:operator-three-part} decomposes each of $\mathcal{S}_\lambda^\varepsilon$, $\mathcal{K}_\lambda^{\varepsilon,*}$, and $\mathcal{T}_\lambda^\varepsilon$ into contributions from $H_0$, $G_{\Lambda_0^*}^\lambda$, and $\widetilde G^\lambda$. For $F\in\{\mathcal{S}_{H_0,\lambda},\mathcal{S}_{\Lambda_0^*,\lambda},\widetilde{\mathcal{S}}_{\lambda}$, $\mathcal{K}_{H_0,\lambda}^*,\mathcal{K}_{\Lambda_0^*,\lambda}^*,\widetilde{\mathcal{K}}_{\lambda}^*$, $\mathcal{T}_{H_0,\lambda},\mathcal{T}_{\Lambda_0^*,\lambda},\widetilde{\mathcal{T}}_{\lambda}\}$ and $p,q\in\{A,B\}$, define
\begin{equation}\label{eq:matrix-element-S}
(F^{pq,\varepsilon})_{m,n} :=\left\langle\varphi_{p,m},F^{pq,\varepsilon}\varphi_{q,n}\right\rangle_{p,\varepsilon}.
\end{equation}

Fix $p=q$, so $\bm x-\bm y=\varepsilon(\bm r(\theta)-\bm r(t))$.
\begin{lemma}
Let $p=q$. As $\varepsilon\to0$, the single-layer entries satisfy
$$
(\mathcal{S}_{H_0,\lambda}^{pp,\varepsilon})_{m,n} =
\begin{cases}
\varepsilon(\ln\varepsilon+\ln\lambda+\gamma_0)+\mathcal{O}(\varepsilon^3|\ln\varepsilon|), & m=n=0,\\
-\dfrac{\varepsilon}{2|n|}+\mathcal{O}(\varepsilon^3|\ln\varepsilon|), & m=n\neq0,\\
\mathcal{O}(\varepsilon^3|\ln\varepsilon|), & m\neq n.
\end{cases}
$$
The free-space Neumann--Poincar\'e entries satisfy
$$
(\mathcal{K}_{H_0,\lambda}^{pp,\varepsilon,*})_{m,n} =
\begin{cases}
\dfrac12+\mathcal{O}(\varepsilon^2|\ln\varepsilon|), & m=n=0,\\
\mathcal{O}(\varepsilon^2|\ln\varepsilon|), & \text{otherwise}.
\end{cases}
$$
With $\bm\tau(\theta)=(\sin\theta,-\cos\theta)$, the free-space tangential entries satisfy
$$
(\mathcal{T}_{H_0,\lambda}^{pp,\varepsilon})_{m,n} =
\begin{cases}
\dfrac{\ii}{2}\operatorname{sgn}(n)+\mathcal{O}(\varepsilon^2|\ln\varepsilon|), & m=n\neq0,\\
\mathcal{O}(\varepsilon^2|\ln\varepsilon|), & \text{otherwise}.
\end{cases}
$$
\end{lemma}

\begin{proof}
The normalized basis satisfies
$$
\langle\varphi_m,\varphi_n\rangle_{L^2(0,2\pi)}=\delta_{mn}, \qquad  \int_0^{2\pi}\varphi_n(t) dt=\sqrt{2\pi} \delta_{n0}.
$$
Equation~\eqref{eq:circle-log-fourier-integral} becomes
$$
\frac1{2\pi}\int_0^{2\pi} \ln|\bm r(\theta)-\bm r(t)|\varphi_n(t) dt =
\begin{cases}
0,&n=0,\\
-\dfrac1{2|n|}\varphi_n(\theta),&n\ne0.
\end{cases}
$$
Substitution into~\eqref{eq:self-S-H0-general} gives
$$
\mathcal{S}_{H_0,\lambda}^{pp,\varepsilon}[\varphi_n] =
\begin{cases}
\varepsilon(\ln\varepsilon+\ln\lambda+\gamma_0)\varphi_0 +\mathcal{O}(\varepsilon^3|\ln\varepsilon|),&n=0,\\
-\dfrac{\varepsilon}{2|n|}\varphi_n +\mathcal{O}(\varepsilon^3|\ln\varepsilon|),&n\ne0.
\end{cases}
$$
The corresponding entries are
$$
(\mathcal{S}_{H_0,\lambda}^{pp,\varepsilon})_{m,n} =
\begin{cases}
\varepsilon(\ln\varepsilon+\ln\lambda+\gamma_0)\delta_{m0} +\mathcal{O}(\varepsilon^3|\ln\varepsilon|),&n=0,\\
-\dfrac{\varepsilon}{2|n|}\delta_{mn} +\mathcal{O}(\varepsilon^3|\ln\varepsilon|),&n\ne0.
\end{cases}
$$

Equation~\eqref{eq:self-K-H0-general} gives
$$
\mathcal{K}_{H_0,\lambda}^{pp,\varepsilon,*}[\varphi_n] =\frac1{4\pi}\sqrt{2\pi} \delta_{n0} +\mathcal{O}(\varepsilon^2|\ln\varepsilon|) =\frac12\delta_{n0}\varphi_0 +\mathcal{O}(\varepsilon^2|\ln\varepsilon|).
$$
Hence,
$$
(\mathcal{K}_{H_0,\lambda}^{pp,\varepsilon,*})_{m,n} =\frac12\delta_{m0}\delta_{n0}+\mathcal{O}(\varepsilon^2|\ln\varepsilon|).
$$
From equation~\eqref{eq:self-T-H0-general}, straightforward calculation gives
$$
\mathcal{T}_{H_0,\lambda}^{pp,\varepsilon}[\varphi_n] =
\begin{cases}
\dfrac{\ii}{2}\operatorname{sgn}(n)\varphi_n +\mathcal{O}(\varepsilon^2|\ln\varepsilon|),&n\ne0,\\
\mathcal{O}(\varepsilon^2|\ln\varepsilon|),&n=0.
\end{cases}
$$
The corresponding entries are
$$
(\mathcal{T}_{H_0,\lambda}^{pp,\varepsilon})_{m,n} =
\begin{cases}
\dfrac{\ii}{2}\operatorname{sgn}(n)\delta_{mn} +\mathcal{O}(\varepsilon^2|\ln\varepsilon|),&n\ne0,\\
\mathcal{O}(\varepsilon^2|\ln\varepsilon|),&n=0.
\end{cases}
$$
\end{proof}

\begin{lemma}
Let $p=q$ and $\Delta_\lambda=\lambda^2-|\bm\alpha^*|^2$. The corresponding single-layer entries satisfy
$$
(\mathcal{S}_{\Lambda_0^*,\lambda}^{pp,\varepsilon})_{m,n} =
\begin{cases}
\dfrac{6\pi\varepsilon}{|Y|\Delta_\lambda} -\dfrac{4\pi Q^2\varepsilon^3}{3|Y|\Delta_\lambda} +\mathcal{O}\left(\dfrac{\varepsilon^4}{|\Delta_\lambda|}\right), & m=n=0,\\
\dfrac{2\pi Q^2\varepsilon^3}{3|Y|\Delta_\lambda} +\mathcal{O}\left(\dfrac{\varepsilon^4}{|\Delta_\lambda|}\right), & m=n=1\text{ or }m=n=-1,\\
\mathcal{O}\left(\dfrac{\varepsilon^4}{|\Delta_\lambda|}\right), & \text{otherwise}.
\end{cases}
$$
The corresponding Neumann--Poincar\'e entries satisfy
$$
(\mathcal{K}_{\Lambda_0^*,\lambda}^{pp,\varepsilon,*})_{m,n} =
\begin{cases}
-\dfrac{4\pi Q^2\varepsilon^2}{3|Y|\Delta_\lambda} +\mathcal{O}\left(\dfrac{\varepsilon^3}{|\Delta_\lambda|}\right), & m=n=0,\\
\dfrac{2\pi Q^2\varepsilon^2}{3|Y|\Delta_\lambda} +\mathcal{O}\left(\dfrac{\varepsilon^3}{|\Delta_\lambda|}\right), & m=n=1\text{ or }m=n=-1,\\
\mathcal{O}\left(\dfrac{\varepsilon^3}{|\Delta_\lambda|}\right), & \text{otherwise}.
\end{cases}
$$
The corresponding tangential entries satisfy
$$
(\mathcal{T}_{\Lambda_0^*,\lambda}^{pp,\varepsilon})_{m,n} =
\begin{cases}
-\dfrac{2\pi \ii Q^2\varepsilon^2}{3|Y|\Delta_\lambda} +\mathcal{O}\left(\dfrac{\varepsilon^3}{|\Delta_\lambda|}\right), & m=n=1,\\
\dfrac{2\pi \ii Q^2\varepsilon^2}{3|Y|\Delta_\lambda} +\mathcal{O}\left(\dfrac{\varepsilon^3}{|\Delta_\lambda|}\right), & m=n=-1,\\
\mathcal{O}\left(\dfrac{\varepsilon^3}{|\Delta_\lambda|}\right), & \text{otherwise}.
\end{cases}
$$
\end{lemma}

\begin{proof}
Set $\bm d(\theta,t)=\bm r(\theta)-\bm r(t)$. Fourier orthogonality gives
$$
\int_0^{2\pi}|\bm d(\theta,t)|^2\varphi_n(t) dt =
\begin{cases}
4\pi\varphi_0(\theta),&n=0,\\
-2\pi\varphi_n(\theta),&n=\pm1,\\
0,&|n|>1,
\end{cases}
$$
because $|\bm d(\theta,t)|^2=2-2\cos(\theta-t)$. Also, $\int_0^{2\pi}\varphi_n(t) dt=\sqrt{2\pi}\delta_{n0}$. Thus~\eqref{eq:self-S-res} becomes
$$
\mathcal{S}_{\Lambda_0^*,\lambda}^{pp,\varepsilon}[\varphi_n] =
\begin{cases}
\left(\dfrac{6\pi\varepsilon}{|Y|\Delta_\lambda} -\dfrac{4\pi Q^2\varepsilon^3}{3|Y|\Delta_\lambda}\right)\varphi_0 +\mathcal{O} \left(\dfrac{\varepsilon^4}{|\Delta_\lambda|}\right),&n=0,\\
\dfrac{2\pi Q^2\varepsilon^3}{3|Y|\Delta_\lambda}\varphi_n +\mathcal{O} \left(\dfrac{\varepsilon^4}{|\Delta_\lambda|}\right),&n=\pm1,\\
\mathcal{O} \left(\dfrac{\varepsilon^4}{|\Delta_\lambda|}\right),&|n|>1.
\end{cases}
$$
Orthogonality yields
$$
(\mathcal{S}_{\Lambda_0^*,\lambda}^{pp,\varepsilon})_{m,n} =
\begin{cases}
\left(\dfrac{6\pi\varepsilon}{|Y|\Delta_\lambda} -\dfrac{4\pi Q^2\varepsilon^3}{3|Y|\Delta_\lambda}\right)\delta_{m0} +\mathcal{O} \left(\dfrac{\varepsilon^4}{|\Delta_\lambda|}\right),&n=0,\\
\dfrac{2\pi Q^2\varepsilon^3}{3|Y|\Delta_\lambda}\delta_{mn} +\mathcal{O} \left(\dfrac{\varepsilon^4}{|\Delta_\lambda|}\right),&n=\pm1,\\
\mathcal{O} \left(\dfrac{\varepsilon^4}{|\Delta_\lambda|}\right),&|n|>1.
\end{cases}
$$

Similarly,
$$
\int_0^{2\pi}(1-\cos(\theta-t))\varphi_n(t) dt =
\begin{cases}
2\pi\varphi_0(\theta),&n=0,\\
-\pi\varphi_n(\theta),&n=\pm1,\\
0,&|n|>1.
\end{cases}
$$
Substituting into~\eqref{eq:self-K-res-general} gives
$$
\mathcal{K}_{\Lambda_0^*,\lambda}^{pp,\varepsilon,*}[\varphi_n] =
\begin{cases}
-\dfrac{4\pi Q^2\varepsilon^2}{3|Y|\Delta_\lambda}\varphi_0 +\mathcal{O} \left(\dfrac{\varepsilon^3}{|\Delta_\lambda|}\right),&n=0,\\
\dfrac{2\pi Q^2\varepsilon^2}{3|Y|\Delta_\lambda}\varphi_n +\mathcal{O} \left(\dfrac{\varepsilon^3}{|\Delta_\lambda|}\right),&n=\pm1,\\
\mathcal{O} \left(\dfrac{\varepsilon^3}{|\Delta_\lambda|}\right),&|n|>1.
\end{cases}
$$

The tangential integrals satisfy
$$
\int_0^{2\pi}\sin(\theta-t)\varphi_n(t) dt =
\begin{cases}
-\ii\pi\varphi_1(\theta),&n=1,\\
\ii\pi\varphi_{-1}(\theta),&n=-1,\\
0,&n=0\text{ or }|n|>1.
\end{cases}
$$
Equation~\eqref{eq:self-T-res-general} therefore gives
$$
\mathcal{T}_{\Lambda_0^*,\lambda}^{pp,\varepsilon}[\varphi_n] =
\begin{cases}
-\dfrac{2\pi \ii Q^2\varepsilon^2}{3|Y|\Delta_\lambda}\varphi_1 +\mathcal{O} \left(\dfrac{\varepsilon^3}{|\Delta_\lambda|}\right),&n=1,\\ 
\dfrac{2\pi \ii Q^2\varepsilon^2}{3|Y|\Delta_\lambda}\varphi_{-1} +\mathcal{O} \left(\dfrac{\varepsilon^3}{|\Delta_\lambda|}\right),&n=-1,\\ 
\mathcal{O} \left(\dfrac{\varepsilon^3}{|\Delta_\lambda|}\right),&n=0\text{ or }|n|>1.
\end{cases}
$$
Orthogonality gives
$$
(\mathcal{K}_{\Lambda_0^*,\lambda}^{pp,\varepsilon,*})_{m,n} =
\begin{cases}
-\dfrac{4\pi Q^2\varepsilon^2}{3|Y|\Delta_\lambda}\delta_{m0} +\mathcal{O} \left(\dfrac{\varepsilon^3}{|\Delta_\lambda|}\right),&n=0,\\
\dfrac{2\pi Q^2\varepsilon^2}{3|Y|\Delta_\lambda}\delta_{mn} +\mathcal{O} \left(\dfrac{\varepsilon^3}{|\Delta_\lambda|}\right),&n=\pm1,\\
\mathcal{O} \left(\dfrac{\varepsilon^3}{|\Delta_\lambda|}\right),&|n|>1,
\end{cases}
$$
and
$$
(\mathcal{T}_{\Lambda_0^*,\lambda}^{pp,\varepsilon})_{m,n} =
\begin{cases}
-\dfrac{2\pi \ii Q^2\varepsilon^2}{3|Y|\Delta_\lambda}\delta_{m1} +\mathcal{O} \left(\dfrac{\varepsilon^3}{|\Delta_\lambda|}\right),&n=1,\\
\dfrac{2\pi \ii Q^2\varepsilon^2}{3|Y|\Delta_\lambda}\delta_{m,-1} +\mathcal{O} \left(\dfrac{\varepsilon^3}{|\Delta_\lambda|}\right),&n=-1,\\
\mathcal{O} \left(\dfrac{\varepsilon^3}{|\Delta_\lambda|}\right),&n=0\text{ or }|n|>1.
\end{cases}
$$
\end{proof}

\begin{lemma}\label{lem:self-interaction-smooth}
Let $p=q$ and write
$$
\widetilde{\bm g}_\lambda=\nabla_{\bm x}\widetilde G^\lambda(\bm0) =(\widetilde g_{\lambda,1},\widetilde g_{\lambda,2}), \qquad  \widetilde g_{\lambda,+}:=\widetilde g_{\lambda,1}+\ii\widetilde g_{\lambda,2},  \qquad \widetilde g_{\lambda,-}:=\widetilde g_{\lambda,1}-\ii\widetilde g_{\lambda,2}.
$$
The smooth single-layer entries satisfy
$$
(\widetilde{\mathcal{S}}_{\lambda}^{pp,\varepsilon})_{m,n} =
\begin{cases}
2\pi\varepsilon \widetilde G^\lambda(0)+\mathcal{O}(\varepsilon^3), & m=n=0,\\
\pi\varepsilon^2\widetilde g_{\lambda,-}+\mathcal{O}(\varepsilon^3), & m=1, n=0,\\
\pi\varepsilon^2\widetilde g_{\lambda,+}+\mathcal{O}(\varepsilon^3), & m=-1, n=0,\\
-\pi\varepsilon^2\widetilde g_{\lambda,+}+\mathcal{O}(\varepsilon^3), & m=0, n=1,\\
-\pi\varepsilon^2\widetilde g_{\lambda,-}+\mathcal{O}(\varepsilon^3), & m=0, n=-1,\\
\mathcal{O}(\varepsilon^3), & \text{otherwise}.
\end{cases}
$$
The smooth Neumann--Poincar\'e entries satisfy
$$
(\widetilde{\mathcal{K}}_{\lambda}^{pp,\varepsilon,*})_{m,n} =
\begin{cases}
\pi\varepsilon \widetilde g_{\lambda,-}+\mathcal{O}(\varepsilon^2), & m=1, n=0,\\
\pi\varepsilon \widetilde g_{\lambda,+}+\mathcal{O}(\varepsilon^2), & m=-1, n=0,\\
\mathcal{O}(\varepsilon^2), & \text{otherwise}.
\end{cases}
$$
The smooth tangential entries satisfy
$$
(\widetilde{\mathcal{T}}_{\lambda}^{pp,\varepsilon})_{m,n} =
\begin{cases}
-\ii\pi\varepsilon \widetilde g_{\lambda,-}+\mathcal{O}(\varepsilon^2), & m=1, n=0,\\
\ii\pi\varepsilon \widetilde g_{\lambda,+}+\mathcal{O}(\varepsilon^2), & m=-1, n=0,\\
\mathcal{O}(\varepsilon^2), & \text{otherwise}.
\end{cases}
$$
\end{lemma}

\begin{proof}
Write
$$
\widetilde{\bm g}_\lambda\cdot\bm r(\theta) =\frac12\widetilde g_{\lambda,-}e^{\ii\theta} +\frac12\widetilde g_{\lambda,+}e^{-\ii\theta},
$$
and
$$
\bm\tau(\theta)\cdot\widetilde{\bm g}_\lambda =-\frac{\ii}{2}\widetilde g_{\lambda,-}e^{\ii\theta} +\frac{\ii}{2}\widetilde g_{\lambda,+}e^{-\ii\theta}.
$$
Using
$$
\int_0^{2\pi}\varphi_n(t) dt=\sqrt{2\pi}\delta_{n0}, \qquad  e^{\pm \ii\theta}=\sqrt{2\pi} \varphi_{\pm1}(\theta).
$$
Equation~\eqref{eq:self-S-smooth} gives, for $n=0,1,-1$,
\begin{align*}
\widetilde{\mathcal{S}}_\lambda^{pp,\varepsilon}[\varphi_0] &=2\pi\varepsilon\widetilde G^\lambda(\bm0)\varphi_0 +\pi\varepsilon^2 (\widetilde g_{\lambda,-}\varphi_1 +\widetilde g_{\lambda,+}\varphi_{-1})+\mathcal{O}(\varepsilon^3),\\
\widetilde{\mathcal{S}}_\lambda^{pp,\varepsilon}[\varphi_1] &=-\pi\varepsilon^2\widetilde g_{\lambda,+}\varphi_0+\mathcal{O}(\varepsilon^3),\\
\widetilde{\mathcal{S}}_\lambda^{pp,\varepsilon}[\varphi_{-1}] &=-\pi\varepsilon^2\widetilde g_{\lambda,-}\varphi_0+\mathcal{O}(\varepsilon^3).
\end{align*}
For $|n|>1$, the constant and linear Fourier integrals vanish, so $\widetilde{\mathcal{S}}_\lambda^{pp,\varepsilon}[\varphi_n]=\mathcal{O}(\varepsilon^3)$. Thus, by $\langle\varphi_m,\varphi_k\rangle=\delta_{mk}$,
\begin{align*}
(\widetilde{\mathcal{S}}_\lambda^{pp,\varepsilon})_{m,0} &=2\pi\varepsilon\widetilde G^\lambda(\bm0)\delta_{m0} +\pi\varepsilon^2\bigl( \widetilde g_{\lambda,-}\delta_{m1} +\widetilde g_{\lambda,+}\delta_{m,-1}\bigr)+\mathcal{O}(\varepsilon^3),\\
(\widetilde{\mathcal{S}}_\lambda^{pp,\varepsilon})_{m,1} &=-\pi\varepsilon^2\widetilde g_{\lambda,+}\delta_{m0}+\mathcal{O}(\varepsilon^3),\\
(\widetilde{\mathcal{S}}_\lambda^{pp,\varepsilon})_{m,-1} &=-\pi\varepsilon^2\widetilde g_{\lambda,-}\delta_{m0}+\mathcal{O}(\varepsilon^3),
\end{align*}
and $(\widetilde{\mathcal{S}}_\lambda^{pp,\varepsilon})_{m,n}=\mathcal{O}(\varepsilon^3)$ for $|n|>1$.

Equation~\eqref{eq:self-K-smooth} gives
$$
\widetilde{\mathcal{K}}_\lambda^{pp,\varepsilon,*}[\varphi_0] =\pi\varepsilon (\widetilde g_{\lambda,-}\varphi_1 +\widetilde g_{\lambda,+}\varphi_{-1})+\mathcal{O}(\varepsilon^2),
$$
with an $\mathcal{O}(\varepsilon^2)$ term for $n\ne0$. Equation~\eqref{eq:self-T-smooth} gives
$$
\widetilde{\mathcal{T}}_\lambda^{pp,\varepsilon}[\varphi_0] =\pi\varepsilon (-\ii\widetilde g_{\lambda,-}\varphi_1 +\ii\widetilde g_{\lambda,+}\varphi_{-1})+\mathcal{O}(\varepsilon^2),
$$
For $n\ne0$, the operator action is again $\mathcal{O}(\varepsilon^2)$. For $n=0$,
\begin{align*}
(\widetilde{\mathcal{K}}_\lambda^{pp,\varepsilon,*})_{m,0} &=\pi\varepsilon\bigl(\widetilde g_{\lambda,-}\delta_{m1} +\widetilde g_{\lambda,+}\delta_{m,-1}\bigr)+\mathcal{O}(\varepsilon^2),\\
(\widetilde{\mathcal{T}}_\lambda^{pp,\varepsilon})_{m,0} &=\pi\varepsilon\bigl(-\ii\widetilde g_{\lambda,-}\delta_{m1} +\ii\widetilde g_{\lambda,+}\delta_{m,-1}\bigr)+\mathcal{O}(\varepsilon^2).
\end{align*}

\end{proof}

Let $p\ne q$ and $\bm d_{pq}=\bm z_p-\bm z_q$. Since the kernels are smooth near $\bm d_{pq}\ne\bm0$, their cross-component entries have the following common form.

\begin{lemma}\label{lem:cross-component}
Let $p\ne q$ and $\bm d_{pq}=\bm z_p-\bm z_q$. Write
$$
\widetilde{\bm g}_{\lambda,pq} =\nabla_{\bm x}\widetilde G^\lambda(\bm d_{pq}) =(\widetilde g_{\lambda,pq,1},\widetilde g_{\lambda,pq,2}), \qquad  \widetilde g_{\lambda,pq,+} :=\widetilde g_{\lambda,pq,1} +\ii\widetilde g_{\lambda,pq,2},  \qquad \widetilde g_{\lambda,pq,-} :=\widetilde g_{\lambda,pq,1} -\ii\widetilde g_{\lambda,pq,2}.
$$
The smooth single-layer entries satisfy
$$
(\widetilde{\mathcal{S}}_\lambda^{pq,\varepsilon})_{m,n} =
\begin{cases}
2\pi\varepsilon \widetilde G^\lambda(\bm d_{pq}) +\mathcal{O}(\varepsilon^3), & m=n=0,\\
\pi\varepsilon^2\widetilde g_{\lambda,pq,-} +\mathcal{O}(\varepsilon^3), & m=1, n=0,\\
\pi\varepsilon^2\widetilde g_{\lambda,pq,+} +\mathcal{O}(\varepsilon^3), & m=-1, n=0,\\
-\pi\varepsilon^2\widetilde g_{\lambda,pq,+} +\mathcal{O}(\varepsilon^3), & m=0, n=1,\\
-\pi\varepsilon^2\widetilde g_{\lambda,pq,-} +\mathcal{O}(\varepsilon^3), & m=0, n=-1,\\
\mathcal{O}(\varepsilon^3), & \text{otherwise}.
\end{cases}
$$
The smooth Neumann--Poincar\'e entries satisfy
$$
(\widetilde{\mathcal{K}}_\lambda^{pq,\varepsilon,*})_{m,n} =
\begin{cases}
\pi\varepsilon \widetilde g_{\lambda,pq,-} +\mathcal{O}(\varepsilon^2), & m=1, n=0,\\
\pi\varepsilon \widetilde g_{\lambda,pq,+} +\mathcal{O}(\varepsilon^2), & m=-1, n=0,\\
\mathcal{O}(\varepsilon^2), & \text{otherwise}.
\end{cases}
$$
The smooth tangential entries satisfy
$$
(\widetilde{\mathcal{T}}_\lambda^{pq,\varepsilon})_{m,n} =
\begin{cases}
-\ii\pi\varepsilon \widetilde g_{\lambda,pq,-} +\mathcal{O}(\varepsilon^2), & m=1, n=0,\\
\ii\pi\varepsilon \widetilde g_{\lambda,pq,+} +\mathcal{O}(\varepsilon^2), & m=-1, n=0,\\
\mathcal{O}(\varepsilon^2), & \text{otherwise}.
\end{cases}
$$
\end{lemma}

\begin{proof}
Set $\bm d(\theta,t)=\bm r(\theta)-\bm r(t)$. Since $\bm d_{pq}\notin\Lambda$, the smooth remainder $\widetilde G^\lambda$ is smooth in a fixed neighborhood of $\bm d_{pq}$. Taylor expansion gives
$$
\widetilde G^\lambda (\bm d_{pq}+\varepsilon\bm d(\theta,t)) = \widetilde G^\lambda(\bm d_{pq}) +\varepsilon\widetilde{\bm g}_{\lambda,pq} \cdot(\bm r(\theta)-\bm r(t)) +\mathcal{O}(\varepsilon^2),
$$
and
$$
\nabla_{\bm x}\widetilde G^\lambda (\bm d_{pq}+\varepsilon\bm d(\theta,t)) = \widetilde{\bm g}_{\lambda,pq}+\mathcal{O}(\varepsilon).
$$
The remainders are uniform under the stated nonresonance and parameter assumptions.

Write
$$
\widetilde{\bm g}_{\lambda,pq}\cdot\bm r(\theta) = \frac12\widetilde g_{\lambda,pq,-}e^{\ii\theta} +\frac12\widetilde g_{\lambda,pq,+}e^{-\ii\theta},
$$
and
$$
\bm\tau(\theta)\cdot\widetilde{\bm g}_{\lambda,pq} = -\frac{\ii}{2}\widetilde g_{\lambda,pq,-}e^{\ii\theta} +\frac{\ii}{2}\widetilde g_{\lambda,pq,+}e^{-\ii\theta}.
$$
Using
$$
\int_0^{2\pi}\varphi_n(t) dt =\sqrt{2\pi}\delta_{n0}, \qquad  e^{\pm \ii\theta} =\sqrt{2\pi} \varphi_{\pm1}(\theta),
$$
the single-layer expansion gives, for $n=0,1,-1$,
\begin{align*}
\widetilde{\mathcal{S}}_\lambda^{pq,\varepsilon}[\varphi_0] &=2\pi\varepsilon \widetilde G^\lambda(\bm d_{pq})\varphi_0 +\pi\varepsilon^2 \bigl(\widetilde g_{\lambda,pq,-}\varphi_1 +\widetilde g_{\lambda,pq,+}\varphi_{-1}\bigr) +\mathcal{O}(\varepsilon^3),\\
\widetilde{\mathcal{S}}_\lambda^{pq,\varepsilon}[\varphi_1] &=-\pi\varepsilon^2 \widetilde g_{\lambda,pq,+}\varphi_0 +\mathcal{O}(\varepsilon^3),\\
\widetilde{\mathcal{S}}_\lambda^{pq,\varepsilon}[\varphi_{-1}] &=-\pi\varepsilon^2 \widetilde g_{\lambda,pq,-}\varphi_0 +\mathcal{O}(\varepsilon^3).
\end{align*}
For $|n|>1$, the constant and linear Fourier integrals vanish, so $\widetilde{\mathcal{S}}_\lambda^{pq,\varepsilon}[\varphi_n] =\mathcal{O}(\varepsilon^3)$. Thus, by $\langle\varphi_m,\varphi_k\rangle=\delta_{mk}$,
\begin{align*}
(\widetilde{\mathcal{S}}_\lambda^{pq,\varepsilon})_{m,0} &=2\pi\varepsilon \widetilde G^\lambda(\bm d_{pq})\delta_{m0} +\pi\varepsilon^2 \bigl( \widetilde g_{\lambda,pq,-}\delta_{m1} +\widetilde g_{\lambda,pq,+}\delta_{m,-1} \bigr)+\mathcal{O}(\varepsilon^3),\\
(\widetilde{\mathcal{S}}_\lambda^{pq,\varepsilon})_{m,1} &=-\pi\varepsilon^2 \widetilde g_{\lambda,pq,+}\delta_{m0} +\mathcal{O}(\varepsilon^3),\\
(\widetilde{\mathcal{S}}_\lambda^{pq,\varepsilon})_{m,-1} &=-\pi\varepsilon^2 \widetilde g_{\lambda,pq,-}\delta_{m0} +\mathcal{O}(\varepsilon^3),
\end{align*}
and $(\widetilde{\mathcal{S}}_\lambda^{pq,\varepsilon})_{m,n} =\mathcal{O}(\varepsilon^3)$ for $|n|>1$.

For the normal-derivative operator, the gradient expansion gives
$$
\widetilde{\mathcal{K}}_\lambda^{pq,\varepsilon,*}[\varphi_0] = \pi\varepsilon \bigl( \widetilde g_{\lambda,pq,-}\varphi_1 +\widetilde g_{\lambda,pq,+}\varphi_{-1} \bigr)+\mathcal{O}(\varepsilon^2),
$$
with an $\mathcal{O}(\varepsilon^2)$ term for $n\ne0$. Similarly, the tangential operator satisfies
$$
\widetilde{\mathcal{T}}_\lambda^{pq,\varepsilon}[\varphi_0] = \pi\varepsilon \bigl( -\ii\widetilde g_{\lambda,pq,-}\varphi_1 +\ii\widetilde g_{\lambda,pq,+}\varphi_{-1} \bigr)+\mathcal{O}(\varepsilon^2),
$$
again with an $\mathcal{O}(\varepsilon^2)$ term for $n\ne0$. For $n=0$, orthogonality yields
\begin{align*}
(\widetilde{\mathcal{K}}_\lambda^{pq,\varepsilon,*})_{m,0} &=\pi\varepsilon \bigl( \widetilde g_{\lambda,pq,-}\delta_{m1} +\widetilde g_{\lambda,pq,+}\delta_{m,-1} \bigr)+\mathcal{O}(\varepsilon^2),\\
(\widetilde{\mathcal{T}}_\lambda^{pq,\varepsilon})_{m,0} &=\pi\varepsilon \bigl( -\ii\widetilde g_{\lambda,pq,-}\delta_{m1} +\ii\widetilde g_{\lambda,pq,+}\delta_{m,-1} \bigr)+\mathcal{O}(\varepsilon^2).
\end{align*}
Together with the estimates for $n\ne0$, these give the stated matrix elements.
\end{proof}

\subsection{Fourier matrix elements}

Define $\Delta_{\mathrm m}:=\omega_{\mathrm m}(\omega)^2-|\bm\alpha^*|^2$, $\Delta_{\mathrm c}:=\omega_{\mathrm c}(\omega)^2-|\bm\alpha^*|^2$. We use $\omega_s$ and $\Delta_s$, $s\in\{\mathrm m,\mathrm c\}$, to denote either medium. Recall that $\widehat{\mathcal{S}}_\lambda^{\bm\alpha,\varepsilon} :=\varepsilon^{-1}\mathcal{S}_\lambda^{\bm\alpha,\varepsilon}$.

The three vectors $\bm\xi_\ell$ defined in Lemma \ref{lem:resonant-green-local} have common length $k_D:=|\bm\xi_\ell|=2Q/3=4\pi/(3a_{\rm lat})$. Let $z_D:=k_D\varepsilon$. Choose $\theta_1=-\pi/6$, $\theta_2=-5\pi/6$, $\theta_3=\pi/2$, so that $\bm\xi_\ell=k_D(\cos\theta_\ell,\sin\theta_\ell)^\top$. For $p,q\in\{A,B\}$ and $r\in\mathbb Z$, define $\bm d_{pq}:=\bm z_p-\bm z_q$ and $\Xi_{pq}(r):=\sum_{\ell=1}^{3} e^{\ii\bm\xi_\ell\cdot\bm d_{pq}}e^{\ii r\theta_\ell}$.

\begin{lemma}
Let $J_n$ be the Bessel function of the first kind and $J_n'$ its derivative. The finite resonant part satisfies
\begin{equation}\label{eq:exact-res-S}
(\widehat{\mathcal{S}}_{\Lambda_0^*,\lambda}^{pq,\varepsilon})_{m,n} = \frac{2\pi}{|Y|\Delta_\lambda} (-1)^n\ii^{m+n}J_m(z_D)J_n(z_D)\Xi_{pq}(n-m),
\end{equation}
and
\begin{align}
(\mathcal{K}_{\Lambda_0^*,\lambda}^{pq,\varepsilon,*})_{m,n} &= \frac{2\pi\varepsilon k_D}{|Y|\Delta_\lambda} (-1)^n\ii^{m+n}J_m'(z_D)J_n(z_D)\Xi_{pq}(n-m), \label{eq:exact-res-K}\\
(\mathcal{T}_{\Lambda_0^*,\lambda}^{pq,\varepsilon})_{m,n} &= -\ii m (\widehat{\mathcal{S}}_{\Lambda_0^*,\lambda}^{pq,\varepsilon})_{m,n}. \label{eq:exact-res-T}
\end{align}
A nonzero entry necessarily satisfies $m+s_p\equiv n+s_q\pmod3$.
\end{lemma}

\begin{proof}
For every $\ell$, apply the Jacobi--Anger identity
$$
e^{\ii z_D\cos(\theta-\theta_{\ell})} = \sum\limits_{r\in\mathbb{Z}}\ii^r J_r(z_D)e^{\ii r(\theta-\theta_\ell)}.
$$
By the orthogonality of $\{\varphi_n\}$,
\begin{align*}
\int_0^{2\pi}\overline{\varphi_m(\theta)} e^{\ii\varepsilon\bm\xi_\ell\cdot\bm r(\theta)} d\theta &=\sqrt{2\pi} \ii^mJ_m(z_D)e^{-\ii m\theta_\ell},\\
\int_0^{2\pi}e^{-\ii\varepsilon\bm\xi_\ell\cdot\bm r(t)} \varphi_n(t) dt &=\sqrt{2\pi} (-\ii)^nJ_n(z_D)e^{\ii n\theta_\ell}.
\end{align*}

On the other hand,
$$
G_{\Lambda_0^*}^{\lambda} \bigl( \bm d_{pq} +\varepsilon(\bm r(\theta)-\bm r(t)) \bigr) = \frac1{|Y|\Delta_\lambda} \sum_{\ell=1}^{3} e^{\ii\bm\xi_\ell\cdot\bm d_{pq}} e^{\ii\varepsilon\bm\xi_\ell\cdot\bm r(\theta)} e^{-\ii\varepsilon\bm\xi_\ell\cdot\bm r(t)}.
$$
Then
\begin{align*}
(\widehat{\mathcal{S}}_{\Lambda_0^*,\lambda}^{pq,\varepsilon})_{m,n} &= \frac1{|Y|\Delta_\lambda} \sum_{\ell=1}^{3} e^{\ii\bm\xi_\ell\cdot\bm d_{pq}} \left( \sqrt{2\pi} \ii^mJ_m(z_D)e^{-\ii m\theta_\ell} \right) \left( \sqrt{2\pi} (-\ii)^nJ_n(z_D)e^{\ii n\theta_\ell} \right)\\
&= \frac{2\pi}{|Y|\Delta_\lambda} \ii^m(-\ii)^nJ_m(z_D)J_n(z_D) \sum_{\ell=1}^{3} e^{\ii\bm\xi_\ell\cdot\bm d_{pq}} e^{\ii(n-m)\theta_\ell}\\
&= \frac{2\pi}{|Y|\Delta_\lambda} \ii^m(-\ii)^nJ_m(z_D)J_n(z_D)\Xi_{pq}(n-m).
\end{align*}
That is~\eqref{eq:exact-res-S}.

For the normal derivative, we have $\ii\bm\xi_\ell\cdot\bm\nu(\theta) e^{\ii\varepsilon\bm\xi_\ell\cdot\bm r(\theta)} =\partial_\varepsilon e^{\ii\varepsilon\bm\xi_\ell\cdot\bm r(\theta)}$, hence
$$
\int_0^{2\pi} \overline{\varphi_m(\theta)}  \ii\bm\xi_\ell\cdot\bm\nu(\theta) e^{\ii\varepsilon\bm\xi_\ell\cdot\bm r(\theta)}  d\theta = \partial_\varepsilon \left( \sqrt{2\pi} \ii^mJ_m(k_D\varepsilon)e^{-\ii m\theta_\ell} \right) = \sqrt{2\pi} \ii^mk_DJ_m'(z_D)e^{-\ii m\theta_\ell},
$$
and
\begin{align*}
(\mathcal{K}_{\Lambda_0^*,\lambda}^{pq,\varepsilon,*})_{m,n} &= \frac{\varepsilon}{|Y|\Delta_\lambda} \sum_{\ell=1}^{3} e^{\ii\bm\xi_\ell\cdot\bm d_{pq}} \left( \sqrt{2\pi} \ii^mk_DJ_m'(z_D)e^{-\ii m\theta_\ell} \right) \left( \sqrt{2\pi} (-\ii)^nJ_n(z_D)e^{\ii n\theta_\ell} \right)\\
&= \frac{2\pi\varepsilon k_D}{|Y|\Delta_\lambda} \ii^m(-\ii)^nJ_m'(z_D)J_n(z_D)\Xi_{pq}(n-m),
\end{align*}
which is \eqref{eq:exact-res-K}.

For the tangential derivative, we have $\ii\bm\xi_\ell\cdot\bm\tau(\theta)e^{\ii\varepsilon\bm\xi_\ell\cdot\bm r(\theta)} =-\varepsilon^{-1}\partial_\theta e^{\ii\varepsilon\bm\xi_\ell\cdot\bm r(\theta)}$. Using integration by parts we get
\begin{align*}
&\int_0^{2\pi} \overline{\varphi_m(\theta)}  \ii\bm\xi_\ell\cdot\bm\tau(\theta) e^{\ii\varepsilon\bm\xi_\ell\cdot\bm r(\theta)}  d\theta = -\frac1{\varepsilon} \int_0^{2\pi} \overline{\varphi_m(\theta)} \partial_\theta e^{\ii\varepsilon\bm\xi_\ell\cdot\bm r(\theta)}  d\theta\\
&= \frac1{\varepsilon} \int_0^{2\pi} \partial_\theta\overline{\varphi_m(\theta)} e^{\ii\varepsilon\bm\xi_\ell\cdot\bm r(\theta)}  d\theta = -\frac{\ii m}{\varepsilon} \sqrt{2\pi} \ii^mJ_m(z_D)e^{-\ii m\theta_\ell}.
\end{align*}
Hence
$$
(\mathcal{T}_{\Lambda_0^*,\lambda}^{pq,\varepsilon})_{m,n} = -\frac{2\pi \ii m}{|Y|\Delta_\lambda} \ii^m(-\ii)^nJ_m(z_D)J_n(z_D)\Xi_{pq}(n-m) = -\ii m  (\widehat{\mathcal{S}}_{\Lambda_0^*,\lambda}^{pq,\varepsilon})_{m,n},
$$
which is~\eqref{eq:exact-res-T}.

From $\theta_2=\theta_1-\frac{2\pi}{3}$ and $\theta_3=\theta_1+\frac{2\pi}{3}$ we have $e^{\ii r\theta_\ell}=e^{\ii r\theta_1}e^{-2\pi \ii(\ell-1)r/3}, \ell=1,2,3$. Using $\bm d_{AA}=\bm d_{BB}=\bm0$, $\bm d_{AB}=-(\bm e_1+\bm e_2)/3$, $\bm d_{BA}=-\bm d_{AB}$ and $s_A=0, s_B=1$, direct calculation gives $e^{\ii\bm\xi_\ell\cdot\bm d_{pq}}=e^{\ii\bm\xi_1\cdot\bm d_{pq}}e^{2\pi \ii(\ell-1)(s_p-s_q)/3}$. As a result,
$$
\Xi_{pq}(r) = e^{\ii\bm\xi_1\cdot\bm d_{pq}}e^{\ii r\theta_1} \sum_{\ell=1}^{3} e^{2\pi \ii(\ell-1)(s_p-s_q-r)/3} = e^{\ii\bm\xi_1\cdot\bm d_{pq}}e^{\ii r\theta_1} \left( 1+e^{2\pi \ii(s_p-s_q-r)/3} +e^{4\pi \ii(s_p-s_q-r)/3} \right),
$$
and $\Xi_{pq}(r)\neq0$ is equivalent to $r\equiv s_p-s_q\pmod3$. $m+s_p\equiv n+s_q\pmod3$ follows by letting $r=n-m$.
\end{proof}

We next compute the low-order self-interaction entries needed for the candidate frequencies. Let $(\widehat{\mathcal{S}}_{\omega_s}^{pp,\varepsilon})_{n,n}$, $(\mathcal{K}_{\omega_s}^{pp,\varepsilon,*})_{n,n}$, and $(\mathcal{T}_{\omega_s}^{pp,\varepsilon})_{n,n}$ denote the self-interaction entries of $\widehat{\mathcal{S}}_{\omega_s}^{\bm\alpha^*,\varepsilon}$, $(\mathcal{K}_{\omega_s}^{-\bm\alpha^*,\varepsilon})^*$, and $\mathcal{T}_{\omega_s}^{\bm\alpha^*,\varepsilon}$, respectively, on the $n$th Fourier mode.

\begin{lemma}\label{lem:low-self-fourier}
Assume $|\Delta_s|\ge c\varepsilon^2$, with $c>0$ independent of $\varepsilon$. Then
\begin{align*}
(\widehat{\mathcal{S}}_{\omega_s}^{pp,\varepsilon})_{0,0} &= \ln\varepsilon+\ln\omega_s+\gamma_0 +\frac{6\pi}{|Y|\Delta_s} -\frac{4\pi Q^2\varepsilon^2}{3|Y|\Delta_s} +2\pi\widetilde G^{\omega_s}(\bm0) +\mathcal{O}(\varepsilon^2|\ln\varepsilon|),\\
(\widehat{\mathcal{S}}_{\omega_s}^{pp,\varepsilon})_{\pm1,\pm1} &= -\frac12+ \frac{2\pi Q^2\varepsilon^2}{3|Y|\Delta_s} +\mathcal{O}(\varepsilon^2|\ln\varepsilon|),\\
(\mathcal{K}_{\omega_s}^{pp,\varepsilon,*})_{0,0} &= \frac12- \frac{4\pi Q^2\varepsilon^2}{3|Y|\Delta_s} +\mathcal{O}(\varepsilon^2|\ln\varepsilon|),\\
(\mathcal{K}_{\omega_s}^{pp,\varepsilon,*})_{\pm1,\pm1} &= \frac{2\pi Q^2\varepsilon^2}{3|Y|\Delta_s} +\mathcal{O}(\varepsilon^2|\ln\varepsilon|),\\
(\mathcal{T}_{\omega_s}^{pp,\varepsilon})_{0,0} &=0,\\
(\mathcal{T}_{\omega_s}^{pp,\varepsilon})_{1,1} &= \frac{\ii}{2} -\frac{2\pi \ii Q^2\varepsilon^2}{3|Y|\Delta_s} +\mathcal{O}(\varepsilon^2|\ln\varepsilon|),\\
(\mathcal{T}_{\omega_s}^{pp,\varepsilon})_{-1,-1} &= -\frac{\ii}{2} +\frac{2\pi \ii Q^2\varepsilon^2}{3|Y|\Delta_s} +\mathcal{O}(\varepsilon^2|\ln\varepsilon|).
\end{align*}
\end{lemma}

\begin{proof}
For the single-layer entry with $n=1$, the free-space and smooth contributions from Lemmas~\ref{lem:appendix-a-1} and~\ref{lem:self-interaction-smooth} are, respectively, $(\widehat{\mathcal{S}}_{H_0,\omega_s}^{pp,\varepsilon})_{1,1} =-1/2+\mathcal{O}(\varepsilon^2|\ln\varepsilon|)$ and $(\widetilde{\widehat{\mathcal{S}}}_{\omega_s}^{pp,\varepsilon})_{1,1} =\mathcal{O}(\varepsilon^2)$. Also, $\Xi_{pp}(0)=3$, $J_1(z_D)=z_D/2+\mathcal{O}(z_D^3)$, and $z_D=2Q\varepsilon/3$ give
\begin{align*}
(\widehat{\mathcal{S}}_{\omega_s}^{pp,\varepsilon})_{1,1} &= -\frac12+ \frac{6\pi}{|Y|\Delta_s}J_1(z_D)^2 +\mathcal{O}(\varepsilon^2|\ln\varepsilon|)\\
&= -\frac12+ \frac{6\pi}{|Y|\Delta_s} \left( \frac{Q^2\varepsilon^2}{9} +\mathcal{O}(\varepsilon^4) \right) +\mathcal{O}(\varepsilon^2|\ln\varepsilon|)\\
&= -\frac12+ \frac{2\pi Q^2\varepsilon^2}{3|Y|\Delta_s} +\mathcal{O}(\varepsilon^2|\ln\varepsilon|).
\end{align*}
Here, $\mathcal{O}(\varepsilon^4/|\Delta_s|)=\mathcal{O}(\varepsilon^2)$. Since $J_{-1}=-J_1$, the $n=-1$ single-layer expansion is identical.

For $n=0$, use $J_0(z_D)=1-z_D^2/4+\mathcal{O}(z_D^4)$,
\begin{align*}
(\widehat{\mathcal{S}}_{\omega_s}^{pp,\varepsilon})_{0,0} &= \ln\varepsilon+\ln\omega_s+\gamma_0 +2\pi\widetilde G^{\omega_s}(\bm0) +\frac{6\pi}{|Y|\Delta_s}J_0(z_D)^2 +\mathcal{O}(\varepsilon^2|\ln\varepsilon|)\\
&= \ln\varepsilon+\ln\omega_s+\gamma_0 +\frac{6\pi}{|Y|\Delta_s} -\frac{4\pi Q^2\varepsilon^2}{3|Y|\Delta_s} +2\pi\widetilde G^{\omega_s}(\bm0) +\mathcal{O}(\varepsilon^2|\ln\varepsilon|).
\end{align*}
Then use $J_0'=-J_1$, $J_1'(z_D)=1/2+\mathcal{O}(\varepsilon^2)$, and the free-space normal entries to obtain
\begin{align*}
(\mathcal{K}_{\omega_s}^{pp,\varepsilon,*})_{0,0} &= \frac12+ \frac{6\pi\varepsilon k_D}{|Y|\Delta_s} J_0'(z_D)J_0(z_D) +\mathcal{O}(\varepsilon^2|\ln\varepsilon|) = \frac12- \frac{4\pi Q^2\varepsilon^2}{3|Y|\Delta_s} +\mathcal{O}(\varepsilon^2|\ln\varepsilon|),\\
(\mathcal{K}_{\omega_s}^{pp,\varepsilon,*})_{1,1} &= \frac{6\pi\varepsilon k_D}{|Y|\Delta_s} J_1'(z_D)J_1(z_D) +\mathcal{O}(\varepsilon^2|\ln\varepsilon|) = \frac{2\pi Q^2\varepsilon^2}{3|Y|\Delta_s} +\mathcal{O}(\varepsilon^2|\ln\varepsilon|).
\end{align*}
The normal entry is the same for $n=-1$. For the full single-layer operator, $\partial_{\bm\tau}=-\varepsilon^{-1}\partial_\theta$ gives $(\mathcal{T}_{\omega_s}^{pp,\varepsilon})_{m,n}=-\ii m(\widehat{\mathcal{S}}_{\omega_s}^{pp,\varepsilon})_{m,n}$. This identity, together with~\eqref{eq:exact-res-T}, gives the expansions of $(\mathcal{T}_{\omega_s}^{pp,\varepsilon})_{0,0}$ and $(\mathcal{T}_{\omega_s}^{pp,\varepsilon})_{\pm1,\pm1}$.
\end{proof}

\begin{lemma} \label{lem:low-resonant-phases}
$ \Xi_{AB}(-1)=-3\ii, \Xi_{BA}(1)=3\ii, \Xi_{AB}(2)=\Xi_{BA}(-2)=-3, $ $\Xi_{AA}(3)=\Xi_{BB}(3)=-3\ii$, $\Xi_{AA}(-3)=\Xi_{BB}(-3)=3\ii$. Moreover, $\Xi_{AB}(0)=\Xi_{BA}(0)=0$.
\end{lemma}

\begin{proof}
Since $\bm d_{AB}=-(\bm e_1+\bm e_2)/3$, direct calculation gives
\begin{equation}\label{eq:AB-resonant-phases}
e^{\ii\bm\xi_1\cdot\bm d_{AB}}=e^{-2\pi \ii/3}, \qquad  e^{\ii\bm\xi_2\cdot\bm d_{AB}}=e^{2\pi \ii/3}, \qquad  e^{\ii\bm\xi_3\cdot\bm d_{AB}}=1.
\end{equation}
Substituting~\eqref{eq:AB-resonant-phases} and $\theta_1=-\pi/6$, $\theta_2=-5\pi/6$, $\theta_3=\pi/2$ into the definition of $\Xi_{pq}(r)$ gives $\Xi_{AB}(-1)$, $\Xi_{AB}(2)$, and $\Xi_{AB}(0)$. The relation $\bm d_{BA}=-\bm d_{AB}$ gives the corresponding $BA$ identities. For $p=q$, $\bm d_{pp}=\bm0$, so $\Xi_{pp}(r)=\sum_{\ell=1}^{3}e^{\ii r\theta_\ell}$. Taking $r=\pm3$ proves the remaining identities.
\end{proof}

\end{appendices}



\begin{thebibliography}{99}

\bibitem{AmmariFitzpatrickKangRuizYuZhang2018} H.~Ammari, B.~Fitzpatrick, H.~Kang, M.~Ruiz, S.~Yu, and H.~Zhang, \emph{Mathematical and Computational Methods in Photonics and Phononics}, Mathematical Surveys and Monographs, vol.~235, American Mathematical Society, Providence, RI, 2018.

\bibitem{AmmariKangLee2009} H.~Ammari, H.~Kang, and H.~Lee, \emph{Layer Potential Techniques in Spectral Analysis}, Mathematical Surveys and Monographs, vol.~153, American Mathematical Society, Providence, RI, 2009.

\bibitem{AmmariKosche2024} H.~Ammari and T.~Kosche, Topological phenomena in honeycomb Floquet metamaterials, \emph{Math. Ann.} \textbf{388} (2024), 2755--2785. \href{https://doi.org/10.1007/s00208-023-02583-0}{doi:10.1007/s00208-023-02583-0}.

\bibitem{ColtonKress2013} D.~Colton and R.~Kress, \emph{Integral Equation Methods in Scattering Theory}, Classics in Applied Mathematics, vol.~72, Society for Industrial and Applied Mathematics, Philadelphia, PA, 2013.

\bibitem{FeffermanLeeThorpWeinstein2016} C.~L.~Fefferman, J.~P.~Lee-Thorp, and M.~I.~Weinstein, Edge states in honeycomb structures, \emph{Ann. PDE} \textbf{2} (2016), Article~12. \href{https://doi.org/10.1007/s40818-016-0015-3}{doi:10.1007/s40818-016-0015-3}.

\bibitem{FeffermanWeinstein2012} C.~L.~Fefferman and M.~I.~Weinstein, Honeycomb lattice potentials and Dirac points, \emph{J. Amer. Math. Soc.} \textbf{25} (2012), 1169--1220. \href{https://doi.org/10.1090/S0894-0347-2012-00745-0}{doi:10.1090/S0894-0347-2012-00745-0}.

\bibitem{FeffermanWeinstein2014} C.~L.~Fefferman and M.~I.~Weinstein, Wave packets in honeycomb structures and two-dimensional Dirac equations, \emph{Comm. Math. Phys.} \textbf{326} (2014), 251--286. \href{https://doi.org/10.1007/s00220-013-1847-2}{doi:10.1007/s00220-013-1847-2}.

\bibitem{John1987} S.~John, Strong localization of photons in certain disordered dielectric superlattices, \emph{Phys. Rev. Lett.} \textbf{58} (1987), 2486--2489. \href{https://doi.org/10.1103/PhysRevLett.58.2486}{doi:10.1103/PhysRevLett.58.2486}.

\bibitem{Kuchment2016} P.~Kuchment, An overview of periodic elliptic operators, \emph{Bull. Amer. Math. Soc. (N.S.)} \textbf{53} (2016), 343--414.

\bibitem{LiLinQiuZhang2025} W.~Li, J.~Lin, J.~Qiu, and H.~Zhang, Mathematical theory for photonic Hall effect in honeycomb photonic crystals, preprint, arXiv:2512.24477, 2025. \url{https://arxiv.org/abs/2512.24477}.

\bibitem{Moroz2006} A.~Moroz, Quasi-periodic Green's functions of the Helmholtz and Laplace equations, \emph{J. Phys. A: Math. Gen.} \textbf{39} (2006), 11247--11282.

\bibitem{Olver2010} F.~W.~J.~Olver, D.~W.~Lozier, R.~F.~Boisvert, and C.~W.~Clark, eds., \emph{NIST Handbook of Mathematical Functions}, Cambridge University Press, Cambridge, 2010.

\bibitem{QiuLinXieZhang2023} J.~Qiu, J.~Lin, P.~Xie, and H.~Zhang, Mathematical theory for the interface mode in an acoustic waveguide bifurcated from a Dirac point, \emph{Commun. Math. Sci.} \textbf{24} (2026), 749--796.

\bibitem{RaghuHaldane2008} S.~Raghu and F.~D.~M.~Haldane, Analogs of quantum-Hall-effect edge states in photonic crystals, \emph{Phys. Rev. A} \textbf{78} (2008), 033834. \href{https://doi.org/10.1103/PhysRevA.78.033834}{doi:10.1103/PhysRevA.78.033834}.

\bibitem{Wang2008} Z.~Wang, Y.~D.~Chong, J.~D.~Joannopoulos, and M.~Solja\v{c}i\'c, Reflection-free one-way edge modes in a gyromagnetic photonic crystal, \emph{Phys. Rev. Lett.} \textbf{100} (2008), 013905. \href{https://doi.org/10.1103/PhysRevLett.100.013905}{doi:10.1103/PhysRevLett.100.013905}.

\bibitem{Wang2009} Z.~Wang, Y.~Chong, J.~D.~Joannopoulos, and M.~Solja\v{c}i\'c, Observation of unidirectional backscattering-immune topological electromagnetic states, \emph{Nature} \textbf{461} (2009), 772--775. \href{https://doi.org/10.1038/nature08293}{doi:10.1038/nature08293}.

\bibitem{Yablonovitch1987} E.~Yablonovitch, Inhibited spontaneous emission in solid-state physics and electronics, \emph{Phys. Rev. Lett.} \textbf{58} (1987), 2059--2062. \href{https://doi.org/10.1103/PhysRevLett.58.2059}{doi:10.1103/PhysRevLett.58.2059}.
    
\bibitem{LiLinZhang2023} W. Li, J. Lin, and H. Zhang, Dirac points for the honeycomb lattice with impenetrable obstacles, \emph{SIAM J. Appl. Math.}, \textbf{83} (2023), no.~4, 1546--1571.

\end{thebibliography}
\end{document}